\documentclass[12pt]{article}
\usepackage{amsmath,amssymb,amsthm,longtable}
\numberwithin{equation}{section}

\newcommand{\zero}{\mathbf {0}}
\newcommand{\1}{\mathbf {1}}
\newcommand{\N}{{\mathbb N}}
\newcommand{\Z}{{\mathbb Z}}
\newcommand{\Q}{{\mathbb Q}}
\newcommand{\R}{{\mathbb R}}
\newcommand{\C}{{\mathbb C}}

\newcommand{\mbbM}{{\mathbb M}}

\newcommand{\K}{{\mathcal K}}
\newcommand{\mcL}{\widehat{(L^{\perp})^{\oplus {\ell}}}}
\newcommand{\M}{{\mathcal M}}

\newcommand{\Klein}{{\mathcal K}}
\newcommand{\Lat}{\Gamma}
\newcommand{\h}{{\mathfrak h}}
\newcommand{\grs}{\bar{H}_{\ell}}
\newcommand{\tgrs}{{H}_{\ell}}
\newcommand{\tgrl}{{G}_{\ell}}

\newcommand{\kz}{\kappa_{36}}
\newcommand{\lweight}{\Delta}
\newcommand{\Leech}{\Lambda}

\newcommand{\LCzero}{L_{C \times \zero }}
\newcommand{\LzeroD}{L_{\zero \times D}}
\newcommand{\LCD}{L_{C \times D}}

\newcommand{\hLCD}{\hat{L}_{C \times D}}
\newcommand{\hLCDtau}{\hat{L}_{C \times D,\, \htau}}
\newcommand{\hLCzerotau}{\hat{L}_{C \times \zero,\, \htau}}

\newcommand{\Tpe}{T_{\psi_{\eta}}}
\newcommand{\Ceta}{\C_{\psi_{\eta}}}

\newcommand{\hMtau}{\hat{M}_{\htau}}
\newcommand{\hRtau}{\hat{R}_{\htau}}

\newcommand{\htheta}{\theta}
\newcommand{\htau}{{\tau}}
\newcommand{\otau}{\tau}
\newcommand{\ova}{\overline{a}}

\newcommand{\VLCD}{V_{L_{C \times D}}}

\newcommand{\VLCDtau}{V_{L_{C \times D}}^{\htau}}

\newcommand{\VLCDTeta}{V_{L_{C \times D}}^{T, \eta}}

\newcommand{\al}{\alpha}
\newcommand{\be}{\beta}
\newcommand{\dl}{\delta}
\newcommand{\gm}{\gamma}
\newcommand{\lm}{\lambda}
\newcommand{\om}{\omega}

\newcommand{\vep}{\varepsilon}

\newcommand{\la}{\langle}
\newcommand{\op}{\oplus}
\newcommand{\ot}{\otimes}
\newcommand{\ra}{\rangle}

\newcommand{\tll}{m}
\newcommand{\os}{S^{*}}
\newcommand{\anx}{y}
\newcommand{\vt}[4]{V_{#1}^{T,#2}(#3)[#4]}
\newcommand{\e}{\mbox{\bf e}}

\DeclareMathOperator{\Aut}{Aut} \DeclareMathOperator{\End}{End}
\DeclareMathOperator{\Hom}{Hom} 
 
 \DeclareMathOperator{\spn}{span}
\DeclareMathOperator{\supp}{supp} 
\DeclareMathOperator{\wt}{wt}
\DeclareMathOperator{\pr}{pr}
\DeclareMathOperator{\module}{mod}

\newtheorem{lem}{Lemma}[section]
\newtheorem{thm}[lem]{Theorem}
\newtheorem{prop}[lem]{Proposition}

\theoremstyle{definition}
\newtheorem{rmk}[lem]{Remark}

\title{Fixed point subalgebras of lattice vertex operator algebras
by an automorphism of order three}
\author{Kenichiro Tanabe\footnote{Corresponding author.
Partially supported by JSPS Grant-in-Aid for
Scientific Research No. 20740002}\\
Department of Mathematics, Hokkaido University,\\
Sapporo, Hokkaido, 060-0810, Japan\\
ktanabe@math.sci.hokudai.ac.jp\\\\
Hiromichi Yamada\footnote{Partially supported by JSPS Grant-in-Aid for
Scientific Research No. 23540009}\\
Department of Mathematics, Hitotsubashi University,\\
Kunitachi, Tokyo 186-8601, Japan\\
yamada@econ.hit-u.ac.jp}
\date{}
\begin{document}
\maketitle

\begin{abstract}
We study the fixed point subalgebra of a certain class of lattice
vertex operator algebras by an automorphism of order $3$, which is a
lift of a fixed-point-free isometry of the underlying lattice. We
classify the irreducible modules for the subalgebra. Moreover, the
rationality and the $C_2$-cofiniteness of the subalgebra are
established. Our result contains the case of the vertex operator
algebra associated with the Leech lattice.
\end{abstract}



\tableofcontents

\section{Introduction}
Let $V$ be a vertex operator algebra. For an automorphism $g$ of
$V$ of finite order, the space $V^g = \{ v \in V\,|\, gv = v\}$ of
fixed points is a subalgebra of $V$ called an orbifold of 
the vertex operator algebra $V$. 
It is conjectured in \cite{DVVV} that
every irreducible $V^g$-module is contained in some irreducible
untwisted or twisted $V$-module. It is also conjectured that 
if $V$ is rational and $C_2$-cofinite, then so 
is $V^g$.
These conjectures have important meanings in the
theory of vertex operator algebras. However, it is difficult to
investigate an orbifold in general, even if the original vertex
operator algebra $V$ is well understood.

In the case where $V$ is the lattice vertex operator algebra
$V_\Gamma$ associated with a positive definite even lattice
$\Gamma$ and the automorphism $g$ is a canonical lift $\theta$ of
the $-1$ isometry $\al \mapsto -\al$ of the lattice $\Gamma$, the
orbifold $V_\Gamma^{\theta} = V_\Gamma^+$ has been studied
extensively. In fact, the representation theory of $V_\Gamma^+$,
that is, the classification of irreducible modules \cite{AD, DN}
and the determination of fusion rules \cite{A, ADL}, together with 
the $C_2$-cofiniteness \cite{ABD, Yamskulna} 
of $V_\Gamma^+$ are established.

In this paper we study an orbifold of a certain class of lattice
vertex operator algebras by an automorphism of order $3$. We start
with a lattice $L \cong \sqrt{2}$ ($A_2$-lattice) and a
fixed-point-free isometry $\otau$ of $L$ of order $3$. There are
12 cosets of $L$ in its dual lattice $L^\perp$. 
Using an even $\Z_2 \times
\Z_2$-code $C$ of length $\ell$ and a self-orthogonal $\Z_3$-code
$D$ of the same length, we construct a positive definite even
lattice $\LCD\subset (L^{\perp})^{\oplus\ell}$ of rank $2\ell$ 
from the 12 cosets of $L$ in
$L^\perp$. We also consider an action of $\otau$ on $\Z_2 \times \Z_2$.
The isometry $\otau$ induces a fixed-point-free
isometry $(\otau,\ldots,\otau)$ of $\LCD$ provided that $C$ is 
invariant under the corresponding action of $(\otau,\ldots,\otau)$ on 
$(\Z_2 \times \Z_2)^\ell$.
For simplicity of notation, we denote $(\otau,\ldots,\otau)$ by $\htau$ also.

Our main concern is to classify the irreducible
modules for the orbifold $\VLCDtau$ of the lattice vertex operator
algebra $\VLCD$ by an automorphism $\htau$ of order $3$ which is a
lift of the isometry $\htau$ of $\LCD$. 
The vertex operator algebra $\VLCD$ is simple, rational,
$C_2$-cofinite, and of CFT type. 
The dual
lattice $(\LCD)^\perp$ of $\LCD$ is equal to $L_{C^\perp \times
D^\perp}$, where $C^\perp$ (resp. $D^\perp$) is the dual code of
$C$ (resp. $D$). Then $V_{L_{(\lambda + C) \times (\gamma + D)}}$,
$\lambda + C \in C^\perp/C$, $\gamma + D \in D^\perp/D$ form a
complete set of representatives of equivalence classes of
irreducible $\VLCD$-modules.
Such a 
$V_{L_{(\lambda + C) \times (\gamma + D)}}$ is $\htau$-stable if
and only if $\lambda \in C$. 
One can also construct irreducible $\htau^i$-twisted
$\VLCD$-modules $\VLCDTeta(\htau^i)$, $\eta \in D^\perp \pmod{D}$
for $i=1,2$ by the method of \cite{DL2, L}.

The orbifold $\VLCDtau$ is a simple vertex operator algebra. The
following is a list of known irreducible $\VLCDtau$-modules. Let
$\zeta_3 = \exp(2\pi\sqrt{-1}/3)$.

\medskip
(1) $V_{L_{C \times (\gamma + D)}}(\vep) = \{ u \in V_{L_{C \times
(\gamma + D)}}\,|\, \htau u = \zeta_3^\vep u\}$, $\gamma + D \in
D^\perp/D$, $\vep \in \Z_3$.

\medskip
(2) $V_{L_{(\lambda + C) \times (\gamma + D)}}$, $0 \ne \lambda +
C \in (C^\perp/C)_{\equiv_\htau}$, $\gamma + D \in D^\perp/D$,
where $(C^\perp/C)_{\equiv_\htau}$ is the set of $\htau$-orbits in
$C^\perp/C$.

\medskip
(3) $\VLCDTeta(\htau^i)[\vep] = \{ u \in \VLCDTeta(\htau^i)\,|\,
\htau^i u = \zeta_3^\vep u \}$, $\eta \in D^\perp  \pmod{D}$, $\vep
\in \Z_3$, $i=1,2$.

\medskip
These irreducible $\VLCDtau$-modules are inequivalent each other
\cite{DY, MT}. The above mentioned conjecture says that any
irreducible $\VLCDtau$-module is isomorphic to one of these.

In our argument we deal with not only simple current extension 
\cite{DLM0} but
also certain nonsimple current extension. Simple current extension
is rather easy, whereas nonsimple current extension is complicated
and difficult to study. In order to avoid the difficulty, we
restrict ourselves to the special case where $C$ is a
$\htau$-invariant self-dual $\Z_2 \times \Z_2$-code with minimum 
weight at least $4$ and $D$ is a self-dual $\Z_3$-code. In
this case the lattice $\LCD$ is unimodular and there is a unique
irreducible $\VLCD$-module, namely, $\VLCD$ itself. Likewise,
there is a unique irreducible $\htau^i$-twisted $\VLCD$-module
$V_{L_{C\times D}}^{T,\zero}(\htau^i)$, $i=1,2$, where $\zero$ is
the zero codeword. Under this hypothesis we have the following
theorem (Theorem \ref{thm:main}).

\medskip\noindent
{\bfseries Theorem.} \emph{Suppose $C$ is a $\htau$-invariant
self-dual $\Z_2 \times \Z_2$-code with minimum weight at least
$4$ and $D$ is a self-dual $\Z_3$-code. Then the vertex operator
algebra $\VLCDtau$ is simple, rational, $C_2$-cofinite, and of CFT
type. Moreover, every irreducible $V_{L_{C\times
D}}^{\htau}$-module is isomorphic to one of $\VLCD(\vep)$,
$\vt{L_{C\times D}}{\zero}{\htau^i}{\vep}$, $\vep \in \Z_3$,
$i=1,2$.}

\medskip
One of the most important examples of orbifold is the fixed point
subalgebra $V_{\Leech}^\theta$ of the Leech lattice vertex operator
algebra $V_{\Leech}$ by the automorphism $\theta$ of order $2$.
This orbifold was first studied by I. Frenkel, J. Lepowsky and A.
Meurman, and in fact it was used for the construction of the
moonshine vertex operator algebra $V^\natural$ \cite{FLM}.
We note that the Leech lattice $\Leech$ can be expressed as
$\LCD$ for some $C$ and $D$ which satisfy the hypothesis of the
theorem (Remark \ref{rmk:Leech}).

A remarkable property of $V^\natural$ is that its automorphism
group $\Aut V^\natural$ is isomorphic to the Monster $\mbbM$. 
The construction of $V^\natural$ in \cite{FLM} is based on 
a $2B$-element of $\mbbM$. In
\cite[Introduction]{FLM}, it is suggested that an analogous
construction may be possible for some appropriate elements in
$\mbbM$ of order $3$, $5$, $7$, and $13$.
The classification of
irreducible modules, the determination of fusion rules, the
rationality and the $C_2$-cofiniteness 
for the orbifold $V_\Lambda^g$ by such an element $g$
should play an important role in those expected
construction. This is the motivation for the present work.

The organization of the paper is as follows. Section 2 is devoted
to the preliminaries. In Section 2.1 we collect basic terminology
for later use. In Section 2.2 we introduce the lattice $\LCD$ and
study its properties. 
In Section \ref{subsetion:group-twist}
we introduce a central extension $\hLCD$ of $\LCD$
by a group $\la \kz\ra$ of order $36$ and discuss an
action of a lift of the isometry $\htau$ of the lattice $\LCD$. 
In Section \ref{subsec:VOA_VLCD} we study the vertex operator algebra 
$V_{\LCD}$ and its irreducible modules. 
The automorphism $\htau$ of $\hLCD$
naturally induces an automorphism of $V_{\LCD}$ of order $3$, which
is again denoted by $\htau$.

In Section 3 we discuss in detail the irreducible $\htau^i$-twisted
$V_{\LCD}$-modules $\VLCDTeta(\htau^i)$, $i=1,2$, which are
obtained by the method of \cite{DL2,L}. We describe those
irreducible $\htau^i$-twisted $V_{\LCD}$-modules as modules for
$(V_L^{\otau})^{\otimes \ell}$ (Theorem \ref{thm:twist-L-st}). 
The classification of irreducible
modules for the orbifold $V_L^{\otau}$ was accomplished in
\cite{TY}. Our argument here is based on the result. 

In Section 4 we determine certain fusion rules for $V_L^{\otau}$ 
(Proposition \ref{prop:fusionL}),
which will be necessary in Section 5. In fact, these fusion rules
are crucial for our arguments.

The proof of the main theorem is divided into three steps. In
Section 5 we begin with the classification of irreducible modules
for the orbifold $V_{L^{\oplus \ell}}^\htau$ 
(Proposition \ref{prop:irr-L-l}). This is the case where both of
$C$ and $D$ are the zero code. The rationality and the
$C_2$-cofiniteness of $V_{L^{\oplus \ell}}^\htau$ are also
obtained. Moreover, some of the fusion rules are computed 
(Proposition \ref{prop:fusion-L-l}).

In Section 6 we classify the irreducible modules for
$V_{L_{\zero \times D}}^\htau$ (Theorem \ref{thm:irr-D}), 
which is the case where $C$ is the
zero code $\{\zero\}$. In this case only simple current extension
is involved and the argument is relatively straightforward. 
The rationality and
the $C_2$-cofiniteness of $V_{L_{\zero  \times D}}^\htau$ 
(Theorem \ref{thm:irr-D}), together
with some of the fusion rules are also obtained 
(Proposition \ref{prop:fusion-D}).

Section 7 consists of two subsections. In Section 7.1 we use Zhu's
theory to study the irreducible $V_{L^{\oplus\ell}}^\htau$-modules
contained in a $V_{L_{C(\mu) \times \zero}}^\htau$-module, where
$C(\mu)$ is the $\Z_2 \times \Z_2$-code generated by $\mu$ and
$\htau(\mu)$. 
The results obtained here will be necessary
in Section 7.2. We do not discuss the classification of
irreducible modules nor the rationality for the vertex operator
algebra $V_{L_{C(\mu) \times \zero}}^\htau$. Note that
$V_{L_{C(\mu) \times \zero}}^\htau$ is a nonsimple current
extension of $V_{L^{\oplus\ell}}^\htau$.

In Section 7.2 we study the orbifold $V_{\LCD}^\htau$ and prove the
main theorem  (Theorem \ref{thm:main}) under the hypothesis that 
$C$ is a $\htau$-invariant
self-dual $\Z_2 \times \Z_2$-code with minimum weight at least
$4$ and $D$ is a self-dual $\Z_3$-code. We need to assume that
$D$ is self-dual for the proof of 
Proposition \ref{prop:CD-untwist}.
Our argument 
fails if the minimum weight of $C$ is $2$ (Remark \ref{rmk:counter}).
The case $\LCD \cong E_8$-lattice is such an example 
(Remark \ref{rmk:e8}).

We should make a few remarks on the simplicity and the CFT type
property. Most of the vertex operator algebras discussed in this
paper are clearly simple and of CFT type. In such a case we omit
the proof of these properties.

This paper is the detailed version of our paper \cite{TY2}.

\section{Preliminaries}
Throughout this paper, $\zeta_n=\exp({2\pi\sqrt{-1}/n})$ is
a primitive $n$-th root of unity for a positive integer $n$. For
simplicity, $0,1$ and $2$ are sometimes understood to be elements
of $\Z_3$.

\subsection{Basic terminology}
Let $g$ be an automorphism of a vertex operator algebra 
$(V,Y,\1,\om)$ of finite order $T$. 
Set $V^r = \{ v \in V\,|\, gv = \zeta_{T}^rv\}$,
so that $V = \op_{r \in \Z/T\Z} V^r$.

For subsets $A, B$ of $V$ and a subset $X$ of a weak $g$-twisted $V$-module $M$,
set $A\cdot B=\spn_{\C}\{u_nv\ |\ u\in A, v\in B,n\in\Z\}$ and
$A\cdot X=\spn_{\C}\{u_nw\ |\ u\in A, w\in X,n\in(1/T)\Z\}$.
Then it follows that $(A\cdot B)\cdot X=A\cdot (B\cdot X)$ by
\cite[Lemma 3.12]{Li3} and \cite[Lemmas 2.5 and 2.6]{TY}.

Let $\N$ be the set of nonnegative integers.
A $(1/T)\N$-graded weak $g$-twisted $V$-module here is called
an admissible $g$-twisted $V$-module in \cite{DLM1}. Without loss we
can shift the grading of a $(1/T)\N$-graded weak $g$-twisted
$V$-module $M$ so that $M{(0)} \ne 0$ if $M \ne 0$. We call such an
$M{(0)}$ the \emph{top level} of $M$.

A vertex operator algebra $V$ is said to be \emph{rational} if every
$\N$-graded weak $V$-module is a direct sum of irreducible
$\N$-graded weak $V$-modules. 
If the dimension of the quotient space 
$V/\spn_{\C} \{ u_{-2}v\,|\, u,v \in V \}$ is finite, $V$ is said to
be $C_2$-\emph{cofinite} \cite{Z}. 
If $V=\op_{n=0}^\infty V_n$ and $V_0 = \C\1$, then $V$ is said to
be of \emph{CFT type}. Here $V_n=\{u\in V\ |\ \omega_1u=nu\}$ is
the homogeneous subspace of weight $n$.
If $V$ is $C_2$-cofinite and of CFT type, then 
the classification of irreducible $V$-modules
means the classification of irreducible weak $V$-modules 
\cite[Proposition 5.6 and Corollary 5.7]{ABD}.

For $h\in\Aut V$ and
a weak (resp. $(1/T)\N$-graded weak) $g$-twisted $V$-module $(M,Y_M)$,
we define a weak (resp. $(1/T)\N$-graded weak) $h^{-1}gh$-twisted
$V$-module $(M\circ h, Y_{M\circ h})$ by 
$M\circ h=M$ as vector spaces and $Y_{M\circ h}(u,x)=Y(hu,x)$.
If $M$ is irreducible, so is $M\circ h$.

Let $G$ be an automorphism group of $V$ and $V^G$ the vertex
operator subalgebra of $G$-invariants of $V$. 
A set ${\mathcal S}$ of irreducible $V$-modules is said to be
$G$-{\it stable} if for any $M\in{\mathcal S}$ and $h\in G$ there
exists $W\in{\mathcal S}$ such that $M\circ h\cong W$. An
irreducible $V$-module $M$ is said to be $G$-stable if $M\circ
g\cong M$ for all $g\in G$. 
It is shown in \cite[Theorem 4.4]{DM}
that if $V$ is simple and $G$ is of finite order, 
then $V^G$ is simple.

We denote by $I_V\binom{M^3}{M^1\ M^2}$ the set of all intertwining operators 
of type $\binom{M^3}{M^1\ M^2}$ \cite{FHL}. 
Let $\M$ be the set of all irreducible $V$-modules up to isomorphism 
and $\Z\M$ be a free $\Z$-module with basis $\M$.
For $M^1, M^2\in \M$,
\begin{equation*}
M^1\times M^2=\sum_{M^3\in\M}
\dim_{\C}I_{V}\binom{M^3}{M^1\ M^2}M^3\in\Z\M
\end{equation*}
is the fusion rule.
We write
$\sum_{M\in\M}S_{M}M\geq
\sum_{M\in\M}T_{M}M$ when $S_{M}\geq T_{M}$ for all $M\in \M$.

\subsection{Lattice $\LCD$}\label{subsec:latticeLCD}
We follow the notation in \cite{DLTYY, KLY1, KLY2,TY}.
Let $(L, \la\,\cdot\,,\,\cdot\,\ra)$ be
$\sqrt{2}$ times an ordinary root lattice of type $A_2$ and let
$\{\beta_1,\beta_2\}$ be a $\Z$-basis of $L$ such that $ \langle
\beta_1,\beta_1\rangle= \langle \beta_2,\beta_2\rangle=4$ and
$\langle \beta_1,\beta_2\rangle=-2$. Set $\beta_0=-\beta_1-\beta_2$.
Let $\otau$ be an isometry of $L$ induced by the permutation
$\beta_1 \mapsto \beta_2 \mapsto \beta_0 \mapsto \beta_1$. Then
$\otau$ is fixed-point-free and of order $3$. 

There are 12 cosets of $L$ in its dual lattice $L^{\perp} = \{ \al
\in \Q \otimes_\Z L\,|\, \la \al,L\ra \subset \Z\}$. These 12 cosets
are parametrized by $\Z_2 \times \Z_2$ and $\Z_3$. Let
$(L^{\perp})^{\oplus {\ell}}$ be an orthogonal sum of $\ell$ copies
of $L^\perp$. We shall construct a lattice $\LCD$ in
$(L^{\perp})^{\oplus {\ell}}$ from those 12 cosets of $L$ by using a
$\Z_2 \times \Z_2$-code $C$ and a $\Z_3$-code $D$.
We shall also introduce certain isometry groups of
$(L^{\perp})^{\oplus {\ell}}$.

First, $\otau$ can be extended to an isometry of $L^{\perp}$. Let
$\tgrs$ be a direct product of  ${\ell}$ copies of the group
$\langle\otau\rangle$ generated by $\otau$. Each element
$g=(g_1,\ldots,g_\ell)$ of $\tgrs$ transforms $\al=\al_1+\cdots +
\al_\ell \in (L^{\perp})^{\oplus {\ell}}$ as $g(\al) = g_1(\al_1) +
\cdots + g_\ell(\al_\ell)$, where $g_s \in \la \otau \ra$ and
$\al_s$ is the $s$-th component of $\al$. 
For convenience, we denote $(\otau,\ldots,\otau)\in \tgrs$ 
simply by $\tau$ also. 
A symmetric group $\mathfrak{S}_{\ell}$ of degree $\ell$ acts
on $(L^{\perp})^{\oplus {\ell}}$ by permuting the components. Let
$\tgrl$ be an isometry group of $(L^{\perp})^{\oplus {\ell}}$
generated by $\tgrs$ and $\mathfrak{S}_{\ell}$, which is a semidirect product
$\tgrs\rtimes \mathfrak{S}_{\ell}$ of $\tgrs$ by $\mathfrak{S}_{\ell}$.

Now, we discuss a $\Z_2 \times \Z_2$-code and a $\Z_3$-code. A $\Z_2
\times \Z_2$-code of length $\ell$ means an additive subgroup of
$\K^\ell$, where $\K = \{0,a,b,c\} \cong \Z_2 \times \Z_2$ is
Klein's four-group.
We call it a $\K$-code also.
Note that $b+c=a$ in $\K$.
For $x, y \in \K$,
define
\begin{equation*}
x \cdot y =
\begin{cases}
1 &\text{if }  x \ne y, x \ne 0, y \ne 0,\\
0 &\text{otherwise}.
\end{cases}
\end{equation*}

We have
\begin{equation}\label{eq:product_in_K}
x \cdot y \equiv m_1 n_2 + m_2 n_1 \quad \pmod{2\Z}
\end{equation}
if $x = m_1 c + m_2 b, y = n_1 c + n_2 b \in \K$ with
$m_1, m_2, n_1, n_2 \in \Z$.

For $\lm =
(\lm_1, \ldots, \lm_\ell)$, $\mu = (\mu_1, \ldots, \mu_\ell) \in
\K^\ell$, let $\langle\lm, \mu\rangle_{\Klein} = \sum_{i=1}^\ell
\lm_i \cdot \mu_i \in \Z_2$. The orthogonal form $(\lm, \mu) \mapsto
\langle\lm, \mu\rangle_{\Klein}$ on $\K^\ell$ was used in
\cite{KLY1, LY}. For a $\K$-code $C$ of length $\ell$, we define its
dual code by
\begin{equation*}
C^\perp = \{\lm \in \K^\ell \,|\, \langle\lm, \mu\rangle_{\Klein} =
0 \text{ for all } \mu \in C\}.
\end{equation*}

A $\K$-code $C$ is said to be {\em self-orthogonal} if $C \subset C^\perp$
and self-dual if $C = C^\perp$. For $\lm = (\lm_1,\ldots,\lm_\ell)
\in \Klein^{\ell}$, its {\em support} is defined to be
$\supp_{\Klein}(\lm) = \{ i\,|\, \lm_i \ne 0\}$. The cardinality of
$\supp_{\Klein}(\lm)$ is called the {\em weight} of $\lm$. We denote the
weight of $\lm$ by $\wt_{\Klein}(\lm)$.
In the case $\ell = 1$, we have $\wt_{\Klein}(x) = 0$ or $1$ according
to $x = 0$ or $x \in \{ a,b,c\}$.
A $\K$-code $C$ is said to
be {\em even} if $\wt_{\Klein}(\lm)$ is even for every $\lm \in C$.

We consider an action of $\otau$ on $\Klein$ such that $\otau(0) =
0$, $\otau(a) = b$, $\otau(b) = c$, and $\otau(c) = a$. Moreover, we
consider a componentwise action of $\tgrs$ on $\K^\ell$, so that
$\tau$ acts on $\Klein^{\ell}$ by $\tau(\lm_1, \ldots,\lm_\ell) =
(\otau(\lm_1),\ldots,\otau(\lm_\ell))$. 
Then $\tgrl$ acts on
$\Klein^{\ell}$ naturally. We denote by
$(\Klein^{\ell})_{\equiv_{\htau}}$ the set of all $\htau$-orbits in
$\Klein^{\ell}$.
For simplicity of notation, we sometimes denote a $\tau$-orbit
in $\Klein^{\ell}$ by its representative $\lambda\in\Klein^{\ell}$.

The first assertion of the next lemma is \cite[Lemma 2.8]{LY}.
The second assertion follows from the fact that
$\langle\lambda,\htau(\lambda)\rangle_{\Klein}
\equiv\wt_{\Klein}(\lambda)\pmod{2\Z}$
for $\lambda\in\Klein^{\ell}$.

\begin{lem}\label{lem:C-property}
Let $C$ be a $\K$-code of length $\ell$.

$(1)$ If $C$ is even, then $C$ is self-orthogonal.

$(2)$ If $C$ is $\tau$-invariant, then $C$ is even if and only if
$C$ is self-orthogonal.
\end{lem}

A $\Z_3$-code of length $\ell$ is a subspace of the vector space
$\Z_3^\ell$. For $\gm = (\gm_1,\ldots,\gm_\ell)$,
$\dl=(\dl_1,\ldots,\dl_\ell) \in \Z_3^\ell$, we consider the
ordinary inner product $\langle\gm, \dl\rangle_{\Z_3} =
\sum_{i=1}^\ell \gm_i \dl_i \in \Z_3$. The dual code $D^\perp$ of a
$\Z_3$-code $D$ is defined to be
\begin{equation*}
D^\perp = \{ \gm \in \Z_3^\ell\,|\, \langle\gm, \dl\rangle_{\Z_3} =
0 \text{ for all } \dl \in D\}.
\end{equation*}
Then $D$ is said to be {\em self-orthogonal} if $D \subset D^\perp$ and
{\em self-dual} if $D=D^\perp$.

We define the support and the weight of $\gm =
(\gm_1,\ldots,\gm_\ell) \in \Z_3^\ell$ in the same way as before.
Thus $\supp_{\Z_3}(\gm) = \{ i\,|\, \gm_i \ne 0\}$ and
$\wt_{\Z_3}(\gm)$ is the cardinality of $\supp_{\Z_3}(\gm)$. Note
that $\wt_{\Z_3}(\gm) \equiv \la \gamma, \gamma\ra \pmod{3\Z}$. 
Then the following lemma holds.

\begin{lem}\label{lem:coset}
Let $D$ be a self-orthogonal $\Z_3$-code of length $\ell$. Then
$\wt_{\Z_3}(\delta - \gamma) \equiv \wt_{\Z_3}(\delta) \pmod{3\Z}$
for any $\gamma \in D$ and $\delta \in D^\perp$.
\end{lem}

We consider the trivial action of $\otau$ on $\Z_3$, that is,
$\otau(j) = j$ for $j \in \Z_3$. Then $\tgrs$ acts trivially on
$\Z_3^\ell$
and 
$\tgrl$ acts on $\Z_3^{\ell}$ naturally. 

Take a $\Z$-basis $\tilde{\beta}_1=\beta_1/2,
\tilde{\beta}_2=(\beta_1-\beta_2)/6$ of $L^{\perp}$. Note that
$\{2\tilde{\beta}_1, 6\tilde{\beta}_2\}$ is a $\Z$-basis of $L$.
For $\alpha=m_1\tilde{\beta}_1+m_2\tilde{\beta}_2,
\beta=n_1\tilde{\beta}_1+n_2\tilde{\beta}_2\in L^{\perp}$,
we have
\begin{align}\label{eq:stan-inner}
\langle
\alpha,\beta\rangle
&=m_1n_1+\dfrac{m_1n_2+m_2n_1}{2}+
\dfrac{m_2n_2}{3}.
\end{align}
We also have
$\otau(\tilde{\beta}_1) = \tilde{\beta}_1 - 3\tilde{\beta}_2$ and
$\otau(\tilde{\beta}_2) = \tilde{\beta}_1 - 2\tilde{\beta}_2$.
We use the same notation as in \cite{DLTYY, KLY1, KLY2,TY} to denote
the $12$ cosets $L^{(x,i)}$, $x \in \K$, $i \in \Z_3$ of $L$ in its
dual lattice $L^\perp$. For each $x \in \K$ we assign $\be(x) \in
L^\perp$ by $\be(0) = 0$, $\be(a) = \be_2/2$, $\be(b) = \be_0/2$,
and $\be(c) = \be_1/2$. Then
\begin{equation}\label{eq:coset-lattice}
L^{(x,i)} = \beta(x) + i\frac{-\beta_1+\beta_2}{3} + L.
\end{equation}

Since $\tilde{\beta}_1 = \beta(c) \in L^{(c,0)}$ and
$\tilde{\beta}_2 = \beta(b) + (-\beta_1+\beta_2)/3 + \beta_1 \in L^{(b,1)}$, 
we can describe $L^{(x,i)}$ by using the basis
$\{\tilde{\beta}_1, \tilde{\beta}_2\}$ of $L^{\perp}$.

\begin{lem}\label{lem:coset-b1b2}
For $x \in \K$ and $i \in \Z_3$,
\begin{equation*}
L^{(x,i)} =\{
m_1\tilde{\beta}_1 + m_2\tilde{\beta}_2 \in L^{\perp}\, |\,
x = m_1c + m_2b \text{ in } \K \text{ and } i = m_2 + 3\Z\}.
\end{equation*}
\end{lem}

We also have the following lemma.

\begin{lem}\label{lem:coset-b1b2-2}
Let $\alpha\in L^{(x,i)}$ and $\beta\in L^{(y,j)}$
with $x,y \in \K$, $i,j \in \Z_3$.

$(1)$ $\langle\alpha,\beta\rangle \equiv x\cdot y/2+ij/3\pmod{\Z}$.

$(2)$ $\langle\alpha,\alpha\rangle \equiv \wt_{\Klein}(x)-2i^2/3\pmod{2\Z}$.
\end{lem}

For $\lm = (\lm_1,\ldots,\lm_\ell) \in \K^\ell$ and $\gm
=(\gm_1,\ldots,\gm_\ell) \in \Z_3^\ell$, let
\begin{equation*}
L_{(\lm,\gm)} = L^{(\lm_1,\gm_1)} \oplus \cdots \oplus
L^{(\lm_\ell,\gm_\ell)} \subset (L^\perp)^{\oplus \ell}.
\end{equation*}
Moreover, for $\mu\in\Klein^{\ell},\delta\in\Z_3^{\ell}$, $P\subset
\Klein^{\ell}$, and $Q\subset \Z_3^{\ell}$, set
\begin{equation*}
L_{\mu\times Q} = \bigcup_{\gm \in Q} L_{(\mu,\gm)},\quad L_{P
\times \delta} = \bigcup_{\lm \in P} L_{(\lm,\delta)},\quad L_{P
\times Q} = \bigcup_{\lm \in P, \gm \in Q} L_{(\lm,\gm)}.
\end{equation*}

For a $\Klein$-code $C$ of length ${\ell}$ and a $\Z_3$-code $D$ of
the same length, $L_{C \times D}$ is an additive subgroup of
$(L^\perp)^{\oplus \ell}$. However, $L_{C \times D}$ is not an
integral lattice in general. In the case where $C=\K^\ell$ and $D =
\Z_3^\ell$, $L_{C \times D}$ coincides with $(L^\perp)^{\oplus
\ell}$. If $C = \{\zero \}$ and $D = \{\zero\}$, then $L_{\{\zero\}
\times \{\zero\}}=L_{(\zero,\zero)}= L^{\oplus \ell}$, where $\zero
= (0,\ldots,0)$. 
In the case of $\ell = 1$, we
note that $L_{\Klein\times\zero}=\Z\tilde{\beta}_1+\Z(3\tilde{\beta}_2)$,
$L_{\zero\times\Z_3}=\Z(2\tilde{\beta}_1)+\Z(2\tilde{\beta}_2)$, and
$L = L_{\zero\times\zero} = \Z(2\tilde{\beta}_1)+\Z(6\tilde{\beta}_2)$.

Let $(L_{C \times D})^\perp =\{ \al \in (\Q \ot_\Z L)^{\oplus \ell}
\,|\, \la \al, L_{C \times D}\ra \subset \Z\}$. 
The following lemma is a consequence of Lemma \ref{lem:coset-b1b2-2} (1).

\begin{lem}\label{lem:dual-lattice}
$(L_{C \times D})^\perp = L_{C^\perp \times D^\perp}$.
\end{lem}

Thus $L_{C \times D}$ is an integral
lattice if and only if both of $C$ and $D$ are self-orthogonal. The
first assertion of the next lemma follows from
Lemma \ref{lem:coset-b1b2-2} (2). The second assertion is a special
case of the above lemma (see also \cite[Theorems 5.6, 5.7]{KLY1}).

\begin{lem}\label{lem:even-lattice}
$(1)$ If $C$ is even and $D$ is self-orthogonal,
then $L_{C \times D}$ is an even lattice.

$(2)$ If $C$ and $D$ are self-dual, then  $L_{C \times D}$ is a
unimodular lattice.
\end{lem}

\subsection{Central extensions $\hat{L}_{C\times D}$,
$\hat{L}_{C\times D,\htau^i}$, $i = 1,2$}
\label{subsetion:group-twist}
Suppose $C$ is a $\tau$-invariant even $\K$-code of length $\ell$
and $D$ is a self-orthogonal $\Z_3$-code of the same length. Then
$\LCD$ is a positive definite even lattice by Lemma \ref{lem:even-lattice}. 
The isometry $\otau$ of $L^\perp$ permutes 
the cosets $L^{(x,i)}$, $x \in \K$, $i \in \Z_3$ of $L$ in $L^\perp$. 
In fact, $\otau(L^{(x,i)}) = L^{(\otau(x),i)}$ 
by our definition of the action of $\otau$ on $L^\perp$, $\K$ and
$\Z_3$ introduced in Section \ref{subsec:latticeLCD}.
In particular, $\tau$ induces an
isometry of $\LCD$, for we are assuming that $C$ is
$\tau$-invariant. Note that $\tau$ is fixed-point-free on $\LCD$. 
We also have $g(L_{C\times D})=L_{g(C)\times g(D)}$ for $g\in\tgrl$. 

For any positive integer $n$, let $\la \kappa_n \ra$ be a cyclic
group of order $n$ with generator $\kappa_n$. We assume that
${\kappa_n}^{n/m} = \kappa_m$ if $m$ is a divisor of $n$.
We shall construct three central extensions $\widehat{L}_{C\times D}$ and
$\hat{L}_{C\times D,\htau^i}$, $i = 1,2$ of
$L_{C\times D}$ by $\langle\kz\rangle$ which will be used
in later sections.
We realize each of these central extensions
as a subgroup of a central extension of $(L^{\perp})^{\oplus \ell}$
by $\langle\kz\rangle$.

Define $\Z$-bilinear forms $\vep_1, \vep_2, \vep_2^{\prime},
c_1, c_2, c_2^{\prime} \colon L^\perp \times L^\perp \rightarrow \Z/36\Z$
as follows.
For $\alpha=m_1\tilde{\beta}_1+m_2\tilde{\beta}_2$,
$\beta=n_1\tilde{\beta}_1+n_2\tilde{\beta}_2\in L^{\perp}$, set
\begin{align}
\varepsilon_1(\alpha,\beta)
& = 27m_1n_1+27m_2n_1+9m_2n_2+36\Z,\label{eq:vep1}\\
\varepsilon_2(\alpha,\beta)
& = 6m_1n_1+6m_2n_1+14m_2n_2+36\Z,\label{eq:vep2}\\
\varepsilon_2^{\prime}(\alpha,\beta)
& = 6m_1n_1+15m_1n_2+27m_2n_1+14m_2n_2+36\Z,\label{eq:vept2}
\end{align}
and
\begin{align}
c_1(\alpha,\beta)
& = \varepsilon_1(\alpha,\beta)-\varepsilon_1(\beta,\alpha) =
9m_1n_2+27m_2n_1+36\Z,\label{eq:c1}\\
c_2(\alpha,\beta)
& = \varepsilon_2(\alpha,\beta)-\varepsilon_2(\beta,\alpha) =
30m_1n_2+6m_2n_1+36\Z,\label{eq:c2}\\
c_2^{\prime}(\alpha,\beta)
& = \varepsilon_2^{\prime}(\alpha,\beta)-\varepsilon_2^{\prime}(\beta,\alpha) =
24m_1n_2+12m_2n_1+36\Z.\label{eq:ct2}
\end{align}
We also set
\begin{align}
c_{0}(\alpha,\beta)&=18\langle\alpha,\beta\rangle+36\Z\nonumber\\
&=18m_1n_1+9m_1n_2+9m_2n_1+6m_2n_2+36\Z.\label{eq:c0}
\end{align}

All of these $\Z$-bilinear forms are $\otau$-invariant.
Since $\vep_1$ is $\Z$-bilinear, it is a 2-cocycle.
Let $\widehat{L^{\perp}} =
\langle\kz\rangle\times L^{\perp}$. 
We simply write $\kz^p e^{\alpha}$ for $(\kz^p,\alpha)\in
\widehat{L^{\perp}}$. 
In particular, $\kz^p = (\kz^p,0)$ and $e^\alpha = (1,\alpha)$. 
Define a multiplication on the set $\widehat{L^{\perp}}$ by 
\begin{equation}\label{eq:mul1}
(\kz^p e^{\alpha}) \cdot (\kz^q e^{\beta}) =
\kz^{p+q+\varepsilon_1(\alpha,\beta)} e^{\alpha+\beta}.
\end{equation}

Take $\varepsilon_2$ (resp. $\varepsilon_2^{\prime}$)
in place of $\varepsilon_1$.
Then we obtain a multiplicative group
${\widehat{L^{\perp}}}_\otau$ (resp. ${\widehat{L^{\perp}}}_{\otau^2}$).
We use the same notation $\kz^p e^{\alpha}$ to denote its element.
As to its multiplication, we write $\times_\otau$
(resp. $\times_{\otau^2}$) so that
\begin{align}
(\kz^pe^{\alpha}) \times_{\otau} (\kz^qe^{\beta}) &=
\kz^{p+q+\varepsilon_2(\alpha,\beta)}e^{\alpha+\beta},\label{eq:mul2}\\
(\kz^pe^{\alpha}) \times_{\otau^2} (\kz^qe^{\beta}) &=
\kz^{p+q+\varepsilon_2^{\prime}(\alpha,\beta)}e^{\alpha+\beta}\label{eq:mul3}.
\end{align}

For $a, b \in \widehat{L^{\perp}}$ or  ${\widehat{L^{\perp}}}_{\otau^i}$,
$i = 1,2$, we simply write $ab$ for the product in the group
when there is no ambiguity.
Define $^{-}\colon \widehat{L^{\perp}} \rightarrow L^\perp$
(resp. ${\widehat{L^{\perp}}}_{\otau^i} \rightarrow L^\perp$) by
$\overline{\kz^{p}e^{\alpha}} = \alpha$.
Then $\widehat{L^{\perp}}$ (resp. ${\widehat{L^{\perp}}}_{\otau}$,
${\widehat{L^{\perp}}}_{\otau^2}$) is a central extension of
$L^{\perp}$ by $\langle\kz\rangle$ with
associated commutator map $c_1$ (resp. $c_2$, $c_2^{\prime}$) 
(\cite[Sections 5.1, 5.2]{FLM}, \cite[Section 6.4]{LL}).

Note that
\begin{equation}\label{eq:prod-in-LD}
e^\alpha e^\beta = e^{\alpha + \beta} \mbox{ in } \widehat{L^{\perp}}
\end{equation}
for $\alpha, \beta \in L_{\zero\times\Z_3} =
\Z(2\tilde{\beta}_1) + \Z(2\tilde{\beta}_2)$ by \eqref{eq:mul1}.

Define an automorphism of the group $\widehat{L^{\perp}}$
(resp. ${\widehat{L^{\perp}}}_{\otau^i}$, $i=1,2$) of order $3$ by
\begin{equation}\label{eq:def-tau}
\begin{split}
\kz &\mapsto \kz,\\
e^{\alpha} &\mapsto
e^{\otau(\alpha)}
\end{split}
\end{equation}
for $\alpha\in L^{\perp}$.
Since $\varepsilon_1$ (resp. $\varepsilon_2$,
$\varepsilon_2^{\prime}$) is $\otau$-invariant,
the map is in fact an automorphism
of the group $\widehat{L^{\perp}}$
(resp. $\widehat{L^{\perp}_{\otau}}$,
${\widehat{L^{\perp}}}_{\otau^2}$) of order $3$. 
By abuse of notation, we denote it by $\otau$ also.

\begin{rmk}\label{rmk:relation_with_Dong_Lepowsky}
In \cite[Remark 2.2]{DL}, three bilinear forms
$\varepsilon_0$, $c_{0}$ and $c_{0}^{\nu}$ were
considered. Apply \cite[(2.9), (2.10), (2.13)]{DL}
to $L^{\perp}$ in place of $L$
with $\nu = \otau$ or $\otau^2$, $p = 3$
and $q = 36$. Then the bilinear form $c^0$ of \cite[(2.9)]{DL}
is identical with our $c^0$. 
Moreover, $\varepsilon_0$ and $c_{0}^{\nu}$ become
\begin{align}
\varepsilon_0(\alpha,\beta)
&=30\langle\alpha,\otau(\beta)\rangle+36\Z\nonumber\\
&=21m_1n_1+21m_2n_1+31m_2n_2+36\Z,\label{eq:vep0}\\
\varepsilon_0^{\prime}(\alpha,\beta)
&=30\langle\alpha,\otau^2(\beta)\rangle+36\Z\nonumber\\
&=21m_1n_1+21m_1n_2+31m_2n_2+36\Z,\label{eq:vept0}\\
c_{0}^{\otau}(\alpha,\beta)
&=12\langle\otau(\alpha)+2\otau^2(\alpha),\beta\rangle+36\Z\nonumber\\
&=18m_1n_1+30m_1n_2+24m_2n_1+30m_2n_2+36\Z, \label{eq:ctau}\\
c_{0}^{\otau^2}(\alpha,\beta)
&=12\langle\otau^2(\alpha)+2\otau^4(\alpha),\beta\rangle+36\Z\nonumber\\
&=18m_1n_1+24m_1n_2+30m_2n_1+30m_2n_2+36\Z \label{eq:ctaut}
\end{align}
for $\alpha=m_1\tilde{\beta}_1+m_2\tilde{\beta}_2$,
$\beta=n_1\tilde{\beta}_1+n_2\tilde{\beta}_2 \in L^{\perp}$.
Here we write $\varepsilon_0^{\prime}$ for $\varepsilon_0$ of
\cite[(2.13)]{DL} in the case $\nu = \otau^2$. 
Note that $q$ of \cite[Remark 2.2]{DL} should be a multiple of $12$ 
by \eqref{eq:c0}. 
We take $q = 36$ so that every coefficient of $m_i n_j$ in
\eqref{eq:vep0} and \eqref{eq:vept0} is an integer.
These bilinear forms are related to our ones as follows.
\begin{align*}
\varepsilon_0(\alpha,\beta) &=
\varepsilon_1(\alpha,\beta)-\varepsilon_2(\alpha,\beta),\\
\varepsilon_0^{\prime}(\alpha,\beta) &=
\varepsilon_1(\alpha,\beta)-\varepsilon_2^{\prime}(\alpha,\beta),\\
c_{0}(\alpha,\beta) &\equiv
c_1(\alpha,\beta)-36\langle \alpha,\htau(\beta)\rangle \pmod{36\Z},\\
c_{0}^{\htau}(\alpha,\beta) &\equiv
c_2(\alpha,\beta)+36\langle \alpha,\htau(\beta)\rangle \pmod{36\Z},\\
c_{0}^{\htau^2}(\alpha,\beta) &\equiv
c_2^{\prime}(\alpha,\beta)+36\langle \alpha,\htau(\beta)\rangle \pmod{36\Z}.
\end{align*}
\end{rmk}

We extend the $\Z$-bilinear forms
$\vep_1, \vep_2, \vep_2^{\prime}, c_1, c_2, c_2^{\prime}, c_0,
\vep_0, \vep_0^{\prime}, c_{0}^{\htau}, c_{0}^{\htau^2}$
on $L^\perp$ to $(L^{\perp})^{\oplus {\ell}}$ naturally. For example,
\begin{equation*}
\vep_1(\alpha,\beta) = \sum_{s=1}^\ell \vep_1(\alpha^{(s)},\beta^{(s)})
\end{equation*}
for $\alpha = \sum_{s=1}^\ell \alpha^{(s)}$, $\beta =
\sum_{s=1}^\ell \beta^{(s)} \in (L^{\perp})^{\oplus {\ell}}$, where
$\alpha^{(s)}$ and $\beta^{(s)}$ are in the $s$-th entry of
$(L^{\perp})^{\oplus {\ell}}$. 
These $\Z$-bilinear forms are all $\htau$-invariant.

\begin{rmk}\label{rmk:relation_with_Dong_Lepowsky_2}
If $\langle\alpha,\htau(\beta)\rangle \in \Z$,
then Remark \ref{rmk:relation_with_Dong_Lepowsky} implies that
$c_1(\alpha,\beta)= c_0(\alpha,\beta)$,
$c_2(\alpha,\beta)=c_{0}^{\htau}(\alpha,\beta)$,
$c_2^{\prime}(\alpha,\beta)=c_{0}^{\htau^2}(\alpha,\beta)$
and
\begin{equation}\label{eq:ve-relation}
\begin{split}
\varepsilon_0(\alpha,\beta)-\varepsilon_0(\beta,\alpha)
&=c_0(\alpha,\beta)-c_0^{\htau}(\alpha,\beta),\\
\varepsilon_0^{\prime}(\alpha,\beta)-\varepsilon_0^{\prime}(\beta,\alpha)
&=c_0(\alpha,\beta)-c_0^{\htau^2}(\alpha,\beta).
\end{split}
\end{equation}
\end{rmk}

Let $(\widehat{L^{\perp}})^{\ell}$ be a direct product of ${\ell}$
copies of $\widehat{L^{\perp}}$ and let
$T$ be a subgroup in the
center of $(\widehat{L^{\perp}})^{\ell}$ generated by
$\kz^{(r)}(\kz^{(s)})^{-1}$, $1\leq r,s\leq \ell$, where $\kz^{(s)}$
denotes $\kz \in L^{\perp}$ in the $s$-th entry of
$(\widehat{L^{\perp}})^{\ell}$. We consider 
$(\widehat{L^{\perp}})^{\ell}/T$. For simplicity of notation, we
write $e^{\alpha_1+\cdots+\alpha_{\ell}}$ for
$(e^{\alpha_1},\ldots,e^{\alpha_{\ell}})T$ and $\kz^p$ for
$(\kz^{(1)})^pT$ in $(\widehat{L^{\perp}})^{\ell}/T$. 
Then any element of $(\widehat{L^{\perp}})^{\ell}/T$ can be
expressed uniquely in the form $\kz^{p}e^{\alpha}$ with
$p\in\Z/36\Z$ and $\alpha\in (L^{\perp})^{\oplus {\ell}}$.

By \eqref{eq:mul1} we have
\begin{equation}\label{eq:multiplication-inmcL}
e^{\alpha}e^{\beta} =
\kz^{\varepsilon_1(\alpha,\beta)}e^{\alpha+\beta}
\end{equation}
in $(\widehat{L^{\perp}})^{\ell}/T$. 
For $\kz^{p}e^{\alpha}\in (\widehat{L^{\perp}})^{\ell}/T$, let
$\overline{\kz^{p}e^{\alpha}}=\alpha\in (L^{\perp})^{\oplus \ell}$.
Then
\begin{align*}
1 \rightarrow \langle\kz \rangle \rightarrow (\widehat{L^{\perp}})^{\ell}/T 
\xrightarrow{-} (L^{\perp})^{\oplus {\ell}} \rightarrow 1
\end{align*}
is a central extension of $(L^{\perp})^{\oplus {\ell}}$ by
$\langle\kz\rangle$ with associated commutator map $c_1$. 
We denote $(\widehat{L^{\perp}})^{\ell}/T$ by $\mcL$ also.

By \eqref{eq:def-tau},
$\tgrl$ acts on the group $\mcL$ naturally. In
particular, $\htau = (\otau,\ldots,\otau)$ acts on $\mcL$ as 
an automorphism of order $3$. 
We have $\overline{g(a)}=g(\bar{a})$ 
for $g\in \tgrl$ and $a\in \mcL$.

By \eqref{eq:c0} and Remark \ref{rmk:relation_with_Dong_Lepowsky}, we have
\begin{equation*}
\kz^{c_1(\alpha,\beta)} = \kz^{c_0(\alpha,\beta)}
= \kappa_2^{\langle \alpha,\beta\rangle}
\end{equation*}
if $\langle \alpha, \htau(\beta)\rangle$ is an integer.
This is the case for $\alpha, \beta \in \LCD$,
since $\LCD$ is a $\tau$-invariant integral lattice.

For any subset $Q$ of $(L^{\perp})^{\oplus {\ell}}$, we set 
$\hat{Q} = \{a\in \mcL\, |\, \bar{a}\in Q\}$. In particular, 
$\hat{L}_{C\times D}=\{a\in \mcL\, |\, \bar{a}\in L_{C\times D}\}$. 
Then
\begin{equation}
1 \rightarrow \langle\kz \rangle \rightarrow \hat{L}_{C\times D}
\xrightarrow{-} \LCD \rightarrow 1
\label{eq:cent36}
\end{equation}
is a central extension of $\LCD$ by $\langle\kz \rangle$ with
associated commutator map $c_1$.

Replace $\widehat{L^\perp}$ with ${\widehat{L^\perp}}_{\otau}$
(resp. ${\widehat{L^\perp}}_{\otau^2}$) and $\vep_1$ with
$\vep_2$ (resp. $\vep_2^{\prime}$) in the above argument. Then we
obtain a central extension $\widehat{(L^{\perp})^{\oplus {\ell}}}_\tau$
(resp. $\widehat{(L^{\perp})^{\oplus {\ell}}}_{\tau^2}$) of 
$(L^{\perp})^{\oplus {\ell}}$ by
$\langle\kz \rangle$ with
associated commutator map $c_2$ (resp. $c_2^{\prime}$). 
We have $e^{\alpha}e^{\beta} =
\kz^{\varepsilon_2(\alpha,\beta)}e^{\alpha+\beta}$ in 
$\widehat{(L^{\perp})^{\oplus {\ell}}}_\tau$ by \eqref{eq:mul2}
(resp. $e^{\alpha}e^{\beta} =
\kz^{\varepsilon_2^{\prime}(\alpha,\beta)}e^{\alpha+\beta}$ in
$\widehat{(L^{\perp})^{\oplus {\ell}}}_{\tau^2}$ by \eqref{eq:mul3}) 
for $\alpha, \beta \in (L^{\perp})^{\oplus {\ell}}$. 
We also consider 
$\hat{Q}_{\tau^i} = 
\{ a \in \widehat{(L^{\perp})^{\oplus {\ell}}}_{\tau^i} \,|\, 
\bar{a}\in Q\}$, $i = 1,2$ 
similarly for a subset $Q$ of $(L^{\perp})^{\oplus {\ell}}$. 

Note that  $\tau$ induces an automorphism of
$\hat{L}_{C\times D}$ of order $3$.
Let $\theta \in \Aut \hat{L}_{C\times D}$ be a distinguished lift of
the isometry $-1$ of $L_{C\times D}$ defined by \cite[(10.3.12)]{FLM}
\begin{equation}\label{eq:theta}
\theta \colon \hat{L}_{C\times D}\rightarrow \hat{L}_{C\times D};\quad
a\mapsto a^{-1}\kappa_2^{\langle\bar{a},\bar{a}\rangle/2}.
\end{equation}
Then $\theta^2=1, \overline{\theta(a)}=-\bar{a}$ for
$a\in\hat{L}_{C\times D}$, and $\theta(\kz)=\kz$. 
Moreover, $\htheta \htau = \htau \htheta$
since $\la \,\cdot\, ,\,\cdot\,\ra$ is $\tau$-invariant. 
Thus we have obtained the following lemma.

\begin{lem}\label{lem:c-ex-La}
$\hat{L}_{C\times D}$ is a central extension of $L_{C\times D}$ by
$\langle\kz\rangle$ with commutation relation
\begin{equation}
ab = \kappa_2^{\langle \bar{a},\bar{b}\rangle}ba,\quad a,b\in
\hat{L}_{C\times D}.
\end{equation}
Moreover, $\tau$ and $\theta$ are automorphisms of $\hat{L}_{C\times
D}$ such that $\tau^3=\theta^2=1$, $\tau(\kz)=\theta(\kz)=\kz$,
$\overline{\tau(a)}=\tau(\bar{a})$, $\overline{\theta(a)}=-\bar{a}$
for $a\in\hat{L}_{C\times D}$, and $\theta\tau=\tau\theta$.
\end{lem}

The sublattice $L_{\zero\times D}$ of $\LCD$ has nice properties.
For $\alpha,\beta\in L_{\zero\times D}$, we have $e^\alpha e^\beta =
e^{\alpha+\beta}$ by \eqref{eq:prod-in-LD} and $\tau(e^\alpha) =
e^{\tau(\alpha)}$ by \eqref{eq:def-tau}. Furthermore,
$\theta(e^{\alpha})=e^{-\alpha}$ by \eqref{eq:theta}, since $\la
\alpha,\alpha\ra \in 4\Z$ for $\alpha \in L_{\zero\times D}$.

Now, set $\C\{(L^{\perp})^{\oplus
\ell}\}=\C[\mcL]/(\kz-\zeta_{36})\C[\mcL]$, which is a twisted group
algebra of $(L^{\perp})^{\oplus \ell}$. By abuse of notation, we
denote the image of $e^\alpha \in \mcL$ in $\C\{(L^{\perp})^{\oplus
\ell}\}$ by the same symbol $e^\alpha$ for $\alpha \in
(L^\perp)^{\oplus \ell}$. The automorphisms $\tau$ and $\theta$ also
induce automorphisms of $\C\{(L^{\perp})^{\oplus \ell}\}$. We use
the same symbols $\tau$ and $\theta$ to denote those automorphisms. 
For any subset $P$ of $(L^{\perp})^{\oplus \ell}$, we set
$\C\{P\} = \spn_{\C}\{e^{\alpha}\, |\, \alpha\in P\}
\subset\C\{(L^{\perp})^{\oplus l}\}$.

The following lemma
is a direct consequence of Lemma \ref{lem:c-ex-La}.

\begin{lem}\label{lem:tw-l-gr}
$\C\{L_{C\times D}\}$ is a twisted group algebra of $L_{C\times D}$
such that
\begin{equation*}
e^{\alpha}e^{\beta}=(-1)^{\langle\alpha,\beta\rangle}e^{\beta}e^{\alpha},
\quad \alpha,\beta\in L_{C\times D}.
\end{equation*}
Moreover, $\tau$ and $\theta$ are automorphisms of $\C\{L_{C\times
D}\}$ such that $\tau^3=\theta^2=1$ and $\theta\tau=\tau\theta$.
\end{lem}

\subsection{Vertex operator algebra $V_{L_{C\times D}}$}\label{subsec:VOA_VLCD}
We use the standard notation for the vertex operator
algebra $(V_{\Lat},Y(\,\cdot\,,x))$ associated with a positive
definite even lattice $\Lat$ and its module $V_{\Lat^\perp}$
(\cite[Chapter 8]{FLM}, \cite[Section 6.4]{LL}).
Let
$C$ be a $\tau$-invariant even $\Klein$-code of length ${\ell}$ and
$D$ be a self-orthogonal $\Z_3$-code of the same length. Thus the
lattice $L_{C\times D}$ is a $\tau$-invariant positive definite even
lattice by Lemma \ref{lem:even-lattice}.
We use the twisted group algebra $\C\{L_{C\times D}\}$ of 
Lemma \ref{lem:tw-l-gr} for the vertex operator algebra
$V_{L_{C\times D}} = M(1) \otimes \C\{L_{C\times D}\}$.
We identify $V_{L^{\oplus\ell}}$ with
$V_{L}^{\otimes\ell}$ and $V_{(L^{\perp})^{\oplus\ell}}$ with
$V_{L^{\perp}}^{\otimes\ell}$.

Recall the action of the group
$\tgrl$ on $(L^\perp)^{\oplus \ell}$,
$\K^\ell$ and $\Z_3^\ell$ discussed in Section 2.2. 
For $g \in \tgrl$, define a linear isomorphism on
$V_{(L^{\perp})^{\oplus\ell}} = M(1) \otimes \C\{(L^{\perp})^{\oplus\ell}\}$
by
\begin{equation*}
\alpha^1(-n_1)\cdots \alpha^k(-n_k)e^{\beta} \mapsto
(g\alpha^1)(-n_1)\cdots (g\alpha^k)(-n_k)g(e^{\beta}).
\end{equation*}
For simplicity of notation, we denote it by $g$ also. Then
\begin{equation*}
g(Y_{(L_{C\times D})^{\perp}}(u,x)v) = Y_{(L_{g(C)\times
g(D)})^{\perp}}(gu,x)gv
\end{equation*}
for $u\in V_{L_{C\times D}}$ and $v\in V_{(L_{C\times D})^{\perp}}$.
Hence $g \colon V_{L_{C\times D}}\mapsto V_{L_{g(C)\times g(D)}}$ is an
isomorphism of vertex operator algebras. 
In particular, $\tau$ is an
automorphism of $V_{L_{C\times D}}$. Our purpose is the
classification of irreducible modules for the fixed point subalgebra
$V_{L_{C\times D}}^{\htau}= \{u\in V_{L_{C\times D}}\, |\, \tau
u=u\}$ of $V_{L_{C\times D}}$ by the automorphism $\tau$.

We also note that
\begin{equation*}
g \colon V_{L_{(\lambda+C)\times (\gamma+D)}}\mapsto
V_{L_{(g(\lambda)+g(C))\times (g(\gamma)+g(D))}}
\end{equation*}
for $\lambda\in C^{\perp}$ and $\gamma\in D^{\perp}$ is a map
from $V_{L_{C\times D}}$-modules to $V_{L_{g(C)\times
g(D)}}$-modules.
In the case where $C$ and $D$ are
$g$-invariant, we have
\begin{equation}\label{eq:actionD-un}
V_{L_{(\lambda+C)\times(\gamma+D)}}~\circ~ g\cong 
g^{-1}\big(V_{L_{(\lambda+C)\times(\gamma+D)}})
= V_{L_{(g^{-1}(\lambda)+C)\times(g^{-1}(\gamma)+D)}}.
\end{equation}

By \cite[Theorem 3.1]{D1} and Lemma \ref{lem:dual-lattice}, we
have the following proposition.
\begin{prop}
$\{V_{L_{(\lambda+C)\times(\gamma+D)}}\, |\, \lambda + C \in
C^{\perp}/C, \gamma + D \in D^{\perp}/D\}$ is a set of all
irreducible $V_{\LCD}$-modules up to isomorphism.
\end{prop}

The following lemma is a straightforward consequence of 
\eqref{eq:actionD-un}.

\begin{lem}
We have $V_{L_{(\lambda+C)\times (\gamma+D)}}\circ
\tau\cong V_{L_{(\tau^{-1}(\lambda)+C)\times (\gamma+D)}}$. In particular,
$V_{L_{(\lambda+C)\times (\gamma+D)}}$ is $\tau$-stable if and only if
$\lambda\in C$.
\end{lem}

For $\varepsilon=0,1,2$, let $V_{L_{C\times(\gamma+D)}}(\varepsilon) 
= \{u\in V_{L_{C\times (\gamma+D)}}\ |\ \htau
u=\zeta_3^{\varepsilon}u\}$. These are irreducible
$\VLCDtau$-modules.

The following proposition is clear.
\begin{lem}\label{lem:lattice-st}
As $(V_{L})^{\otimes {\ell}}$-modules, we have
\begin{equation*}
V_{L_{(\lambda+C)\times(\gamma+D)}} =
\bigoplus_{\mu\in\lambda+C,\delta\in \gamma+D}V_{L_{(\mu,\delta)}}.
\end{equation*}
\end{lem}

The fusion rules for $V_{L_{C\times D}}$ are 
known by \cite[Corollary 12.10]{DL}.

\begin{lem}\label{lem:fusion-lattice}
For $\lambda^1,\lambda^2\in C^{\perp}$ and $\gamma^1,\gamma^2\in
D^{\perp}$, we have
\begin{align*}
V_{L_{(\lambda^1+C)\times (\gamma^1+D)}}\times
V_{L_{(\lambda^2+C)\times (\gamma^2+D)}}&=
V_{L_{(\lambda^1+\lambda^2+C)\times (\gamma^1+\gamma^2+D)}}.
\end{align*}
\end{lem}

\section{Irreducible $\htau^i$-twisted $V_{\LCD}$-modules, $i=1,2$}
As before, we assume that $C$ is a $\tau$-invariant even
$\K$-code of length $\ell$ and $D$ is a self-orthogonal
$\Z_3$-code of the same length. 
We shall describe a decomposition of every irreducible $\htau^i$-twisted
$V_{\LCD}$-module constructed by the method of \cite{DL2, L} 
into a direct sum of
irreducible $(V_L^{\otau})^{\otimes \ell}$-modules, $i=1,2$.
The argument in the $\tau^2$-twisted
case is parallel to that in the $\tau$-twisted case. Thus we deal
with mainly the $\tau$-twisted ones.

By our construction 
$\widehat{L^{\oplus\ell}}$ 
(resp. ${\widehat{L^{\oplus\ell}}_{\htau}}$) is a subgroup
of $\hat{L}_{C\times D}$ (resp. $\hat{L}_{C\times D,\htau}$).
In \cite{DLTYY,TY}, we have
considered irreducible $\otau$-twisted $V_L$-modules
$V_{L}^{T_{\chi_j}}(\otau),j=0,1,2$. 
In order to apply the
results obtained in these previous papers, we need to examine the
relation between $\hat{L}$ (resp. $\hat{L}_{\tau}$) of
\cite{DLTYY} and $\widehat{L^{\oplus\ell}}$ (resp.
${\widehat{L^{\oplus\ell}}_{\htau}}$).

In \cite[(2.1)]{DLTYY}, $\hat{L}$ was a central extension of $L$
by $\langle\kappa_{6}\rangle$ with trivial associated commutator
map $L\times L\rightarrow \Z/6\Z$ and a section
$L\rightarrow \hat{L}$; $\alpha\mapsto e^{\alpha}$ was chosen so that
$e^{\alpha}e^{\beta}=e^{\alpha+\beta}$  and
$\otau(e^{\alpha})=e^{\otau(\alpha)}$. In our case we have
$e^{\alpha}e^{\beta}=e^{\alpha+\beta}$  and
$\otau(e^{\alpha})=e^{\otau(\alpha)}$ in $\widehat{L^{\perp}}$ for
$\alpha,\beta\in L$ by \eqref{eq:prod-in-LD} and \eqref{eq:def-tau}.
Thus for each $1\leq s\leq \ell$, the map
\begin{align*}
\kappa_6 &\mapsto \kz^{6}=\kappa_6,\\
e^{\alpha} &\mapsto (1,\ldots,e^{\alpha},\ldots,1) T
\end{align*}
is an injective group homomorphism of $\hat{L}$ to
$\widehat{L^{\oplus\ell}}$, where $(1,\ldots,e^{\alpha},\ldots,1)$
is the element of $(\widehat{L^\perp})^\ell$ whose $s$-th
component is $e^{\alpha}$ and the other components are $1$.
This injective homomorphism is compatible with
the action of $\otau$.

The embedding of $\hat{L}$ into $\widehat{L^{\oplus\ell}}$ 
gives rise to
an embedding $v\mapsto 1\otimes\cdots\otimes v\otimes\cdots
\otimes 1$ of the vertex operator algebra $V_L$ into
$V_L^{\otimes\ell}\cong V_{L^{\oplus\ell}}$ which maps $V_L$
isomorphically to the $s$-th component of $V_L^{\otimes\ell}$ for
each $1\leq s\leq \ell$. This embedding is again compatible with
the action of $\otau$.

We denote the bilinear form $\varepsilon_0$ on $L$ of
\cite[(4.4)]{DLTYY} by $\varepsilon^{\prime}$ for a while. Thus
$\varepsilon^{\prime}(\alpha,\beta) =
5\langle\otau^2\alpha,\beta\rangle+6\Z$. In \cite{DLTYY}, the
multiplications $a\times b$ in $\hat{L}$ and $a\times_{\otau}b$ in
$\hat{L}_{\otau}$ are related as $a\times
b=\kappa_6^{\varepsilon^{\prime}(\bar{a},\bar{b})}a\times_{\otau}b$.
Since $\kz^{\varepsilon_0(\alpha,\beta)}=
\kappa_{6}^{\varepsilon^{\prime}(\alpha,\beta)}$ for
$\alpha, \beta \in L$ by
\eqref{eq:vep0}, it follows from
\eqref{eq:ve-relation} that the map
$\kappa_6 \mapsto \kz^{6}=\kappa_6$, 
$e^{\alpha} \mapsto e^{\alpha}$ 
for $\alpha\in L$ is an injective group homomorphism of
$\hat{L}_{\otau}$ to the $s$-th component of
${\widehat{L^{\oplus\ell}}_{\htau}}$ for each $1\leq s\leq \ell$.

Now, 
$V_L^{\otimes \ell} \cong V_{L^{\oplus \ell}} \subset \VLCD$. 
Since $\htau=(\otau,\ldots,\otau)$ and since the irreducible
$\otau$-twisted $V_{L}$-modules $V_{L}^{T_{\chi_j}}(\otau)$, 
$j=0,1,2$ of \cite{DLTYY,TY} were
constructed by the same method as in \cite{DL2, L}, 
the above argument shows that the 
action of $V_{L}$ on 
$V_{L}^{T_{\chi_j}}(\otau)$ 
is realized in the action of the $s$-th component
of $V_{L}^{\otimes\ell}$ on the irreducible $\tau$-twisted
$\VLCD$-modules $\VLCDTeta(\htau)$ constructed by \eqref{eq:def-tw} 
below.

We can verify the following properties of the $\Z$-bilinear form $c_2$. 
In fact, it is sufficient to show the assertions for the case $\ell = 1$. 
Note that Lemma \ref{lem:coset-b1b2} implies
\begin{equation*}
L_{\Klein^{\ell}\times\zero} =
\{ \sum_{s=1}^{\ell}(m_1^{(s)}\tilde{\beta}_1^{(s)}
+ 3 m_2^{(s)}\tilde{\beta}_2^{(s)})\,|\,
m_1^{(s)},  m_2^{(s)} \in \Z\}.
\end{equation*}

\begin{lem}\label{lem:pre-dual}
$(1)$
For $\alpha\in (L^{\perp})^{\oplus\ell}$, we have
$c_2(\alpha,\beta)=0$ for all $\beta\in L^{\oplus\ell}$
if and only if $\alpha\in L_{\Klein^{\ell}\times\zero}$.

$(2)$
For $\alpha=\sum_{s=1}^{\ell}
\big( m_1^{(s)}\tilde{\beta}_1^{(s)} + 3m_2^{(s)}\tilde{\beta}_2^{(s)} \big)
\in L_{\Klein^{\ell}\times\zero}$
and $\beta\in L_{(\zero, \gamma)}$ with $\gamma \in \Z_3^\ell$, we have
\begin{equation*}
c_2(\alpha,\beta) = 12\langle(m_1^{(s)})_{s=1}^{\ell},\gamma\rangle_{\Z_3}
+ 36\Z.
\end{equation*}

$(3)$
For $\alpha\in L_{(\lambda,\zero)},
\beta\in L_{(\mu,\zero)}$ with $\lambda, \mu\in \Klein^{\ell}$, we have
\begin{equation*}
c_2(\alpha,\beta) = 18\langle \lambda,\mu\rangle_{\Klein} + 36\Z.
\end{equation*}
\end{lem}

We now follow \cite{L}. The commutator map $C(\al,\be)$ of
\cite{L} is ${\kz}^{c_2(\al,\be)}$ in our notation. 
Let $\h = \C\otimes_\Z \LCD$, 
so that $\h = (\C \otimes_\Z L)^{\oplus \ell}$. 
We extend $\tau$ to an isometry of $\h$ linearly. 
Then $\tau$ is fixed-point-free on $\h$ and $N$ of
\cite{L} is identical with $\LCD$ in our case.

Let $R = \{ \al \in \LCD \,|\, c_2(\al, \be) = 0 \text{ for
all } \be \in \LCD\}$ be the radical of the alternating
$\Z$-bilinear form $c_2$ on $\LCD$, which is identical with
the $R$ of \cite[Section 6]{L}.
Since $C$ is self-orthogonal, Lemma \ref{lem:pre-dual} 
implies the following assertion.

\begin{lem}\label{lem:R-explicit}
The radical $R$ of the alternating $\Z$-bilinear form
$c_2$ on $\LCD$ consists of the elements
\begin{equation*}
\sum_{s=1}^{\ell} \big( m_1^{(s)}\tilde{\beta}_1^{(s)} +
3 m_2^{(s)}\tilde{\beta}_2^{(s)} \big) \in L_{C \times \zero}
\end{equation*}
with
$m_1^{(s)}, m_2^{(s)} \in \Z$ such that
$(m_1^{(s)} + 3\Z)_{s=1}^{\ell} \in D^{\perp}$.
\end{lem}

By Lemma \ref{lem:pre-dual} we also have the following lemma.
Thus we can choose $L_{C \times\zero}$ as the group $A$ of
\cite[Proposition 6.2]{L}.

\begin{lem}\label{lem:A}
$L_{C \times \zero}$ is a subgroup of $\LCD$ which is maximal
subject to the condition that the alternating $\Z$-bilinear form
$c_2$ is trivial on it.
\end{lem}

We shall consider $(1-\htau)\LCD=\cup_{(\lambda,\gamma)
\in C\times D}(1-\htau)L_{(\lambda,\gamma)}$,
which corresponds to the subgroup denoted by $M$ in
\cite[Section 6]{L}.
For $m_1,m_2\in\Z$, we have
\begin{equation*}
(1-\otau)(m_1\tilde{\beta}_1 + m_2\tilde{\beta}_2) =
-m_2\tilde{\beta}_1 + 3(m_1+m_2)\tilde{\beta}_2
\end{equation*}
and hence
$(1-\otau)L^\perp = L_{\Klein \times \zero}$ and
$(1-\otau)L = \Z(6\tilde{\beta}_1)+\Z(6\tilde{\beta}_2)$.
More precisely,
\begin{align*}
(1-\otau)(\dfrac{-\beta_1+\beta_2}{3})
&=2\tilde{\beta}_1-6\tilde{\beta}_2,\\
(1-\otau)(\beta(a))&=\beta(c)+2\tilde{\beta}_1-6\tilde{\beta}_2,\\
(1-\otau)(\beta(b))&=\beta(a)+2\tilde{\beta}_1-6\tilde{\beta}_1
+6\tilde{\beta}_2,\\
(1-\otau)(\beta(c))&=\beta(b)+2\tilde{\beta}_1.
\end{align*}
Then we see from \eqref{eq:coset-lattice} that
\begin{equation}\label{eq:1-tau-Lxi}
(1-\otau)L^{(x,i)} =
\wt_{\Klein}(x) (\beta(\otau^{2}(x)) + 2\tilde{\beta}_1)
+ 2i\tilde{\beta}_1
+ \Z(6\tilde{\beta}_1)+\Z(6\tilde{\beta}_2)
\end{equation}
for $x \in \Klein$ and $i \in \Z_3$,
where $\wt_{\Klein}(x) = 1$ if $x \in \{ a,b,c\}$ and $0$
otherwise.
Thus,
\begin{equation}\label{eq:1-tau-Llg}
(1-\tau)L_{(\lambda,\gamma)} = \sum_{s=1}^{\ell}
\Big( \wt_{\Klein}(\lambda_{s})(\beta(\otau^{2}(\lambda_{s}))
+ 2\tilde{\beta}_1^{(s)})
+ 2\gamma_s \tilde{\beta}_1^{(s)}
+ \Z(6\tilde{\beta}_1^{(s)})+\Z(6\tilde{\beta}_2^{(s)})\Big)
\end{equation}
for $\lambda=(\lambda_{s})_{s=1}^{\ell}\in \Klein^\ell$ and
$\gamma=(\gamma_{s})_{s=1}^{\ell}\in \Z_3^\ell$.
We also note that
\begin{equation}\label{eq:1-tau-dec}
L^{(\otau^2(x),0)} = 
(1-\otau)L^{(x,0)}\cup(1-\otau)L^{(x,1)}\cup(1-\otau)L^{(x,2)}
\end{equation}
is a disjoint union for $x\in\Klein$ by \eqref{eq:coset-lattice}
and \eqref{eq:1-tau-Lxi}.
Thus,
\begin{equation}\label{eq:1-tau-dec3}
L_{C \times \zero} = \bigcup_{\substack{\lm \in C\\
\gm \in \Z_3^\ell}} (1-\tau)L_{(\lm,\gm)}; \quad \text{disjoint}.
\end{equation}

Define a $\Z$-linear map $\varphi\colon L^\perp \rightarrow \Z_3$ by
\begin{equation}\label{eq:varphi-L-perp}
\varphi(m_1\tilde{\beta}_1 + m_2\tilde{\beta}_2) = m_1+3\Z
\end{equation}
for $m_1, m_2 \in \Z$.
We can verify that
$\varphi(\beta(\otau^2(x))+2\tilde{\beta}_1)=0$ if $x \in \{ a,b,c\}$. 
Hence \eqref{eq:1-tau-Lxi} implies the following lemma.

\begin{lem}\label{lem:varphi-sur}
$\varphi((1 - \otau)L^{(x,i)}) = \{2i\}$ for $x\in\Klein$, $i \in \Z_3$.
\end{lem}

We extend $\varphi\colon L^\perp \rightarrow \Z_3$ to a
homomorphism of additive groups $\varphi\colon (L^\perp)^{\oplus \ell}
\rightarrow \Z_3^\ell$ componentwise, so that it maps
the $s$-th component $L^\perp$ to $\Z_3$ by \eqref{eq:varphi-L-perp}.
Set $M_0 = (1-\tau)\LCzero$ and $M = (1-\tau)\LCD$. 
By Lemma \ref{lem:varphi-sur}, we have 
$\varphi((1-\htau)L_{(\lambda,\gamma)})=\{2\gamma\}$ 
for $\lambda \in \Klein^\ell$ and $\gamma \in \Z_3^\ell$. 
Thus the following lemma holds by \eqref{eq:1-tau-dec3} 
and Lemma \ref{lem:R-explicit}.

\begin{lem}\label{lem:A-structure}
The restriction $\varphi|_{L_{C \times \zero}} \colon L_{C \times \zero}
\rightarrow \Z_3^\ell$ of $\varphi$ to $L_{C \times \zero}$
is a surjective homomorphism and its kernel is $M_0$. Moreover,
$\varphi(M) = D$ and $\varphi(R) = D^\perp$. That is, $\varphi$
gives the following surjections.
\begin{equation}
\begin{array}{cccccccc}
\vspace{2mm}
 & M_0 & \subset & M & \subset
& R & \subset & L_{C \times\zero}\\
\vspace{2mm} \varphi \colon \quad & \downarrow & & \downarrow & &
\downarrow & & \downarrow \\
& \{\zero\} & \subset & D & \subset & D^\perp & \subset &
\Z_3^\ell.
\end{array}
\end{equation}
\end{lem}

Since $6\tilde{\beta}_1^{(s)} = 3\beta_1^{(s)}$ and
$6\tilde{\beta}_2^{(s)} = \beta_1^{(s)} - \beta_2^{(s)}$,
$M_0$ contains $\be_1^{(s)} - \be_2^{(s)}$,
$\be_2^{(s)} - \be_0^{(s)}$ and $3\be_i^{(s)}$, $i=0,1,2$
by \eqref{eq:1-tau-Llg}.
Let $\gm = (\gm_1,\ldots,\gm_\ell) \in \Z_3^\ell$.
Then the inverse image of $\{ 2\gm \}$ under
$\varphi|_{L_{C \times \zero}}$ is
$\sum_{s=1}^\ell \gm_s \beta_i^{(s)} + M_0$.
By Lemma \ref{lem:A-structure},
$\varphi$ induces an isomorphism
$L_{C \times\zero}/M_0 \cong \Z_3^\ell$. Taking the inverse image
of $D$, $D^\perp$ and $\Z_3^\ell$, respectively, 
we have the following coset decompositions.
\begin{gather}
M = \bigcup_{\gm \in D} (\gm_1\be_i^{(1)} + \cdots +
\gm_\ell\be_i^{(\ell)} + M_0),\label{eq:coset-0}\\
R = \bigcup_{\gm \in D^\perp} (\gm_1\be_i^{(1)}
+ \cdots + \gm_\ell\be_i^{(\ell)} + M_0),\label{eq:coset-1}\\
L_{C \times \zero} = \bigcup_{\gm \in \Z_3^\ell}
(\gm_1\be_i^{(1)} + \cdots + \gm_\ell\be_i^{(\ell)} + M_0).
\label{eq:coset-2}
\end{gather}

Recall that $\hat{Q}_{\tau}$ denotes the inverse image of $Q$ under 
the homomorphism $\hLCDtau \xrightarrow{-} \LCD$ for a subset 
$Q$ of $\LCD$. Lemma \ref{lem:A} implies that the inverse
image $\hLCzerotau$ of $\LCzero$ is isomorphic to $\LCzero \times
\la \kz \ra$, which is a maximal abelian subgroup of
$\hLCDtau$. The inverse image $\hat{R}_\tau$ of $R$ is the center of
the group $\hLCDtau$.

A central subgroup $K$ defined in \cite[Remark 4.2]{DL2} is crucial
for the construction of a certain class of irreducible $\hLCDtau$-modules
(see also \cite[Section 7.4]{FLM}, \cite[Section 6]{L}). Let
$K = \{a\htau(a)^{-1} \,|\, a \in \hLCDtau \}$. Then $K$ is a subgroup of
the center  $\hat{R}_\tau$ of $\hLCDtau$ and $K \cap \la \kz\ra
= 1$ \cite[Remark 4.2]{DL2}. Indeed,
$\overline{a\htau(a)^{-1}} = \overline{a} - \tau(\overline{a}) \in M$.
If $a\htau(a)^{-1} \in \la \kz \ra$, then $\ova =
\tau(\ova)$ and so $\ova =0$. Hence $a \in \la \kz \ra$ and
$a\htau(a)^{-1} = 1$. Thus $K \cap \la \kz\ra = 1$. 
Since $K$ lies in $\hat{R}_\tau$, 
$b\htau(b)^{-1}$ commutes with $\htau(a)^{-1}$ 
for $a, b \in \hLCDtau$ and 
\begin{equation}\label{eq:prod-in-K}
a\htau(a)^{-1}b\htau(b)^{-1} =
ab\htau(b)^{-1}\htau(a)^{-1} = ab\htau(ab)^{-1}.
\end{equation}
Thus $K$ is a group. Now the inverse image $\hMtau$
of $M$ in $\hLCDtau$ is $K \times \la \kz \ra \cong M \times
\la \kz \ra$. Clearly, $K$ is $\htau$-invariant. Moreover, $K$
is $\htheta$-invariant since $\htheta$ commutes with $\htau$.

We shall construct an irreducible $\hLCDtau$-module $T_\psi$ as in
\cite[Proposition 6.2]{L}. Since $\hMtau = K \times \la \kz
\ra$, there is a unique group homomorphism $\rho \colon \hMtau
\rightarrow \C^\times$ such that $\rho(\kz) = \zeta_{36}$ and
$\rho(a) = 1$ for $a \in K$.
Note that $(1+ \tau+\tau^2)\al = 0$ for $\al \in \LCD$. Thus $\rho$
is the homomorphism denoted by $\tau$ in \cite[Proposition 6.1]{L}.
Let $\chi \colon \hRtau \rightarrow \C^\times$ be a homomorphism
extending $\rho$ and $\psi \colon \hLCzerotau \rightarrow \C^\times$ be a
homomorphism extending $\chi$. Then $\psi(\kz) =\zeta_{36}$ and
$\psi$ is $1$ on $K$. Such an extension $\psi$ exists, since in the
central extension
\begin{equation*}
1 \rightarrow \la \kz\ra \rightarrow \hLCDtau/K \rightarrow
\LCD/M \rightarrow 1
\end{equation*}
with associated commutator map $\overline{c}_2$ defined by
$\overline{c}_2(\al+M, \beta+M) =
c_2(\al ,\beta)$, the subgroup $\hLCzerotau/K$ splits by Lemma
\ref{lem:A}. That is, $\hLCzerotau/K \cong (\LCzero/M) \times
\la\kz\ra$ and $\hRtau/K \cong (R/M) \times \la \kz \ra$.
Let $\C_\psi$ be a one dimensional $\hLCzerotau$-module with
character $\psi$ and $T_\psi = \C[\hLCDtau]
\otimes_{\C[\hLCzerotau]} \C_\psi$ be the $\hLCDtau$-module induced
from $\C_\psi$.

We need to know $\psi$ and $T_\psi$ in detail. For this purpose, set
$K_0 = \{ a\htau(a)^{-1}\,|\, a \in \hLCzerotau\}$. Then $K_0$ is a
subgroup of $K$ with $\hat{M}_{0,\tau} = K_0 \times \la\kz\ra$,
where $\hat{M}_{0,\tau}$ denotes the inverse image of $M_0$ in
$\hLCDtau$. 
Moreover, $K_0$ is $\htheta$- and
$\htau$-invariant. We shall describe the group
$\hLCDtau/K_0$ explicitly.

We can verify that
$\vep_2(\alpha,\otau(\alpha)) = \vep_2(\alpha,\alpha)$ and
$e^{\alpha}\otau(e^{\alpha})^{-1}=e^{(1-\otau)\alpha}$ in
${\widehat{L^{\perp}}}_{\otau}$ for $\alpha\in L^{\perp}$
by \eqref{eq:vep2}, \eqref{eq:mul2} and \eqref{eq:def-tau}.
Hence
\begin{equation}\label{eq:1-tau:general}
e^\beta\htau(e^\beta)^{-1} = e^{(1-\tau)\beta} \quad \text{in\ }
\widehat{(L^{\perp})^{\oplus\ell}}_{\htau}
\end{equation}
for $\beta \in (L^{\perp})^{\oplus\ell}$.
In the case of $\beta = -\beta_1^{(s)} +
\beta_2^{(s)}$, we have
\begin{equation}\label{eq:1-tau:special}
e^{-\beta_1^{(s)} + \beta_2^{(s)}} \htau(e^{-\beta_1^{(s)} +
\beta_2^{(s)}})^{-1} = e^{3\beta_2^{(s)}} \quad \text{in\ }
\hLCDtau.
\end{equation}

For $\gm = (\gm_1,\ldots,\gm_\ell) \in D$, set
\begin{equation}\label{eq:a-gamma}
a(\gm) = \sum_{s=1}^\ell j_s\frac{-\beta_1^{(s)}+\beta_2^{(s)}}{3}
\in \LzeroD,
\end{equation}
where $j_s = 0,1,2$ such that $\gm_s = j_s + 3\Z$.
These $a(\gm), \gm \in D$ form a complete set of coset
representatives of $\LCzero$ in $\LCD$,
and so
\begin{equation}\label{eq:cosets-in-hL}
\hLCDtau = \bigcup_{\gm \in D} e^{a(\gm)} \hLCzerotau.
\end{equation}
Then using \eqref{eq:prod-in-K} we see that
\begin{equation}\label{eq:coset-in-K-1}
K = \bigcup_{\gm \in D} e^{a(\gm)}\htau(e^{a(\gm)})^{-1} K_0.
\end{equation}
Moreover, it follows from \eqref{eq:1-tau:general} that
\begin{equation}\label{eq:1-tau:a-gamma}
e^{a(\gm)}\htau(e^{a(\gm)})^{-1} = e^{\sum_{s=1}^\ell
j_s\beta_2^{(s)}} \quad \text{in\ } \hLCDtau.
\end{equation}

Now, using \eqref{eq:vep2} and \eqref{eq:mul2} we can
verify that
\begin{equation}\label{eq:power-in-twisted}
(e^{\beta_i})^m = \kappa_{3}^{m(m-1)}e^{m\beta_i}
\quad \text{in\ } \widehat{L^\perp}_{\otau}
\end{equation}
for $m \in \Z$, $i=0,1,2$.

By \eqref{eq:def-tau}, \eqref{eq:1-tau:special} and 
\eqref{eq:power-in-twisted}, we have the following lemma.

\begin{lem}\label{lem:mod-K_0}
The following assertions hold in $\hLCDtau$ for $1 \le s \le \ell$.

$(1)$ $e^{\beta_1^{(s)}} \equiv e^{\beta_2^{(s)}} \equiv
e^{\beta_0^{(s)}} \pmod{K_0}$.

$(2)$ $(\kappa_3e^{\beta_i^{(s)}})^3 = e^{3\beta_i^{(s)}} \in K_0$,
$i=0,1,2$.

$(3)$ $(\kappa_3e^{\beta_i^{(s)}})^{-1} =
\kappa_3e^{-\beta_i^{(s)}}$, $i=0,1,2$.
\end{lem}

By \eqref{eq:power-in-twisted},
$(\kappa_3e^{\beta_i^{(s)}})^m = \kappa_3^{m^2}e^{m\beta_i^{(s)}}$
in $\hLCDtau$ for any integer $m$.
Now, let $m_s\in\Z, 1\leq s\leq \ell$. Then
$e^{\sum_{s=1}^\ell m_s\beta_i^{(s)}} = e^{m_1\beta_i^{(1)}}
\cdots e^{m_\ell\beta_i^{(\ell)}}$ in $\hLCDtau$,
since $\vep_2(\beta_i^{(s)}, \beta_i^{(t)}) = 0$ if $s \ne t$.
Thus
\begin{equation}\label{eq:general-product}
(\kappa_3e^{\beta_i^{(1)}})^{m_1} \cdots
(\kappa_3e^{\beta_i^{(\ell)}})^{m_\ell} 
= \kappa_3^{\sum_{s=1}^\ell m_s^2}
e^{\sum_{s=1}^\ell m_s\beta_i^{(s)}} \quad \text{in\ } \hLCDtau
\end{equation}
for any $(m_1,\ldots,m_\ell) \in \Z^\ell$.
The above lemma implies that  $(\kappa_3{\beta_i^{(s)}})^{m}K_0$ and
$e^{m\beta_i^{(s)}}K_0$ depend only on $m\pmod{3\Z}$. Hence
(\ref{eq:general-product}) is reduced to
\begin{equation}
(\kappa_3e^{\beta_i^{(1)}})^{\gamma_1} \cdots
(\kappa_3e^{\beta_i^{(\ell)}})^{\gamma_\ell}K_0 =
\kappa_3^{\langle\gamma,\gamma\rangle_{\Z_3}}
e^{\sum_{s=1}^\ell \gamma_s\beta_i^{(s)}}K_0
\end{equation}
modulo $K_0$ for $\gamma=(\gamma_1,\ldots,\gamma_{\ell})\in\Z_3^{\ell}$.
If $\gamma\in D$, then $\langle\gamma,\gamma\rangle_{\Z_3}=0$ since
$D$ is self-orthogonal. Therefore,
\eqref{eq:coset-in-K-1} and
\eqref{eq:1-tau:a-gamma} give that
\begin{equation}\label{eq:coset-in-K-2}
K = \bigcup_{\gm \in D} (\kappa_3e^{\beta_i^{(1)}})^{\gm_1} \cdots
(\kappa_3e^{\beta_i^{(\ell)}})^{\gm_\ell} K_0.
\end{equation}

Motivated by the above result, we set
\begin{align}
K_1 &= \bigcup_{\gm \in D^\perp} (\kappa_3e^{\beta_i^{(1)}})^{\gm_1}
\cdots (\kappa_3e^{\beta_i^{(\ell)}})^{\gm_\ell} K_0,
\label{eq:def-K1}\\
K_2 &= \bigcup_{\gm \in \Z_3^\ell}
(\kappa_3e^{\beta_i^{(1)}})^{\gm_1} \cdots
(\kappa_3e^{\beta_i^{(\ell)}})^{\gm_\ell} K_0 \label{eq:def-K2}
\end{align}
with $\gm = (\gm_1,\ldots,\gm_\ell)$. 
Then the following lemma holds.

\begin{lem}\label{lem:K1K2}
$(1)$ $K_2$ is a subgroup of $\hLCzerotau$ such that $K_2 \cap
\la\kz\ra = 1$ and $\hLCzerotau = K_2 \times \la\kz\ra$.
Moreover,
\begin{equation*}
K_2/K_0 = \la \kappa_3e^{\beta_i^{(1)}} K_0/K_0\ra \times \cdots
\times \la \kappa_3e^{\beta_i^{(\ell)}} K_0/K_0\ra,
\end{equation*}
which is isomorphic to $\LCzero/M_0 \cong \Z_3^\ell$.

$(2)$ $K_1$ is a subgroup of $K_2$ such that $\hat{R}_\tau = K_1
\times \la\kz\ra$. Moreover, $K_1/K_0$ is isomorphic to $R/M_0
\cong D^\perp$.
\end{lem}

Let $\psi \colon \hLCzerotau \rightarrow \C^\times$ be a homomorphism of
abelian groups such that $\psi(\kz) =\zeta_{36}$ and $\psi(a) = 1$
for $a \in K_0$. Then $\psi(\kappa_3e^{\beta_i^{(s)}}) =
\zeta_3^{\eta_s}$, $1 \le s \le \ell$ for some $\eta =
(\eta_1,\ldots,\eta_\ell) \in \Z_3^\ell$ by Lemma \ref{lem:K1K2}.
We denote such a homomorphism $\psi$ by $\psi_\eta$. In fact, $\eta
\mapsto \psi_\eta$ is an isomorphism of the additive group
$\Z_3^\ell$ onto the multiplicative group of all homomorphisms
$\psi \colon \hLCzerotau \rightarrow \C^\times$ with $\psi(\kz) =
\zeta_{36}$ and $\psi(a) = 1$ for $a \in K_0$. The homomorphism
$\psi_\eta$ is determined by the three conditions (i)
$\psi_\eta(\kz) =\zeta_{36}$, (ii) $\psi_\eta$ is $1$ on $K_0$, and
(iii) $\psi_\eta(\kappa_3e^{\beta_i^{(s)}}) = \zeta_3^{\eta_s}$. 

\begin{rmk}
The conditions {\rm (i)}, {\rm (ii)}, and {\rm (iii)} for
$\psi_\eta$ are consistent with the conditions for $\chi_j$ in
\cite[Section 4]{DLTYY}.
\end{rmk}

As before, 
let $\Ceta$ be a one dimensional
$\hLCzerotau$-module affording the character $\psi_\eta$ and
$\Tpe = \C[\hLCDtau] \otimes_{\C[\hLCzerotau]}
\Ceta$ be the $\hLCDtau$-module induced from
$\Ceta$.  It follows from \eqref{eq:cosets-in-hL} that $\{
e^{a(\gm)} \otimes 1_{\eta}\,|\, \gm \in D\}$ is a basis of
$\Tpe$, where $1_{\eta}$ denotes a fixed nonzero vector in
$\Ceta$.
For $b \in \hLCzerotau$, we have 
$be^{a(\gm)} = \kz^{c_2(\overline{b},a(\gm))} e^{a(\gm)}b$ 
and the action of $b$ on
$e^{a(\gm)} \otimes 1_{\eta}$ is
\begin{equation*}
b\cdot (e^{a(\gm)} \otimes 1_{\eta}) =
\zeta_{36}^{c_2(\overline{b},a(\gm))} \psi_\eta(b)
(e^{a(\gm)}\otimes 1_{\eta}).
\end{equation*}

For $\dl \in D$, we have 
$e^{a(\dl)} e^{a(\gm)} \in e^{a(\dl+\gm)}\hLCzerotau$ by 
\eqref{eq:mul2} since $a(\dl)+a(\gm) \equiv
a(\dl+\gm) \pmod{L^{\oplus \ell}}$. 
Then $\Tpe$ is an irreducible $\hLCDtau$-module 
and the following lemma holds.

\begin{lem}\label{lem:irreducible-module-T}
$(1)$ $\kz$ and $K_0$ act on $\Tpe$ as $\zeta_{36}$ and
$1$, respectively. Moreover, $K$ $($resp. $K_1$$)$ acts on
$\Tpe$ as $1$ if and only if $\eta \in D^\perp$ $($resp.
$\eta \in D$$)$.

$(2)$ For $\eta, \eta' \in \Z_3^\ell$, the $\hLCDtau$-modules
$\Tpe$ and $T_{\psi_{\eta'}}$ are equivalent if and only if
$\eta \equiv \eta' \pmod{D}$, which is also equivalent to the
condition that $\psi_\eta$ and $\psi_{\eta'}$ agree on $K_1$.

$(3)$ The action of $\kappa_3e^{\pm\beta_i^{(s)}}$ on
$e^{a(\gm)}\otimes 1_{\eta}$ is such that
\begin{equation*}
\kappa_3e^{\pm\beta_i^{(s)}} \cdot (e^{a(\gm)} \otimes 1_{\eta}) =
\zeta_3^{\pm(\eta_s - \gm_s)} e^{a(\gm)} \otimes 1_{\eta}.
\end{equation*}
That is, $\C e^{a(\gm)} \otimes 1_{\eta}$ is a one dimensional
$\hLCzerotau$-module with character $\psi_{\eta-\gm}$.
\end{lem}

By the above lemma, $e^{\pm \beta_i^{(s)}}$ acts 
on $e^{a(\gm)} \otimes 1_{\eta} \in \Tpe$ as
\begin{equation}\label{eq:action-on-basis}
e^{\pm \beta_i^{(s)}} \cdot (e^{a(\gm)} \otimes 1_{\eta}) =
\zeta_3^{-1\pm(\eta_s - \gm_s)}e^{a(\gm)} \otimes 1_{\eta}.
\end{equation}

\begin{rmk}
$\Tpe$, $\eta
\in D^\perp$ are exactly the modules $T$ of
\cite[Proposition 6.2]{L} in our case.
\end{rmk}

Recall that $\h = \C \otimes_{\Z} L_{C \times D} 
= (\C \otimes_\Z L)^{\oplus \ell}$.
As before, we use $\al^{(s)}$ to denote the element $\al \in \C
\otimes_\Z L$ in the $s$-th entry of $(\C \otimes_\Z L)^{\oplus \ell}$.
Let
\begin{equation*}
h_1^{(s)} = \frac{1}{3}(\be_1^{(s)} + \zeta_3^2\be_2^{(s)} +
\zeta_3\be_0^{(s)}), \quad h_2^{(s)} = \frac{1}{3}(\be_1^{(s)} +
\zeta_3\be_2^{(s)} + \zeta_3^2\be_0^{(s)}).
\end{equation*}
Then $\tau h_j^{(s)} = \zeta_3^j h_j^{(s)}$, $\la h_j^{(s)},h_j^{(t)}
\ra = 0$, and $\la h_1^{(s)},h_2^{(t)} \ra = 2\delta_{s,t}$. Set
\begin{equation*}
\h_{(n)} = \{ \al \in \h \,|\, \tau\al = \zeta_3^n \al \}
\end{equation*}
for $n \in \Z$. The index $n$ of $\h_{(n)}$ is considered to be
modulo $3$. Then $\h_{(0)} = 0$ and $\h = \h_{(1)} \oplus \h_{(2)}$
with $\h_{(n)} = \C h_n^{(1)} \oplus \cdots \oplus \C h_n^{(\ell)}$,
$n=1,2$. If $\al \in \h$, we write $\al_{(n)}$ for the component of
$\al$ in $\h_{(n)}$. In this notation we have $(\be_i^{(s)})_{(1)} =
\zeta_3^{i-1}h_1^{(s)}$ and $(\be_i^{(s)})_{(2)} =
\zeta_3^{2(i-1)}h_2^{(s)}$, $i=0,1,2$.

The $\tau$-twisted affine Lie algebra $\hat{\h}[\tau]$ is defined to
be
\begin{equation*}
\hat{\h}[\tau] = \Big( \bigoplus_{n \in \Z} \h_{(n)} \otimes t^{n/3}
\Big) \oplus \C c
\end{equation*}
with the bracket
\begin{equation*}
[x \otimes t^m, y \otimes t^n] = m \la x,y \ra \dl_{m+n,0}c
\end{equation*}
for $x \in \h_{(3m)}$, $y \in \h_{(3n)}$, $m$, $n \in (1/3)\Z$ and
$[c, \hat{\h}[\tau]] = 0$. The isometry $\tau$ acts on
$\hat{\h}[\tau]$ by $\tau( x \otimes t^{n/3}) = \zeta_3^{n} x \otimes
t^{n/3}$ and $\tau(c) = c$. 
Set
\begin{equation*}
\hat{\h}[\tau]^+ = \bigoplus_{n>0} \h_{(n)} \otimes t^{n/3}, \quad
\hat{\h}[\tau]^- = \bigoplus_{n<0} \h_{(n)} \otimes t^{n/3}, \quad
\hat{\h}[\tau]^0 = \C c
\end{equation*}
and consider the $\hat{\h}[\tau]$-module
\begin{equation*}
S[\tau] = U(\hat{\h}[\tau]) \otimes_{U(\hat{\h}[\tau]^+ \oplus
\hat{\h}[\tau]^0)} \C
\end{equation*}
induced from the $\hat{\h}[\tau]^+ \oplus \hat{\h}[\tau]^0$-module
$\C$, where $\hat{\h}[\tau]^+$ acts as $0$ and $\hat{\h}[\tau]^0$
acts as $1$ on $\C$. 
The weight gradation on $S[\tau]$ is given by 
$\wt(x\otimes t^n) = -n$ and $\wt(1) = \ell/9$ 
for $n \in (1/3)\Z$ and $x \in \h_{(3n)}$ 
\cite[(4.6), (4.10)]{DL2}. 
For $\al \in \h$ and $n
\in (1/3)\Z$, we write $\al(n)$ for the operator on $S[\tau]$
induced by the action of $\al_{(3n)} \otimes t^n$. 
The weight of the operator $h_i^{(s)}(i/3+n)$ is $-i/3-n$.
The group $\tgrs$ acts as 
\begin{equation*}
(\otau^{j_1},\ldots,\otau^{j_{\ell}})(h_i^{(s)}(\frac{i}{3}+n))
=\zeta_3^{j_{s}i} h_i^{(s)}(\frac{i}{3}+n).
\end{equation*}

Set
\begin{equation}\label{eq:def-tw}
\VLCDTeta(\htau) = S[\tau]\otimes \Tpe
\end{equation}
for $\eta \in D^\perp$. 
By \cite[Theorem 7.1]{DL2} and \cite[Proposition 6.2]{L}, we can define 
a $\tau$-twisted vertex operator $Y^\tau(\,\cdot\,,x)\colon V_{\LCD} 
\rightarrow \End(\VLCDTeta(\htau))\{ x\}$ so that 
$(\VLCDTeta(\htau), Y^\tau)$, $\eta \in D^\perp$ is an irreducible 
$\tau$-twisted $V_{\LCD}$-module. 
The weight of any element in $\Tpe$
is defined to be $0$. Hence the weight of elements in $\VLCDTeta(\htau)$ is
given by $\wt(u\otimes v) = \wt(u)$ for $u \in S[\htau]$ and $v \in
\Tpe$.

We define an action of $\tgrs$ on $\C e^{a(\gm)} \otimes 1_{\eta}$ by
\begin{equation*}
(\otau^{j_1},\ldots,\otau^{j_{\ell}})(e^{a(\gamma)}\otimes 1_{\eta})
= \zeta_3^{2\langle 
(j_s)_{s=1}^{\ell},\eta-\gamma\rangle_{\Z_3}}e^{a(\gamma)}\otimes 1_{\eta}
\end{equation*}
and extend to $\Tpe=\oplus_{\gamma\in D}\C e^{a(\gm)} \otimes 1_{\eta}$ 
by linearity.
Note that
Lemma \ref{lem:coset} implies
$\htau(e^{a(\gamma)}\otimes 1_{\eta})
=\zeta_3^{2\wt_{\Z_3}(\eta)}e^{a(\gamma)}\otimes 1_{\eta}$ 
for $\gamma\in D$.
Thus $\htau$ acts on $\Tpe$ as a scalar $\zeta_3^{2\wt_{\Z_3}(\eta)}$,
which depends only on the coset $\eta+D\in D^{\perp}/D$. 
The group $\tgrs$ acts on the vector space $\VLCDTeta(\htau)$ by 
\begin{equation}\label{eq:Hlacts-1}
g(u \otimes v) = g(u) \otimes g(v)
\end{equation}
for $g\in \tgrs$, $u \in S[\tau]$ and $v \in \Tpe$.
Then, $\htau(Y^{\htau}(u,x)w)=Y^{\htau}(\htau u,x)\htau w$
for $u \in V_{L_{C\times D}}$ and $w\in \VLCDTeta(\tau)$ by 
\cite[Section 4]{DL}.

We have discussed only irreducible $\tau$-twisted $\VLCD$-modules
so far. Now, we deal with irreducible $\tau^2$-twisted ones.
Actually, we can construct $|D^{\perp}/D|$ inequivalent
irreducible $\htau^2$-twisted $\VLCD$-modules
$(\VLCDTeta(\htau^2), Y^{\htau^2})$, $\eta \in D^\perp\pmod{D}$
similarly. Indeed, replace $\tau$ with $\tau^2$ in the above
argument and proceed in the same way. 
We can construct a class of irreducible $\hat{L}_{C
\times D, \tau^2}$-modules $\Tpe'$, $\eta \in D^\perp$.
Let $h_1'^{(s)} = h_2^{(s)}$ and $h_2'^{(s)} = h_1^{(s)}$, $1 \le
s \le \ell$. Set $\h'_{(n)} = \{ \al \in \h \,|\, \tau^2 \al =
\zeta_3^n \al \}$ for $n \in \Z$ (see \cite[Section 4.3]{DLTYY}).
Take $h_1'^{(s)}$ and $h_2'^{(s)}$ instead of $h_1^{(s)}$ and
$h_2^{(s)}$, respectively and consider $S[\tau^2]$. 
Then
\begin{equation*}
\VLCDTeta(\tau^2) = S[\tau^2] \otimes \Tpe'.
\end{equation*}

We define an action of $\tgrs$ on $\C e^{a(\gm)} \otimes 1_{\eta}$ by
\begin{equation*}
(\otau^{j_1},\ldots,\otau^{j_{\ell}})(e^{a(\gamma)}\otimes 1_{\eta})
= \zeta_3^{\langle (j_s)_{s=1}^{\ell},\eta-\gamma\rangle_{\Z_3}}
e^{a(\gamma)}\otimes 1_{\eta}
\end{equation*}
and extend to $\Tpe'=\oplus_{\gamma\in D}\C e^{a(\gm)} \otimes 1_{\eta}$ 
by linearity. 
Thus $\tau^2(v) = \zeta_3^{2\wt_{\Z_3}(\eta)}v$ for $v \in \Tpe'$. 
Now, $\tgrs$ acts on the vector space
$\VLCDTeta(\tau^2)$ by 
\begin{equation}\label{eq:Hlacts-2}
g(u \otimes v) =g(u) \otimes g(v)
\end{equation}
for $g\in \tgrs, u \in S[\tau^2]$ and $v \in \Tpe'$. 
We have $\tau(Y^{\tau^2}(u,x)w) = Y^{\tau^2}(\tau u,x) \tau w$ for 
$u \in V_{\LCD}$ and $w \in \VLCDTeta(\tau^2)$.

Since $V_{L_{C\times D}}$ is rational and $C_2$-cofinite, the
number of irreducible $\htau^i$-twisted $V_{\LCD}$-modules is
bounded above by the number of $\htau$-stable irreducible
$V_{L_{C\times D}}$-modules by \cite[Theorem 10.2]{DLM2} for each
$i=1,2$. Now, $\{V_{L_{C\times(\eta+D)}}\, |\, \eta \in
D^{\perp}\pmod{D}\}$ is the set of all $\htau$-stable irreducible
$V_{L_{C\times D}}$-modules up to isomorphism. Hence we have the
following theorem.

\begin{thm}
For $i=1,2$, there are exactly $|D^\perp/D|$ inequivalent
irreducible $\htau^i$-twisted $V_{\LCD}$-modules. They are
represented by $(\VLCDTeta(\htau^i), Y^{\htau^i})$, $\eta \in
D^\perp \pmod{D}$.
\end{thm}

The map $\widehat{L^{\oplus\ell}}_{\tau} \rightarrow 
\widehat{L}_{C\times \zero,\tau}/K_0;\ a \mapsto aK_0$ is surjective 
by \eqref{eq:coset-2} and \eqref{eq:1-tau:general}.
For $\alpha \in L_{\Klein\times\zero}$, note that 
$\alpha \in L$ if $(1-\otau)\alpha \in L$. 
Then 
\begin{equation*}
\{a\times_{\otau}\otau(a)^{-1}\ |\ a\in \hat{L}_{\Klein\times\zero,\otau}\}
\cap \hat{L}_{\otau}=
\{a\times_{\otau}\otau(a)^{-1}\ |\ a\in \hat{L}_{\otau}\}
\end{equation*}
and the following lemma holds.

\begin{lem}\label{lem:L0-Lc}
The map $\widehat{L^{\oplus\ell}}_{\tau} \rightarrow
\widehat{L}_{C\times \zero,\tau}$; $a \mapsto a$ 
induces an isomorphism
$\widehat{L^{\oplus\ell}}_{\tau}/\{a\times_{\tau}\tau(a)^{-1}\ |\
a\in \widehat{L^{\oplus\ell}}_{\tau}\}\cong
\widehat{L}_{C\times \zero,\tau}/K_0$.
\end{lem}

For $i=1,2$ and $\vep \in \Z_3$, set
\begin{equation*}
\vt{L_{C\times D}}{\eta}{\htau^i}{\varepsilon} =\{u\in
V_{L_{C\times D}}^{T,\eta}(\tau^i)\, |\, \tau^i u =
\zeta_3^{\varepsilon} u\}.
\end{equation*}
These are irreducible $\VLCDtau$-modules.

In the case where $\ell=1$ with $C = \{\zero \}$ and $D = \{\zero \}$,
$V_{L_{C \times D}}^{T,\eta}(\tau^i)[\varepsilon]$ reduces to
$V_L^{T,j}(\otau^i)[\varepsilon]$, $j=0,1,2$.
The relation between our $\vt{L}{j}{\otau^i}{\varepsilon}$ and
$V_L^{T_{\chi_j}}(\tau)(\varepsilon)$,
$V_L^{T_{\chi^{\prime}_j}}(\tau^2)(\varepsilon)$ in \cite{DLTYY,TY}
is as
follows (see \cite[(1-1)]{TY} also).
\begin{align*}
\vt{L}{0}{\otau}{\varepsilon}
&=V_L^{T_{\chi_0}}(\tau)(\varepsilon),\\
\vt{L}{0}{\otau^2}{\varepsilon}
&=V_L^{T_{\chi^{\prime}_0}}(\tau^2)(\varepsilon),\\
\vt{L}{j}{\otau}{\varepsilon}
&=V_L^{T_{\chi_j}}(\tau)(\varepsilon+1), \quad j=1,2,\\
\vt{L}{j}{\otau^2}{\varepsilon}
&=V_L^{T_{\chi^{\prime}_j}}(\tau^2)(\varepsilon+1),
\quad j=1,2
\end{align*}
for $\varepsilon\in\Z_3$. Recall that the action of $\tau$ on
$T_{\chi_j}$ and $T_{\chi^{\prime}_j}$ was defined to be $1$ in
\cite{DLTYY,TY},
while $\tau$ acts on $\Tpe$ (resp.
$\Tpe'$) as $\zeta_3^{2\wt_{\Z_3}(\eta)}$ (resp.$\zeta_3^{\wt_{\Z_3}(\eta)}$). The new notation is suitable for
the description of the fusion rules in later sections.

Since
$V_{L^{\oplus\ell}}=V_{L}^{\otimes\ell}$ is a vertex operator subalgebra of
$V_{L_{C\times D}}$,
$V_{L_{C\times D}}^{T,\eta}(\tau)$ is a
$\tau$-twisted $V_{L^{\oplus\ell}}$-module
and for each $\gamma\in D$,
$S[\tau] \otimes (e^{a(\gm)} \otimes 1_{\eta})$
is a $\tau$-twisted $V_{L^{\oplus\ell}}$-submodule
of $V_{L_{C\times D}}^{T,\eta}(\tau)$.
By Lemma \ref{lem:irreducible-module-T} (3) and Lemma \ref{lem:L0-Lc},
the $\hat{L}_{C\times\zero,\htau}$-module
$\C e^{a(\gm)} \otimes 1_{\eta}$ is isomorphic to
$\C_{\psi_{\eta-\gamma}}$ as
$\widehat{L^{\oplus\ell}}_{\htau}$-modules.
This implies that
$S[\tau] \otimes (e^{a(\gm)} \otimes 1_{\eta})\cong
V_{L^{\oplus\ell}}^{T,\eta-\gamma}(\tau)$
as $\htau$-twisted $V_{L^{\oplus\ell}}$-modules.
Thus,
\begin{equation}\label{eq:twist-D-L}
V_{L_{C\times D}}^{T,\eta}(\tau) 
\cong\bigoplus_{\gamma\in D}V_{L^{\oplus\ell}}^{T,\eta-\gamma}(\tau)
\end{equation}
as $\htau$-twisted $V_{L^{\oplus\ell}}$-modules.
For $\rho = (\rho_1,\ldots,\rho_\ell) \in\Z_3^{\ell}$, 
$\tgrs$ acts on $V_{L^{\oplus\ell}}^{T,\rho}(\tau)$.
Note that
$\tgrs$ is an automorphism group of $V_{L^{\oplus\ell}}$
and $g(Y^{\htau}(u,x)w) = Y^{\htau}(gu,x)gw$
for $g\in\tgrs$, $u\in V_{L^{\oplus\ell}}$ and 
$w \in V_{L^{\oplus\ell}}^{T,\rho}(\tau)$ by the definition of 
$(V_{L^{\oplus\ell}}^{T,\rho}(\tau),Y^\tau)$.
Thus, $V_{L^{\oplus\ell}}^{T,\rho}(\tau)\circ g\cong 
g^{-1}(V_{L^{\oplus\ell}}^{T,\rho}(\tau))=V_{L^{\oplus\ell}}^{T,\rho}(\tau)$
for $g\in\tgrs$.
Note that $(V_{L^{\oplus\ell}})^{\tgrs}=(V_L^{\otau})^{\otimes\ell}$. 
We have the following  decomposition of 
$V_{L^{\oplus\ell}}^{T,\rho}(\htau)$ into a direct sum of 
irreducible $(V_L^{\otau})^{\otimes\ell}$-modules.
\begin{align}
V_{L^{\oplus\ell}}^{T,\rho}(\htau)&\cong\bigoplus_{(\vep_1, \ldots, \vep_\ell)
\in \Z_3^\ell} \vt{L}{\rho_1}{\htau}{\vep_1} \otimes \cdots
\otimes \vt{L}{\rho_\ell}{\htau}{\vep_\ell}.
\label{eq:deco-T}
\end{align}
It follows from \cite[Theorem 2]{MT}
that $\vt{L}{\rho_1}{\htau}{\vep_1} \otimes \cdots
\otimes \vt{L}{\rho_\ell}{\htau}{\vep_\ell}, (\vep_1, \ldots, \vep_\ell)
\in \Z_3^\ell$ in \eqref{eq:deco-T} are all inequivalent 
irreducible $(V_L^{\otau})^{\otimes\ell}$-modules.

The corresponding results for $\htau^2$-twisted 
$V_{L_{C\times D}}$-modules can be verified by a similar 
argument as above. Thus we have
obtained the following theorem.

\begin{thm}\label{thm:twist-L-st}
For $i=1,2$ and $\eta = (\eta_1,\ldots,\eta_\ell) \in D^\perp
\pmod{D}$, the irreducible $\htau^i$-twisted $\VLCD$-module
$(\VLCDTeta(\htau^i), Y^{\htau^i})$ is decomposed into a direct
sum of irreducible $(V_{L}^{\otau})^{\otimes {\ell}}$-modules as
follows.
\begin{equation*}
V_{L_{C\times D}}^{T,\eta}(\htau^i) \cong \bigoplus_{(\gm_1,\ldots,\gm_\ell) \in D}
\bigoplus_{(\vep_1, \ldots, \vep_\ell) \in \Z_3^\ell} \vt{L}{\eta_1-i\gm_1}{\otau^i}{\vep_1}
\otimes \cdots \otimes \vt{L}{\eta_\ell-i\gm_\ell}{\otau^i}{\vep_\ell}.
\end{equation*}
Moreover, for the irreducible $V_{L_{C\times D}}^{\htau}$-module
$\vt{L_{C\times D}}{\eta}{\htau^i}{r} $, $r=0,1,2$ we have
\begin{equation*}
\vt{L_{C\times D}}{\eta}{\htau^i}{r} \cong \bigoplus_{(\gm_1,\ldots,\gm_\ell) \in D}
\bigoplus_{\vep_1 + \cdots + \vep_\ell \equiv r \pmod{3}} \vt{L}{\eta_1-i\gm_1}{\otau^i}{\vep_1}
\otimes \cdots \otimes \vt{L}{\eta_\ell-i\gm_\ell}{\otau^i}{\vep_\ell}.
\end{equation*}
\end{thm}

\section{Modules of $V_L^{\otau}$}
\label{section:VLt}

In this section we recall the classification of irreducible
$V_L^{\otau}$-modules in \cite{TY} and compute some fusion rules for $V_{L}^{\otau}$.

\begin{prop}\label{prop:irr-L}{\rm \cite{TY}}
$V_{L}^{\otau}$ is a simple, rational, $C_2$-cofinite, and CFT type vertex operator algebra.
There are exactly $30$ inequivalent irreducible $V_{L}^{\otau}$-modules.
Their representatives are
$V_{L^{(0,j)}}(\varepsilon)$,
$V_{L^{(c,j)}}$ and 
$\vt{L}{k}{\otau^i}{\varepsilon}$ for $i=1,2$ and $j,k,\varepsilon=0,1,2$.
\end{prop}

We need the structure of each irreducible $V_{L}^{\otau}$-module to compute 
certain fusion rules. Let $M_k^i,W_k^i,M_t^j,W_t^{j},M_k^0(\varepsilon)$ and $W_k^0(\varepsilon)$
be as in \cite{DLTYY,TY}.
Then $M_k^0, M_k^0(0)$ and $M_t^{0}$ are simple vertex operator algebras.
Set $M^0=M_k^0(0)\otimes M_t^0$ and $W^0=W_k^0(0)\otimes W_t^0$.
Then $V_{L}^{\otau}=M^0\oplus W^0$ and
\begin{align}
V_{L^{(0,j)}}(\varepsilon)&\cong M_k^{0}(\varepsilon)\otimes
M_t^{j}\oplus W_k^{0}(\varepsilon)\otimes W_t^{j},
\nonumber\\
V_{L^{(c,j)}}&\cong M_k^{c}\otimes M_t^{j}\oplus
W_k^{c}\otimes W_t^{j},\quad j,\varepsilon=0,1,2
\label{eq:untwist-L}
\end{align}
as $M^0$-modules \cite[Section 4]{TY}.

Moreover, let $M_T(\otau^i),W_T(\otau^i),M_T(\otau^i)(\varepsilon)$ and $W_T(\otau^i)(\varepsilon)$
be as in \cite{DLTYY,TY}.
Then, for $j,\varepsilon\in\Z_3$,
\begin{align}
\vt{L}{j}{\otau}{\varepsilon}&\cong
M_T(\otau)(\varepsilon)\otimes M^{-j}_t\oplus W_T(\otau)(\varepsilon)\otimes W^{-j}_t,\nonumber\\
\vt{L}{j}{\otau^2}{\varepsilon}&\cong
M_T(\otau^2)(\varepsilon)\otimes M^{j}_t\oplus W_T(\otau^2)(\varepsilon)\otimes W^{j}_t
\label{eq:twist-L}
\end{align}
as $M^{0}$-modules \cite[Section 4]{TY}.

\begin{prop}\label{prop:5/6}
{\rm\cite{DLTYY}}
$M_k^{0}(0)$ is a  simple, rational, $C_2$-cofinite, and CFT type vertex operator algebra.
There are exactly $20$ inequivalent irreducible $M_k^{0}(0)$-modules.
Their representatives are
$M_k^{0}(\varepsilon),
W_k^{0}(\varepsilon),
M_k^{c},
W_k^{c},
M_T(\otau^i)(\varepsilon)$,
and
$W_T(\otau^i)(\varepsilon)$ for $\varepsilon=0,1,2$ and $i=1,2$.
\end{prop}

\begin{prop}\label{prop:4/5}{\rm \cite{M}}
$M_t^{0}$ is a  simple, rational, $C_2$-cofinite, 
and  CFT type vertex operator algebra.
There are exactly $6$ inequivalent irreducible $M_t^{0}$-modules.
Their representatives are
$M_t^{j}$ and $W_t^{j}$ for $j=0,1,2$.
The fusion rules for $M_t^0$ are as follows.
\begin{align}
M^i_t\times M^j_t&=M^{i+j}_{t},\nonumber\\
M^i_t\times W^j_t&=W^{i+j}_{t},\nonumber\\
W^i_t\times W^j_t&=M^{i+j}_{t}+W^{i+j}_{t}\label{eq:fusionM_t}
\end{align}
for $i,j=0,1,2$.
\end{prop}

We compute some fusion rules for $V_{L}^{\otau}$.
\begin{lem}\label{lem:fusion-leq}
Let $\varepsilon,\varepsilon_1,\varepsilon_2,j,j_1,j_2,k\in\Z_3$ and $i=1,2$.
Then
\begin{align}
V_{L^{(0,j_1)}}(\varepsilon_1)\times V_{L^{(0,j_2)}}(\varepsilon_2)
& \leq 
V_{L^{(0,j_1+j_2)}}(\varepsilon_1+\varepsilon_2),\label{eq:fusion-leq-L-1}\\
V_{L^{(0,j_1)}}(\varepsilon)\times V_{L^{(c,j_2)}}
& \leq V_{L^{(c,j_1+j_2)}},\label{eq:fusion-leq-L-2}\\
V_{L^{(c,j_1)}}\times V_{L^{(c,j_2)}}
& \leq \sum_{\rho=0}^{2}V_{L^{(0,j_1+j_2)}}(\rho)+2V_{L^{(c,j_1+j_2)}},
\label{eq:fusion-leq-L-3}\\
V_{L^{(0,j)}}(\varepsilon_1)\times \vt{L}{k}{\otau^i}{\varepsilon_2}
& \leq 
\vt{L}{k-ij}{\otau^i}{i\varepsilon_1+\varepsilon_2},\label{eq:fusion-leq-L-4}\\
V_{L^{(c,j)}}\times \vt{L}{k}{\otau^i}{\varepsilon}
& \leq \sum_{\rho=0}^{2}\vt{L}{k-ij}{\otau^i}{\rho}.\label{eq:fusion-leq-L-5}
\end{align}
\end{lem}
\begin{proof}
We have the following fusion rules of irreducible $M_k^{0}(0)$-modules.
\begin{align}
M_k^{0}(\varepsilon_1)\times M_k^{0}(\varepsilon_2)&=
M_k^{0}(\varepsilon_1+\varepsilon_2),\nonumber\\
M_k^{0}(\varepsilon)\times M_k^{c}&=M_k^{c},\nonumber\\
M_k^{c}\times M_k^{c}&=\sum_{\rho=0}^{2}M_{k}^{0}(\rho)+2M_k^{c},\nonumber\\
M_k^{0}(\varepsilon_1)\times
M_{T}(\otau^i)(\varepsilon_2)&\leq 
M_{T}(\otau^i)(i\varepsilon_1+\varepsilon_2),\nonumber\\
M_k^{c}\times M_{T}(\otau^i)(\varepsilon)&\leq
\sum_{\rho=0}^{2}M_{T}(\otau^i)(\rho).\label{eq:fusion-M}
\end{align}

The first three fusion rules can be found in \cite[Theorem 4]{Tanabe}
and we can show the last two formulas by applying
the same method used there. We shall sketch the proof.
In \cite{Tanabe}, $M_k^{0}(0)$
and $M_k^{0}(\varepsilon)$ are denoted by  ${\mathcal W}$ and
$M_k^{0 (\varepsilon)}$, respectively, and $M_k^a$ is used instead of $M_k^{c}$.
Let $A(M_k^{0}(0))$ be the Zhu algebra of $M_k^{0}(0)$
and let
$A(M_k^{0}(\varepsilon_1))$, $A(M_k^{c})$ be the $A(M_k^{0}(0))$-bimodules
introduced in \cite{FZ}. In \cite{DLTYY}, it is shown that
$A(M_k^{0}(0))$ is generated  by two elements $[\omega]$ and $[J]$.
Their action on the top level of every  irreducible 
$M_k^{0}(0)$-module are also
computed there.
Using these data and \cite[Proposition 2.10]{Li0}, 
the same argument as in
\cite[Theorem 4]{Tanabe} shows the last two formulas in (\ref{eq:fusion-M}).

By (\ref{eq:fusionM_t}), (\ref{eq:fusion-M}), and
\cite[Proposition 2.10]{DMZ}, we have fusion rules for $M^{0}$ as
follows.
\begin{align}
\label{eq:fusion-Mkt}
M_k^{0}(\varepsilon_1)\otimes M_t^{k_1}\times
M_k^{0}(\varepsilon_2) \otimes M_t^{k_2}&=
M_k^{0}(\varepsilon_1+\varepsilon_2)\otimes M_t^{k_1+k_2},\nonumber\\
M_k^{0}(\varepsilon)\otimes M_t^{k_1}\times M_k^{c}\otimes M_t^{k_2}&=
M_k^{c}\otimes M_t^{k_1+k_2},\nonumber\\
M_k^{c}\otimes M_t^{k_1}\times M_k^{c}\otimes M_t^{k_2}&=
\sum_{\rho=0}^{2}M_{k}^{0}(\rho)\otimes M_t^{k_1+k_2}+2M_k^{c} 
\otimes M_t^{k_1+k_2},\nonumber\\
M_k^{0}(\varepsilon_1)\otimes M_t^{k_1}\times
M_{T}(\otau^i)(\varepsilon_2)\otimes M_t^{k_2}&\leq
M_{T}(\otau^i)(i\varepsilon_1+\varepsilon_2)\otimes M_t^{k_1+k_2},\nonumber\\
M_k^{c}\otimes M_t^{k_1}\times M_{T}(\otau^i)(\varepsilon)\otimes M_t^{k_2}&\leq
\sum_{\rho=0}^{2}M_{T}(\otau^i)(\rho)\otimes M_t^{k_1+k_2},
\end{align}
where $k_1,k_2\in\Z_3$. Let $N$ be an irreducible
$V_L^{\otau}$-module.
By Propositions \ref{prop:irr-L}--\ref{prop:4/5},
(\ref{eq:untwist-L}), (\ref{eq:twist-L}), and \cite[(3.25)]{TY},
there exist irreducible
$M^0$-modules $M_N$ and $W_N$ such that
\begin{align*}
N&=M_N\oplus W_N,\\
W^0\times M_N&=W_N,\\
W^0\times W_N&=M_N+W_N
\end{align*}
as $M^0$-modules. These $M_N$ and $W_N$ are uniquely determined by $N$.

For $V_L^{\otau}$-modules $N^1,N^2$ and $N^3$, 
\begin{equation}
\dim_{\C}I_{V_L^{\otau}}\binom{N^3}{N^1\ N^2}\leq
\dim_{\C}I_{M^{0}}\binom{N^3}{M_{N^1}\ M_{N^2}}\label{eq:upper-M}
\end{equation}
by \cite[Proposition 11.9]{DL} and
\begin{equation}
I_{M^0}\binom{N^3}{M_{N^1}\ M_{N^2}}\cong
I_{M^{0}}\binom{M_{N^3}}{M_{N^1}\ M_{N^2}}
\oplus
I_{M^{0}}\binom{W_{N^3}}{M_{N^1}\ M_{N^2}}\label{eq:plus-M}
\end{equation}
as vector spaces.
The assertion follows from (\ref{eq:untwist-L}), (\ref{eq:twist-L}),
(\ref{eq:fusion-Mkt}), (\ref{eq:upper-M}), and
(\ref{eq:plus-M}).
\end{proof}

For $\mu\in \Klein^{\ell}$, $C(\mu)$ denotes the $\Klein$-code
generated by $\mu$ and $\htau(\mu)$. Note that $C(\mu)$ is
$\htau$-invariant since $\mu + \htau(\mu)+\htau^2(\mu)={\mathbf 0}$,
where ${\mathbf 0} = (0,\ldots,0)$.
For $\gamma\in\Z_3^{\ell}$, $D(\gamma)$ denotes the $\Z_3$-code
generated by $\gamma$.
These symbols will be used in this section, 
Sections 5, and 7.

\begin{prop}\label{prop:fusionL}
Let $\varepsilon,\varepsilon_1,\varepsilon_2,j,j_1,j_2,k\in\Z_3$ and $i=1,2$.
Then
\begin{align}
V_{L^{(0,j_1)}}(\varepsilon_1)\times V_{L^{(0,j_2)}}(\varepsilon_2)
& = V_{L^{(0,j_1+j_2)}}(\varepsilon_1+\varepsilon_2),\label{eq:fusion-L-1}\\
V_{L^{(0,j_1)}}(\varepsilon)\times V_{L^{(c,j_2)}}
& = V_{L^{(c,j_1+j_2)}},\label{eq:fusion-L-2}\\
V_{L^{(c,j_1)}}\times V_{L^{(c,j_2)}}
& = \sum_{\rho=0}^{2}V_{L^{(0,j_1+j_2)}}(\rho)+2V_{L^{(c,j_1+j_2)}},
\label{eq:fusion-L-3}\\
V_{L^{(0,j)}}(\varepsilon_1)\times \vt{L}{k}{\otau^i}{\varepsilon_2}
& = \vt{L}{k-ij}{\otau^i}{i\varepsilon_1+\varepsilon_2},\label{eq:fusion-L-4}\\
V_{L^{(c,j)}}\times \vt{L}{k}{\otau^i}{\varepsilon}
& = \sum_{\rho=0}^{2}\vt{L}{k-ij}{\otau^i}{\rho}.\label{eq:fusion-L-5}
\end{align}
\end{prop}
\begin{proof}
Restricting intertwining operators for $V_{L}$ in 
Lemma \ref{lem:fusion-lattice} to 
$V_{L}^{\htau}$-modules,
we have
\begin{align}
V_{L^{(0,j_1)}}(\varepsilon_1)\times V_{L^{(0,j_2)}}(\varepsilon_2)
& \geq V_{L^{(0,j_1+j_2)}}(\varepsilon_1+\varepsilon_2),\nonumber\\
V_{L^{(0,j_1)}}(\varepsilon)\times V_{L^{(c,j_2)}}
& \geq V_{L^{(c,j_1+j_2)}},\nonumber\\
V_{L^{(c,j_1)}}\times V_{L^{(c,j_2)}}
& \geq \sum_{\rho=0}^{2}V_{L^{(0,j_1+j_2)}}(\rho)+2V_{L^{(c,j_1+j_2)}},
\label{eq:geq-untwist}
\end{align}
where 
$\dim_{\C}
I_{V_{L}^{\htau}}\binom{V_{L^{(c,j_1+j_2)}}}{V_{L^{(c,j_1)}}\ V_{L^{(c,j_2)}}}\geq 2$
follows from the same arguments as in the proof of \cite[Lemma 6 (2)]{Tanabe}.
By Lemma \ref{lem:fusion-leq} and \eqref{eq:geq-untwist}, we have
(\ref{eq:fusion-L-1})--(\ref{eq:fusion-L-3}).

We shall show \eqref{eq:fusion-L-4} and \eqref{eq:fusion-L-5} for $i=1$. 
Note that 
$L_{\zero\times D(1^6)}$ and $L_{C(c^6)\times D(1^6)}$ are even lattices 
by Lemma \ref{lem:even-lattice},
where $(c^6)=(c,c,c,c,c,c)\in\Klein^{6}$ and $(1^6)=(1,1,1,1,1,1)\in\Z_3^{6}$.
We use the lattice vertex operator algebras
$V_{L_{\zero\times D(1^6)}}$ and $V_{L_{C(c^6)\times D(1^6)}}$
instead of $V_{L_{\zero\times D(1)}}$ and $V_{L_{C(c)\times D(1)}}$
since the lattices $L_{\zero\times D(1)}$ and 
$L_{C(c)\times D(1)}$ are not even.
By Theorem \ref{thm:twist-L-st},
\begin{equation}
V_{L_{\zero\times D(1^6)}}^{T,{\bf 0}}(\htau)\cong
\bigoplus_{k=0}^{2}
\bigoplus_{(\rho_1,\ldots,\rho_6)\in\Z_3^{6}}
\bigotimes_{s=1}^{6}
\vt{L}{-k}{\otau}{\rho_s}.
\label{eq:st-six}
\end{equation}
For $j,k,\varepsilon_1,\varepsilon_2\in \Z_3$,
\begin{equation*}
V_{L^{(0,j)}}(\varepsilon_1)^{\otimes 6}\cdot
\vt{L}{k}{\otau}{\varepsilon_2}^{\otimes 6}\subset
\vt{L}{k-j}{\otau}{\varepsilon_1+\varepsilon_2}^{\otimes 6}
\end{equation*}
in $V_{L_{\zero\times D(1^6)}}^{T,{\bf 0}}(\htau)$ by
\eqref{eq:fusion-leq-L-4} and \eqref{eq:st-six}. Since
$V_{L_{\zero\times D(1^6)}}^{T,{\bf 0}}(\htau)$ is
irreducible, we have $V_{L_{\zero\times D(1^6)}}\cdot
\vt{L}{k}{\otau}{\varepsilon_2}^{\otimes 6}= V_{L_{\zero\times
D(1^6)}}^{T,{\bf 0}}(\htau)$ and
\begin{equation}
V_{L^{(0,j)}}(\varepsilon_1)^{\otimes 6}\cdot
\vt{L}{k}{\otau}{\varepsilon_2}^{\otimes 6}=
\vt{L}{k-j}{\otau}{\varepsilon_1+\varepsilon_2}^{\otimes 6}
\label{eq:eq-int}
\end{equation}
in $V_{L_{\zero\times D(1^6)}}^{T,{\bf 0}}(\htau)$. Let
$\pr \colon V_{L_{\zero\times D(1^6)}}^{T,{\bf 0}}(\htau)
\rightarrow
\vt{L}{k-j}{\otau}{\varepsilon_1+\varepsilon_2}^{\otimes 6}$ be
a projection. For $u\in V_{L^{(0,j)}}(\varepsilon_1)^{\otimes 6}, 
v\in \vt{L}{k}{\otau}{\varepsilon_2}^{\otimes 6}$, set
$f(u,x)v=\pr Y_{V_{L_{\zero\times D(1^6)}}^{T,\zero}(\htau)}(u,x)v$. Then $f(\cdot,x)$ is a nonzero intertwining
operator of type
$\binom{\vt{L}{k-j}{\otau}{\varepsilon_1+\varepsilon_2}^{\otimes 6}} {V_{L^{(0,j)}}(\varepsilon_1)^{\otimes 6}\
\vt{L}{k}{\otau}{\varepsilon_2}^{\otimes 6}}$ for $(V_L^{\otau})^{\otimes 6}$
by (\ref{eq:eq-int}).
It follows from \cite[Proposition 11.9]{DL} and \cite[Proposition 2.10]{DMZ}
that \eqref{eq:fusion-L-4} holds for $i=1$.

By Theorem \ref{thm:twist-L-st},
\begin{equation}
V_{L_{C(c^6)\times D(1^6)}}^{T,{\bf 0}}(\htau)\cong
\bigoplus_{k=0}^{2}
\bigoplus_{(\rho_1,\ldots,\rho_6)\in\Z_3^{6}}
\bigotimes_{m=1}^{6}
\vt{L}{-k}{\otau}{\rho_{m}}.\label{eq:st-six-2}
\end{equation}
Since $V_{L_{C(c^6)\times D(1^6)}}$ is simple , 
it
follows from Lemma \ref{lem:lattice-st} and \eqref{eq:fusion-leq-L-2} that
\begin{equation*}
\big(\bigotimes_{m=1}^{6}V_{L^{(0,0)}}(\nu_{m})\big)\cdot
V_{L^{(c,j)}}^{\otimes 6}=
V_{L^{(c,j)}}^{\otimes 6}
\end{equation*}
in $V_{L_{C(c^6)\times D(1^6)}}$ for
$\nu_1,\ldots,\nu_6\in\Z_3$. Therefore,
\begin{align}
&V_{L^{(c,j)}}^{\otimes 6}\cdot
\vt{L}{k}{\otau}{\varepsilon}^{\otimes 6}\nonumber\\
&=
\big(\big(\bigotimes_{m=1}^{6}V_{L^{(0,0)}}(\nu_m)\big)\cdot V_{L^{(c,j)}}^{\otimes 6}\big)
\cdot \vt{L}{k}{\otau}{\varepsilon}^{\otimes 6}\nonumber\\
&=
\big(\bigotimes_{m=1}^{6}V_{L^{(0,0)}}(\nu_m)\big)\cdot \big(V_{L^{(c,j)}}^{\otimes 6}
\cdot \vt{L}{k}{\otau}{\varepsilon}^{\otimes 6}\big)\label{eq:action}
\end{align}
in $V_{L_{C(c^6)\times D(1^6)}}^{T,{\bf 0}}(\htau)$.
For $j,k,\varepsilon\in\Z_3$, \eqref{eq:fusion-leq-L-5} and \eqref{eq:st-six-2}
imply
\begin{equation}
V_{L^{(c,j)}}^{\otimes 6}\cdot
\vt{L}{k}{\otau}{\varepsilon}^{\otimes 6}\subset
\bigoplus_{(\rho_1,\ldots,\rho_6)\in\Z_3^6}
\bigotimes_{m=1}^{6}
\vt{L}{k-j}{\otau}{\rho_m}\label{eq:include}
\end{equation}
and for $\nu_1,\ldots,\nu_6,\rho_1,\ldots,\rho_6\in\Z_3$,
\eqref{eq:fusion-leq-L-4} implies
\begin{align}
\big(\bigotimes_{m=1}^{6}V_{L^{(0,0)}}(\nu_m)\big)\cdot
\big(\bigotimes_{m=1}^{6}
\vt{L}{k-j}{\otau}{\rho_m}\big)
\subset &\bigotimes_{m=1}^{6}\vt{L}{k-j}{\otau}{\nu_m+\rho_m}
\label{eq:action-twist}
\end{align}
in $V_{L_{C(c^6)\times D(1^6)}}^{T,{\bf 0}}(\htau)$.
Since
$V_{L_{C(c^6)\times D(1^6)}}^{T,{\bf 0}}(\htau)$ is irreducible,
$V_{L^{(c,j)}}^{\otimes 6}\cdot
\vt{L}{k}{\otau}{\varepsilon}^{\otimes 6}$ is a nonzero
$(V_L^{\otau})^{\otimes 6}$-module.
Since $\bigotimes_{m=1}^{6}
\vt{L}{k-j}{\otau}{\rho_m},
(\rho_1,\ldots,\rho_6)\in\Z_3^6,$ are all inequivalent
irreducible $(V_L^{\otau})^{\otimes 6}(=V_{L^{(0,0)}}(0)^{\otimes 6})$-modules,
there exists
$(\rho_1^{\prime},\ldots,\rho_6^{\prime})\in\Z_3^6$ such that
\begin{equation*}
V_{L^{(c,j)}}^{\otimes 6}\cdot
\vt{L}{k}{\otau}{\varepsilon}^{\otimes 6}\supset
\bigotimes_{m=1}^{6}
\vt{L}{k-j}{\otau}{\rho_m^{\prime}}
\end{equation*}
by (\ref{eq:include}).
By (\ref{eq:action}) and (\ref{eq:action-twist}), we have
\begin{equation}
V_{L^{(c,j)}}^{\otimes 6}\cdot
\vt{L}{k}{\otau}{\varepsilon}^{\otimes 6}=
\bigoplus_{(\rho_1,\ldots,\rho_6)\in\Z_3^6}
\bigotimes_{m=1}^{6}
\vt{L}{k-j}{\otau}{\rho_m}.\label{eq:twist-2-6}
\end{equation}
For $\rho=(\rho_1,\ldots,\rho_{6})\in\Z_3^{6}$, let $\pr_{\rho} \colon V_{L_{C(c^6)\times D(1^6)}}^{T,{\bf 0}}(\htau)
\rightarrow \bigotimes_{m=1}^{6}\vt{L}{k-j}{\otau}{\rho_m}$
be a projection.
For
$u\in V_{L^{(c,j)}}^{\otimes 6},
v\in \vt{L}{k}{\otau}{\varepsilon}^{\otimes 6}$,
set $f_{\rho}(u,x)v=\pr_{\rho}Y_{V_{L_{C(c^6)\times D(1^6)}}^{T,{\bf 0}}(\htau)}(u,x)v$.
Then $f_{\rho}(\cdot,x)$ is a nonzero intertwining operator
\begin{align*}
f_{\rho}(\cdot,x) :
V_{L^{(c,j)}}^{\otimes 6}\rightarrow
\Hom_{\C}\big(\vt{L}{k}{\otau}{\varepsilon}^{\otimes 6},
\bigotimes_{m=1}^{6}\vt{L}{k-j}{\otau}{\rho_m}\big)\{x\}
\end{align*}
for $(V_L^{\otau})^{\otimes 6}$ by (\ref{eq:twist-2-6}). 
Thus, 
\begin{align*}
V_{L^{(c,j)}}\times
\vt{L}{k}{\otau}{\varepsilon}&\geq
\sum_{\rho=0}^{2}\vt{L}{k-j}{\otau}{\rho}
\end{align*}
holds by \cite[Proposition 11.9]{DL} and \cite[Proposition 2.10]{DMZ}
and hence \eqref{eq:fusion-L-5} holds by \eqref{eq:fusion-leq-L-5}.
We can show \eqref{eq:fusion-L-4} and \eqref{eq:fusion-L-5} for $i=2$ similarly.
\end{proof}
\begin{rmk}
We can show that the equalities hold in the last two formulas in \eqref{eq:fusion-M}
by using \eqref{eq:twist-L}, \eqref{eq:fusion-Mkt},
Proposition \ref{prop:fusionL} and \cite[Proposition 11.9]{DL}.
\end{rmk}

\section{Modules of $V_{L^{\oplus {\ell}}}^{\htau}$}

Let ${\ell}$ be a positive integer.
In this section we discuss $V_{L^{\oplus {\ell}}}^{\htau}$-modules,
namely the case $C=\{{\zero}\}$ and $D=\{{\zero}\}$. 
We shall determine some fusion rules for $V_{L^{\oplus {\ell}}}^{\htau}$.

In view of Proposition \ref{prop:fusionL}, we introduce a new index set $\tilde{\Klein}=\{0,1,2,a,b,c\}$
and define a new commutative binary operation on $\tilde{\Klein}$
by
\begin{align*}
\begin{array}{ll}
i+j=i+j\ (\module {3})&\mbox{for } i,j=0,1,2,\\
j+x=x&\mbox{for } j=0,1,2,\ x=a,b,c,\\
x+x=0&\mbox{for } x=a,b,c,\\
\multicolumn{2}{l}{a+b=c,\quad b+c=a,\quad c+a=b.}
\end{array}
\end{align*}
Then, $\tilde{\Klein}$ contains $\Z_3$ and $\Klein$.
Note that this binary operation is not associative.
We use $\tilde{\Klein}$ to describe fusion rules for $(V_{L}^{\otau})^{\otimes\ell}$
in \eqref{eq:local-fusion-1}.
Define an action of $\otau$ on $\tilde{\Klein}$ by
$\otau(a) = b$, $\otau(b) = c$, $\otau(c) = a$, and $\otau(j) = j$,
$j=0,1,2$, which is compatible with the action of $\otau$ on $\Z_3$ and $\Klein$.
This action of $\otau$ preserves the binary operation on $\tilde{\Klein}$.
The set of $\otau$-orbits on $\tilde{\Klein}$ is $\{0,1,2,c\}$.
We consider the componentwise action of $\tgrs$ on $\tilde{\Klein}^{\ell}$ and the
componentwise binary operation on $\tilde{\Klein}^{\ell}$. 
The symmetric group $\mathfrak{S}_{\ell}$ acts on $\tilde{\Klein}^{\ell}$ 
by permuting the components
and so $\tgrl$ acts on $\tilde{\Klein}^{\ell}$ naturally.
For $\lambda = (\lambda_1,\ldots,\lambda_\ell)
\in \tilde{\Klein}^{\ell}$, its {\em support} is defined to be
$\supp_{\tilde{\Klein}}(\lambda) = \{ i\,|\, \lambda_i\in\{a,b,c\}\}$. 
The cardinality of
$\supp_{\tilde{\Klein}}(\lambda)$ is called the {\em weight} of $\lambda$. We denote the
weight of $\lm$ by $\wt_{\tilde{\Klein}}(\lambda)$. 
For $\lambda^1,\lambda^2\in \tilde{\Klein}^{\ell}$,
we write $\lambda^1\equiv_{\htau}\lambda^2$ if
$\lambda^1$ and $\lambda^2$ belong to the same orbit of $\htau= (\otau, \ldots, \otau)$ in $\tilde{\Klein}^{\ell}$.
We denote by $(\tilde{\Klein^{\ell}})_{\equiv_{\htau}}$ the set of all orbits  of $\htau$ in
$\tilde{\Klein}^{\ell}$.
For a $\htau$-invariant subset $P$ of $\tilde{\Klein}^{\ell}$,
$P_{\equiv_{\htau}}$ denotes the set of all orbits  of $\htau$ in
$P$.

By Proposition \ref{prop:irr-L} and \cite[Proposition 2.7]{DMZ},
$(V_L^{\otau})^{\otimes {\ell}}$ is a rational and $C_2$-cofinite
vertex operator algebra. Moreover,
\begin{equation}
\{U^1\otimes\cdots\otimes U^{\ell}\ |\ U^1,\ldots, U^{\ell}
\mbox{ are irreducible $V_L^{\otau}$-modules}\}
\label{eq:irr-te-l}
\end{equation}
is a complete list of
irreducible $(V_L^{\otau})^{\otimes {\ell}}$-modules up to isomorphism.
Set
\begin{align}
\label{eq:mathcalP}
{\mathcal P}_0&=\{U^1\otimes\cdots\otimes U^{\ell}\ |\ U^1,\ldots,U^{\ell}\in
\{V_{L^{(0,j)}}(\varepsilon), V_{L^{(c,j)}}\ |\ j,\varepsilon\in\Z_3\}\},\nonumber\\
{\mathcal P}_1&=\{U^1\otimes\cdots\otimes U^{\ell}\ |\ U^1,\ldots,U^{\ell}\in
\{\vt{L^{\oplus {\ell}}}{k}{\otau}{\varepsilon}\ |\ k,\varepsilon\in\Z_3\}\},\nonumber\\
{\mathcal P}_2&=\{U^1\otimes\cdots\otimes U^{\ell}\ |\ U^1,\ldots,U^{\ell}\in
\{\vt{L^{\oplus {\ell}}}{k}{\otau^2}{\varepsilon}\ |\ k,\varepsilon\in\Z_3\}\},\nonumber\\
{\mathcal P}&={\mathcal P}_0\cup{\mathcal P}_1\cup{\mathcal P}_2.
\end{align}
Set $\grs=\{(\otau^{i_1},\ldots,\otau^{i_{\ell-1}},1)\in\tgrs\ |\
i_1,\ldots,i_{\ell-1}\in\Z\}$. Then $\grs$ acts on $V_{L^{\oplus {\ell}}}^{\htau}$ naturally and
$(V_{L^{\oplus {\ell}}}^{\htau})^{\grs}=(V_L^{\otau})^{\otimes {\ell}}$.

For $i\in\tilde{\Klein}$ and $j=0,1,2$, set
\begin{equation*}
X_{i,j} = \left\{
\begin{array}{ll}
V_{L^{(0,j)}}(i) &\mbox{if }i=0,1,2,\\
V_{L^{(i,j)}} &\mbox{if }i=a,b,c.
\end{array}\right.
\end{equation*}
For $\xi=(\xi_1,\ldots,\xi_{\ell})\in\tilde{\Klein}^{\ell}$ and
$\gamma=(\gamma_1,\ldots,\gamma_{\ell})\in\Z_3^{\ell}$,
set
\begin{equation*}
X_{\xi,\gamma}=\bigotimes_{i=1}^{\ell}X_{\xi_i,\gamma_i}.
\end{equation*}
Then, 
for $\lambda\in\Klein^{\ell}$ and
$\gamma\in\Z_3^{\ell}$ Lemma \ref{lem:lattice-st} implies that 
$V_{L_{(\lambda,\gamma)}}=\bigoplus_{\xi}X_{\xi,\gamma}$ and
\begin{align}
\label{eq:deco-L-l-0}
V_{L_{(g(\lambda),\gamma)}}&=\bigoplus_{\xi}
X_{g(\xi),\gamma},
\end{align}
where $g\in\tgrs$ and $\xi$ runs over the set
$\{\xi=(\xi_1,\ldots,\xi_{\ell})\in\tilde{\Klein}^{\ell}\ |\
\xi_k=\lambda_k \mbox{ for all }k\in \supp_{\Klein}(\lambda)\}$.
This observation will be used in the argument just after \eqref{eq:schur-weyl}.

We have 
$X_{\xi,\gamma}\cong X_{g(\xi),\gamma}$
as
$(V_L^{\otau})^{\otimes {\ell}}$-modules for $g\in \tgrs$
since $V_{L^{(\tau^{i}(c),j)}}\cong V_{L^{(c,j)}}$ as 
$V_L^{\otau}$-modules for $i,j\in\Z_3$.
Thus, we can choose $\xi$ to be an element of $\{0,1,2,c\}^{\ell}$
when we deal with 
$X_{\xi,\gamma}$ as $(V_L^{\otau})^{\otimes {\ell}}$-modules.
Using this notation,
we can describe some fusion rules for $(V_L^{\otau})^{\otimes {\ell}}$
by Proposition \ref{prop:fusionL} and \cite[Proposition 2.10]{DMZ} as follows:
\begin{align}
X_{\rho,{\gamma^1}}\times
X_{\xi,\gamma^2}=X_{\rho+\xi,\gamma^1+\gamma^2},
\label{eq:local-fusion-1}
\end{align}
for $\rho\in\Z_3^{\ell}, \xi\in\{0,1,2,c\}^{\ell}$,
and $\gamma^1,\gamma^2\in\Z_3^{\ell}$.

For any $
{\mathbf 0}\neq\lambda\in\Klein^{\ell},
\gamma\in\Z_3^{\ell}$, and $\varepsilon=0,1,2$,
set 
\begin{align}
P(V_{L_{(\zero,\gamma)}}(\varepsilon))
&=
\{\xi=(\xi_k)\in\Z_3^{\ell}\ |\ \sum_{k=1}^{\ell}\xi_k\equiv \varepsilon\pmod{3}\},\nonumber\\
P(V_{L_{(\lambda,\gamma)}})
&=\{\xi\in\{0,1,2,c\}^{\ell}\ |\
\supp_{\tilde{\Klein}}(\xi)=\supp_{\Klein}(\lambda)\}.
\label{eq:P-xi}
\end{align}
Then, Lemma \ref{lem:lattice-st} implies that
\begin{align}
V_{L_{({\mathbf 0},\gamma)}}(\varepsilon)&\cong
\bigoplus_{\xi\in P(V_{L_{(\zero,\gamma)}}(\varepsilon))}X_{\xi,\gamma},\nonumber\\
V_{L_{(\lambda,\gamma)}}&\cong\bigoplus_{\xi\in P(V_{L_{(\lambda,\gamma)}})}
X_{\xi,\gamma}
\label{eq:deco-L-l}
\end{align}
as $(V_L^{\otau})^{\otimes {\ell}}$-modules.
In particular,
we have
\begin{equation}
V_{L^{\oplus \ell}}^{\htau}\cong
\bigoplus_{\begin{subarray}{l}
\rho=(\rho_i)\in\Z_3^{\ell}\\
\rho_1+\cdots+\rho_{\ell}=0
\end{subarray}}
X_{\rho,{\mathbf 0}}\label{eq:sim-v}
\end{equation}
as $(V_{L}^{\otau})^{\otimes \ell}$-modules.

We have already seen in \eqref{eq:actionD-un} and Section 3
that
for $\lambda\in\Klein^{\ell}$, $\gamma\in\Z_3^{\ell}$,
$\eta\in\Z_3^{\ell}$, and $g\in\tgrs$
\begin{align}
V_{L_{(\lambda,\gamma)}}\circ g\cong
g^{-1}\big(V_{L_{(\lambda,\gamma)}}\big)=V_{L_{(g^{-1}(\lambda),\gamma)}},\nonumber\\
V_{L^{\oplus {\ell}}}^{T,\eta}(\htau^i)\circ g\cong
g^{-1}(V_{L^{\oplus {\ell}}}^{T,\eta}(\htau^i))
=V_{L^{\oplus {\ell}}}^{T,\eta}(\htau^i)
\label{eq:acH}
\end{align}
as $V_{L^{\oplus {\ell}}}$-modules or $\htau^i$-twisted $V_{L^{\oplus {\ell}}}$-modules.
Hence for any ${\mathbf 0}\neq \lambda\in \Klein^{\ell},
\gamma\in\Z_3^{\ell}, \varepsilon=0,1,2$, and $g\in \grs$,
\begin{gather}
V_{L_{({\mathbf 0},\gamma)}}(\varepsilon)\circ g\cong
V_{L_{({\mathbf 0},\gamma)}}(\varepsilon),
\qquad
V_{L_{(\lambda,\gamma)}}\circ g\cong
V_{L_{(g^{-1}(\lambda),\gamma)}},\nonumber\\
\vt{L^{\oplus {\ell}}}{\eta}{\htau^i}{\varepsilon}\circ g\cong
\vt{L^{\oplus {\ell}}}{\eta}{\htau^i}{\varepsilon}
\label{eq:G-stable}
\end{gather}
as $V_{L^{\oplus {\ell}}}^{\htau}$-modules.

\begin{lem}\label{lem:res-fusion}
Let $N$ be an $\N$-graded weak $V_{L^{\oplus {\ell}}}^{\htau}$-module.
Then any irreducible $(V_{L}^{\otau})^{\otimes {\ell}}$-submodule of $N$
is isomorphic to an element of ${\mathcal P}$
as defined in \eqref{eq:mathcalP}.
\end{lem}
\begin{proof}
Let $U$ be an irreducible $(V_{L}^{\otau})^{\otimes
{\ell}}$-submodule of $N$. By (\ref{eq:irr-te-l}), there are
irreducible $V_{L}^{\otau}$-modules $U^1,\ldots,U^{\ell}$ such
that $U\cong U^1\otimes \cdots\otimes U^{\ell}$. Set
$S=V_{L^{\oplus {\ell}}}^{\htau}\cdot U$. For the same reason as
in \cite[Proof of Lemma 5.2]{TY}, $S$ is an ordinary $V_{L^{\oplus
{\ell}}}^{\htau}$-module. Moreover, $Y_N(v,x)U\neq 0$  for any
nonzero $v\in V_{L^{\oplus {\ell}}}^{\htau}$.

Set
\begin{align*}
Q_0&=\{i\in \{1,\ldots,\ell\}\ |\ U^{i}\in
\{V_{L^{(0,j)}}(\varepsilon), V_{L^{(c,j)}}\ |\ j,\varepsilon\in\Z_3\}\},\\
Q_1&=\{i\in \{1,\ldots,\ell\}\ |\ U^{i}\in
\{\vt{L^{\oplus {\ell}}}{k}{\otau}{\varepsilon}\ |\ k,\varepsilon\in\Z_3\}\},\\
Q_2&=\{i\in \{1,\ldots,\ell\}\ |\ U^{i}\in
\{\vt{L^{\oplus {\ell}}}{k}{\otau^2}{\varepsilon}\ |\ k,\varepsilon\in\Z_3\}\}.
\end{align*}

Let $\omega_L$ be the Virasoro element of $V_{L}^{\otau}$.
By \cite[Section 4]{TY}, the eigenvalues of $(\omega_L)_1$ on the top levels of
irreducible $V_{L}^{\otau}$-modules are
\begin{equation}
\begin{array}{rll}
0&\mbox{for }V_{L^{(0,0)}}(0),\\
1&\mbox{for }V_{L^{(0,0)}}(\varepsilon), & \varepsilon=1,2,\\
2/3&\mbox{for }V_{L^{(0,j)}}(\varepsilon), & j=1,2,\varepsilon=0,1,2,\\
1/2&\mbox{for }V_{L^{(c,0)}},\\
1/6&\mbox{for }V_{L^{(c,j)}}, & j=1,2,\\
1/9&\mbox{for }\vt{L}{0}{\otau^i}{0}\mbox{ and }\vt{L}{j}{\otau^i}{2},& i=1,2,j=1,2,\\
4/9&\mbox{for }\vt{L}{0}{\otau^i}{2}\mbox{ and }\vt{L}{j}{\otau^i}{1},& i=1,2,j=1,2,\\
7/9&\mbox{for }\vt{L}{0}{\otau^i}{1}\mbox{ and }\vt{L}{j}{\otau^i}{0},& i=1,2,j=1,2.
\end{array}\label{eq:eigen}
\end{equation}

Let $W^1$ be an irreducible $V_{L}^{\otau}$-module and let
$W^{1,r}, r=0,1,2$ be the irreducible $V_{L}^{\otau}$-module
determined by the fusion rule $V_{L^{(0,0)}}(r) \times
W^1=W^{1,r}$ in Proposition \ref{prop:fusionL}. Let
$\lambda_{1}$ and $\lambda_{1,r}$ be the eigenvalues of
$(\omega_L)_1$ on the top levels of $W^{1}$ and $W^{1,r}$,
respectively. By (\ref{eq:eigen}), we have
\begin{align}
\lambda_{1,r}-\lambda_1\equiv&
\left\{
\begin{array}{ll}
0&\mbox{if }W^1\in\{V_{L^{(0,j)}}(\varepsilon), V_{L^{(c,j)}}\ |\ j,\varepsilon=0,1,2\},\\
2r/3&\mbox{if }W^1\in\{\vt{L}{j}{\otau}{\varepsilon}\ |\ j,\varepsilon=0,1,2\},\\
r/3&\mbox{if }W^1\in\{\vt{L}{j}{\otau^2}{\varepsilon}\ |\ j,\varepsilon=0,1,2\}.
\end{array}
\right.\pmod{\Z}
\label{eq:eigen-dif}
\end{align}

Let $\omega$ be the Virasoro element of $V_{L^{\oplus
\ell}}^{\htau}$. Assume that $Q_s, Q_t\neq\varnothing, s\neq t$.
Take $i_s\in Q_s$ and $i_t\in Q_t$ and define
$\rho=(\rho_i)\in\Z_3^{\ell}$ by
\begin{equation*}
\rho_i =
\begin{cases}
1 &\text{if }i=i_s,\\
2 &\text{if }i=i_t,\\
0 &\text{otherwise}.
\end{cases}
\end{equation*}
By \eqref{eq:sim-v}, $X_{\rho,\zero}$ is an irreducible
$(V_{L}^{\otau})^{\otimes {\ell}}$-submodule of $V_{L^{\oplus
\ell}}^{\htau}$. Using (\ref{eq:eigen-dif}), Proposition
\ref{prop:fusionL}, and \cite[Proposition 2.10]{DMZ}, 
one can show that $S$ has a $(V_{L}^{\otau})^{\otimes
{\ell}}$-submodule $W$ such that the difference of the minimal
eigenvalues of $\omega_1$ in $W$ and in $U$ is not an integer
since $0\neq X_{\rho,\zero}\cdot U\subset S$. This is a
contradiction. Hence the assertion holds.
\end{proof}

\begin{lem}\label{lem:generate}
Let $N$ be an $\N$-graded weak $V_{L^{\oplus {\ell}}}^{\htau}$-module.
Let $M$ be an irreducible $(V_{L}^{\otau})^{\otimes \ell}$-submodule
of $N$ and $N^1$ the $V_{L^{\oplus \ell}}^{\htau}$-submodule of $N$
generated by $M$.
Then $N^1$ is isomorphic to one of the following
inequivalent irreducible $V_{L^{\oplus \ell}}^{\htau}$-modules.

$(1)$ $V_{L_{({\mathbf 0},\gamma)}}(\varepsilon)$, $\gamma\in\Z_3^{\ell}, \varepsilon=0,1,2$.

$(2)$ $V_{L_{(\lambda,\gamma)}}$, ${\mathbf 0}\neq\lambda\in (\Klein^{\ell})_{\equiv_{\htau}}$,
$\gamma\in\Z_3^{\ell}$.

$(3)$ $\vt{L^{\oplus \ell}}{\eta}{\htau^i}{\varepsilon}$, $\eta\in \Z_3^{\ell}, i=1,2, \varepsilon=0,1,2$.
\end{lem}
\begin{proof}
By \eqref{eq:acH}, we have
\begin{gather*}
V_{L_{({\mathbf 0},\gamma)}}\circ \htau\cong
V_{L_{({\mathbf 0},\gamma)}}
,
\qquad
V_{L_{(\lambda,\gamma)}}\circ \htau\cong V_{L_{(\htau^{-1}(\lambda),\gamma)}}\not\cong V_{L_{(\lambda,\gamma)}}
\end{gather*}
as $V_{L^{\oplus\ell}}$-modules
for $\zero\neq\lambda\in\Klein^{\ell}$ and $\gamma\in\Z_3^{\ell}$
and 
\begin{gather*}
V_{L^{\oplus {\ell}}}^{T,\eta}(\htau^i)\circ \htau\cong
\htau^{-1}(V_{L^{\oplus {\ell}}}^{T,\eta}(\htau^i))=V_{L^{\oplus {\ell}}}^{T,\eta}(\htau^i)
\end{gather*}
as $\htau^i$-twisted $V_{L^{\oplus\ell}}$-modules
for $\eta\in\Z_3^{\ell}$.
It follows from \cite[Theorem 2]{MT}  that
the $V_{L^{\oplus {\ell}}}^{\htau}$-modules in the above list are irreducible and inequivalent.

By Lemma \ref{lem:res-fusion}, $M$ is an element of ${\mathcal P}$ in \eqref{eq:mathcalP}.
Suppose $M\in{\mathcal P}_0$, that is,
$M\cong X_{\xi,\gamma}, \xi\in\{0,1,2,c\}^{\ell},\gamma\in\Z_3^{\ell}$.
Set $\Xi=\{\rho+\xi\in\{0,1,2,c\}^{\ell}\ |\
\rho=(\rho_i)\in\Z_3^{\ell}, \sum_{i=1}^{\ell}\rho_i=0\}$.
Since $(V_L^{\otau})^{\otimes \ell}$ is a rational
vertex operator algebra,
$N^1$ is a direct sum of $(V_{L}^{\otau})^{\otimes\ell}$-modules.
By \eqref{eq:local-fusion-1} and \eqref{eq:sim-v},
we can write $N^1=\oplus_{j\in {\mathcal J}}M^j$
where each $M^j$ is isomorphic to $X_{\nu^j,\gamma},
\nu^j\in\Xi$.
We can take $M^{j_1}=M$ for some $j_1\in {\mathcal J}$. Let
$\pr_{j} \colon N^1\rightarrow M^j, j\in{\mathcal J}$ be
projections. For any $j\in{\mathcal J}, u\in X_{\rho,{\mathbf 0}}\subset V_{L^{\oplus \ell}}^{\htau}, v\in M$, define
$f_j(u,x)v=\pr_j(Y_N(u,x)v)$. Then $f_j\in
I_{(V_L^{\otau})^{\otimes \ell}}\binom{M^j}{X_{\rho,{\mathbf 0}}\
M}$. For each $\nu\in\Xi$, we see from (\ref{eq:local-fusion-1})
that there is at most one $j\in{\mathcal J}$ such that $M^{j}\cong
X_{\nu,\gamma}$ (cf. \cite[Proof of Lemma 5.6]{TY}).

Assume that $X_{\rho,{\mathbf 0}}\cdot M=0$ for
$X_{\rho,{\mathbf 0}}\subset V_{L^{\oplus \ell}}^{\htau}$. Then
\begin{align*}
0&=V_{L^{\oplus \ell}}^{\htau}\cdot (X_{\rho,{\mathbf 0}}\cdot M)
=(V_{L^{\oplus \ell}}^{\htau}\cdot X_{\rho,{\mathbf 0}})\cdot M\\
&=V_{L^{\oplus \ell}}^{\htau}\cdot M\supset M
\end{align*}
since $V_{L^{\oplus \ell}}^{\htau}$ is simple. This is a
contradiction. Hence $0\neq X_{\rho,{\mathbf 0}}\cdot M$, and
consequently $X_{\rho,{\mathbf 0}}\cdot M\cong
X_{\rho+\xi,\gamma}$ as $(V_L^{\otau})^{\otimes \ell}$-modules.
Therefore, we have
\begin{equation}
N^1\cong\bigoplus_{\nu\in\Xi}
X_{\nu,\gamma}\label{eq:simple}
\end{equation}
as $(V_L^{\otau})^{\otimes \ell}$-modules. Applying the above
arguments to $V_{L^{\oplus \ell}}^{\htau}$-module $N^1$, we conclude that
$N^1$ is irreducible.

By \cite[Theorem 6.14]{DY}, if two irreducible $V_{L^{\oplus{\ell}}}^{\htau}$-modules $W^1, W^2$ have an isomorphic
irreducible $(V_L^{\otau})^{\otimes {\ell}}$-submodule, then there
exists $g\in \grs$ such that $W^1\circ g\cong W^2$. Hence by
\eqref{eq:deco-L-l} and \eqref{eq:G-stable}, $N^1$ is
isomorphic to $V_{L_{({\mathbf 0},\gamma)}}(\varepsilon),
\varepsilon\in\Z_3$ or $V_{L_{(\lambda,\gamma)}}, {\mathbf0}\neq\lambda\in\Klein^{\ell}$.

For $i=1,2$, we see from Theorem \ref{thm:twist-L-st} that every
irreducible $(V_L^{\otau})^{\otimes {\ell}}$-module in ${\mathcal
P}_i$ appears in  the irreducible $V_{L^{\oplus {\ell}}}^{\htau}$-modules 
listed in (3). Hence one can show that if $M\in {\mathcal P}_i$, then 
$N^1\cong \vt{L^{\oplus\ell}}{\eta}{\htau^i}{\varepsilon}, \eta\in\Z_3^{\ell},\varepsilon\in\Z_3$ similarly.
\end{proof}

\begin{prop}\label{prop:irr-L-l}
$V_{L^{\oplus {\ell}}}^{\htau}$ is a simple, rational, $C_2$-cofinite, and CFT type  vertex operator algebra.
The following is a complete set of representatives
of equivalence classes of irreducible $V_{L^{\oplus {\ell}}}^{\htau}$-modules.

$(1)$ $V_{L_{({\mathbf 0},\gamma)}}(\varepsilon)$, $\gamma\in\Z_3^{\ell}, \varepsilon=0,1,2$.

$(2)$ $V_{L_{(\lambda,\gamma)}}$, ${\mathbf 0}\neq\lambda\in (\Klein^{\ell})_{\equiv_{\htau}}$,
$\gamma\in\Z_3^{\ell}$.

$(3)$ $\vt{L^{\oplus\ell}}{\eta}{\htau^i}{\varepsilon}$, $\eta\in \Z_3^{\ell}, i=1,2, \varepsilon=0,1,2$.
\end{prop}
\begin{proof}
By \eqref{eq:sim-v} and \cite{Buhl}, $V_{L^{\oplus{\ell}}}^{\htau}$ is a $C_2$-cofinite vertex operator algebra. The
classification of irreducible $V_{L^{\oplus\ell}}^{\htau}$-modules follows from Lemma \ref{lem:generate}.
Since $(V_{L}^{\otau})^{\otimes \ell}$ is rational, the rationality of $V_{L^{\oplus \ell}}^{\htau}$ follows
from Lemma  \ref{lem:generate}.
\end{proof}

The following lemma gives lower bounds for some fusion rules for $V_{L^{\oplus {\ell}}}^{\htau}$.
\begin{lem}\label{lem:fusion-L-l-geq}
Let $\lambda,\lambda^1,\lambda^2$ be nonzero elements of $\Klein^{\ell}$ such that
$\lambda^1\not\equiv_{\htau}\lambda^2$, $\gamma,\gamma^1,\gamma^2,\eta\in\Z_3^{\ell}$, $i=1,2$,
and $\varepsilon,\varepsilon_1,\varepsilon_2=0,1,2$.
Then
\begin{align}
V_{L_{({\mathbf 0},\gamma^1)}}(\varepsilon_1)\times V_{L_{({\mathbf 0},\gamma^2)}}(\varepsilon_2)
&\geq
V_{L_{({\mathbf 0},\gamma^1+\gamma^2)}}(\varepsilon_1+\varepsilon_2),
\label{eq:fusion-Ll-1-geq}
\\
V_{L_{({\mathbf 0},\gamma^1)}}(\varepsilon)\times V_{L_{(\lambda,\gamma^2)}}
&\geq
V_{L_{(\lambda,\gamma^1+\gamma^2)}},
\label{eq:fusion-Ll-2-geq}\\
V_{L_{(\lambda^1,\gamma^1)}}\times
V_{L_{(\lambda^2,\gamma^2)}}
&\geq
\sum_{j=0}^{2}V_{L_{(\lambda^1+\htau^j(\lambda^2),\gamma^1+\gamma^2)}},
\label{eq:fusion-Ll-3-geq}\\
V_{L_{(\lambda,\gamma^1)}}\times
V_{L_{(\lambda,\gamma^2)}}
&\geq
\sum_{\rho=0}^{2}V_{L_{({\mathbf 0},\gamma^1+\gamma^2)}}(\rho)+2
V_{L_{(\lambda,\gamma^1+\gamma^2)}},
\label{eq:fusion-Ll-4-geq}
\\
V_{L_{({\mathbf 0},\gamma)}}(\varepsilon_1)\times
\vt{L^{\oplus {\ell}}}{\eta}{\htau^i}{\varepsilon_2}
&\geq
\vt{L^{\oplus {\ell}}}{\eta-i\gamma}{\htau^i}{i\varepsilon_1+\varepsilon_2},
\label{eq:fusion-Ll-t1-geq}
\\
\label{eq:fusion-Ll-t2-geq}
V_{L_{(\lambda,\gamma)}}\times
\vt{L^{\oplus {\ell}}}{\eta}{\htau^i}{\varepsilon}
&\geq
\sum_{\rho=0}^{2}\vt{L^{\oplus {\ell}}}{\eta-i\gamma}{\htau^i}{\rho}.
\end{align}
\end{lem}
\begin{proof}
Restricting intertwining operators for $V_{L^{\oplus\ell}}$ in Lemma \ref{lem:fusion-lattice} to 
$V_{L^{\oplus\ell}}^{\htau}$-modules, we have
\eqref{eq:fusion-Ll-1-geq}--\eqref{eq:fusion-Ll-4-geq},
where 
$\dim_{\C}
I_{V_{L^{\oplus\ell}}^{\htau}}
\binom{V_{L_{(\lambda,\gamma^1+\gamma^2)}}}{V_{L_{(\lambda,\gamma^1)}}\ V_{L_{(\lambda,\gamma^2)}}}
\geq 2$
follows from the same arguments as in the proof of \cite[Lemma 6 (2)]{Tanabe}.

We shall show \eqref{eq:fusion-Ll-t2-geq} for $i=1$. 
\eqref{eq:fusion-Ll-t1-geq} and \eqref{eq:fusion-Ll-t2-geq} for $i=2$
can be proved by a similar argument.
It is easy to see that
\begin{equation}
I_{V_{L^{\oplus {\ell}}}^{\htau}}\binom{M}
{V_{L_{(\lambda,\gamma)}}\ \vt{L^{\oplus\ell}}{\eta}{\htau}{\varepsilon}}=0
\end{equation}
for all $M\not\cong\vt{L^{\oplus\ell}}{\eta-\gamma}{\htau}{r}, r=0,1,2$,
by Proposition \ref{prop:fusionL} and \cite[Proposition 2.10]{DMZ}.
Set $\tilde{\lambda}=(\lambda,\lambda,\lambda,\lambda,\lambda,\lambda)\in \Klein^{6\ell}$ and
$\tilde{\gamma}=(\gamma,\gamma,\gamma,\gamma,\gamma,\gamma)\in \Z_3^{6\ell}$.
Recall that $C(\tilde{\lambda})$ is the $\Klein$-code generated by $\tilde{\lambda}$ and $\htau(\tilde{\lambda})$ and that
$D(\tilde{\gamma})$ is the $\Z_3$-code generated by $\tilde{\gamma}$ (cf. Section \ref{section:VLt}).
Lemma \ref{lem:even-lattice} implies
$L_{C(\tilde{\lambda})\times D(\tilde{\gamma})}$ is a $\htau$-invariant even lattice.
To obtain (\ref{eq:fusion-Ll-t2-geq}),
we use the lattice vertex operator algebra
$V_{L_{C(\tilde{\lambda})\times D(\tilde{\gamma})}}$ instead of $V_{L_{(\lambda,\gamma)}}$
since the lattice $L_{(\lambda,\gamma)}$ is not even.
Let $\eta\in \Z_3^{\ell}$ and set $\tilde{\eta}=(\eta,\eta,\eta,\eta,\eta,\eta)\in\Z_3^{6\ell}$.
Consider a $\htau$-twisted $V_{L_{C(\tilde{\lambda})\times D(\tilde{\gamma})}}$-module $V_{L_{C(\tilde{\lambda})\times D(\tilde{\gamma})}}^{T,\tilde{\eta}}(\htau)$.
It follows from 
\eqref{eq:twist-D-L} that
\begin{align}
V_{L_{C(\tilde{\lambda})\times D(\tilde{\gamma})}}^{T,\tilde{\eta}}(\htau)
\cong&
\bigoplus_{j=0}^{2}
V_{L^{\oplus 6\ell}}^{T,{\tilde{\eta}-j\tilde{\gamma}}}(\htau)
\label{eq:st-t-L-6}
\end{align}
as $\htau$-twisted $V_{L^{\oplus 6\ell}}$-modules.
We have
\begin{align}
V_{L_{C(\tilde{\lambda})\times D(\tilde{\gamma})}}^{T,\tilde{\eta}}(\htau)
\cong&\bigoplus_{j=0}^{2}
\bigoplus_{\rho_1,\ldots, \rho_{6}\in\Z_3}
\bigotimes_{i=1}^{6}\vt{L^{\oplus {\ell}}}{\eta-j\gamma}{\htau}{ \rho_i}
\label{eq:six-l}
\end{align}
as $(V_{L^{\oplus {\ell}}}^{\htau})^{\otimes 6}$-modules
by the same argument as was used in the proof of Theorem \ref{thm:twist-L-st} 
by replacing $V_{L}$, $V_{L^{\oplus\ell}}$, $C$, $D$, and $\eta$ by $V_{L^{\oplus\ell}}$, $V_{L^{\oplus6\ell}}$,
$C(\tilde{\lambda})$, $D(\tilde{\gamma})$, and $\tilde{\eta}$, respectively.
Since $V_{L_{C(\tilde{\lambda})\times D(\tilde{\gamma})}}$ is simple,
it follows from Lemma \ref{lem:lattice-st} and \eqref{eq:fusion-L-2} that
\begin{equation*}
\big(\bigotimes_{m=1}^{6}V_{L_{({\mathbf 0},{\mathbf 0})}}(\nu_m)\big)\cdot
V_{L_{(\lambda,\gamma)}}^{\otimes 6}=
V_{L_{(\lambda,\gamma)}}^{\otimes 6}
\end{equation*}
in $V_{L_{C(\tilde{\lambda})\times D(\tilde{\gamma})}}$ for
$\nu_1,\ldots,\nu_6\in\Z_3$. Therefore,
\begin{align}
V_{L_{(\lambda,\gamma)}}^{\otimes 6}\cdot
\vt{L^{\oplus {\ell}}}{\eta}{\htau}{\varepsilon}^{\otimes 6}
&=
\big(\big(\bigotimes_{m=1}^{6}V_{L_{({\mathbf 0},{\mathbf 0})}}(\nu_m)\big)
\cdot V_{L_{(\lambda,\gamma)}}^{\otimes 6}\big)
\cdot \vt{L^{\oplus {\ell}}}{\eta}{\htau}{\varepsilon}^{\otimes 6}\nonumber\\
&=
\big(\bigotimes_{m=1}^{6}V_{L_{({\mathbf 0},{\mathbf 0})}}(\nu_m)\big)
\cdot \big(V_{L_{(\lambda,\gamma)}}^{\otimes 6}
\cdot \vt{L^{\oplus {\ell}}}{\eta}{\htau}{\varepsilon}^{\otimes 6}
\big)\label{eq:action-L-l}
\end{align}
in $V_{L_{C(\tilde{\lambda})\times D(\tilde{\gamma})}}^{T,\tilde{\eta}}(\htau)$.
By Proposition \ref{prop:fusionL} and \eqref{eq:six-l},
\begin{align*}
V_{L_{(\lambda,\gamma)}}^{\otimes 6}\cdot
\big(\vt{L^{\oplus {\ell}}}{\eta}{\htau}{\varepsilon}\big)^{\otimes 6}
\subset&\bigoplus_{\rho_1,\ldots, \rho_{6}\in\Z_3}
\bigotimes_{m=1}^{6}\vt{L^{\oplus {\ell}}}{\eta-\gamma}{\htau}{\rho_m}
\end{align*}
and for
$\nu_1,\ldots,\nu_6,\rho_1,\ldots,\rho_6 \in\Z_3$,
\begin{equation}
\big(\bigotimes_{m=1}^{6}V_{L_{({\mathbf 0},{\mathbf 0})}}(\nu_{m})\big)\cdot
\big(\bigotimes_{m=1}^{6}\vt{L^{\oplus {\ell}}}{\eta-\gamma}{\htau}{\rho_m}\big)
\subset \bigotimes_{m=1}^{6}\vt{L^{\oplus {\ell}}}{\eta-\gamma}{\htau}{\nu_{m}+\rho_{m}}
\label{eq:action-L-l-twist}
\end{equation}
in $V_{L_{C(\tilde{\lambda})\times D(\tilde{\gamma})}}^{T,\tilde{\eta}}(\htau)$.
Since $V_{L_{C(\tilde{\lambda})\times D(\tilde{\gamma})}}^{T,\tilde{\eta}}(\htau)$
is a $\htau$-twisted irreducible $V_{L_{C(\tilde{\lambda})\times D(\tilde{\gamma})}}$-module,
$V_{L_{(\lambda,\gamma)}}^{\otimes 6}\cdot
\big(\vt{L^{\oplus {\ell}}}{\eta}{\htau}{\varepsilon}\big)^{\otimes 6}$
is a nonzero $(V_{L^{\oplus {\ell}}}^{\htau})^{\otimes 6}$-module.
Since $\otimes_{m=1}^{6}\vt{L^{\oplus {\ell}}}{\eta-\gamma}{\htau}{\rho_m}$, 
$(\rho_1,\ldots,\rho_6)\in\Z_3^6$,
are all inequivalent irreducible $(V_{L^{\oplus {\ell}}}^{\htau})^{\otimes 6}$-modules,
there exists $(\rho_1^{\prime},\ldots,\rho_6^{\prime})\in\Z_3^{6}$ such that
\begin{align}
V_{L_{(\lambda,\gamma)}}^{\otimes 6}\cdot
\big(\vt{L^{\oplus {\ell}}}{\eta}{\htau}{\varepsilon}\big)^{\otimes 6}
\supset \bigotimes_{i=1}^{6}\vt{L^{\oplus {\ell}}}{\eta-\gamma}{\htau}{\rho_i^{\prime}}.
\label{eq:fusion-inc}
\end{align}

By \eqref{eq:action-L-l}--\eqref{eq:fusion-inc}, we have
\begin{equation}
V_{L_{(\lambda,\gamma)}}^{\otimes 6}\cdot
\vt{L^{\oplus {\ell}}}{\eta}{\htau}{\varepsilon}^{\otimes 6}
=\bigoplus_{\rho_1,\ldots,\rho_6\in \Z_3}\bigotimes_{i=1}^{6}
\vt{L^{\oplus {\ell}}}{\eta-\gamma}{\htau}{\rho_i}.
\label{L-l-twist-pro}
\end{equation}
Using the same argument as in the proof of \eqref{eq:fusion-L-5}, we have
\eqref{eq:fusion-Ll-t2-geq}.
\end{proof}

We want to use the results in \cite{DY} and \cite{Tanabe}.
We follow the notation of \cite{DY}.
Note that we can take all $2$-cocyles in \cite{DY} to be trivial in our setting.
Let ${\mathcal S}$ be a finite $\grs$-stable set of irreducible $V_{L^{\oplus\ell}}^{\htau}$
-modules (cf. Section 2.1).
Set
${\mathcal M}=\oplus_{M\in{\mathcal S}}M$.
Note that $\grs$ acts on ${\mathcal M}$ by \eqref{eq:actionD-un},
\eqref{eq:Hlacts-1}, and \eqref{eq:Hlacts-2}.
Define a vector space $\C{\mathcal S}=\oplus_{M\in{\mathcal S}}e(M)$ with formal basis $e(M)$, $M\in{\mathcal S}$. The space $\C{\mathcal S}$ is
an associative algebra under the product $e(M)e(N)=\delta_{M,N}e(M)$.
Define the vector space ${\mathcal A}(\grs,{\mathcal S})=\C[\grs]\otimes \C{\mathcal S}$
with basis $g\otimes e(M)$ for $g\in \grs$ and $M\in{\mathcal S}$,
and a multiplication on it by:
\begin{equation*}
g\otimes e(M)\cdot h\otimes e(N)=gh\otimes e(h^{-1}(M))e(N).
\end{equation*}
Then ${\mathcal A}(\grs,{\mathcal S})$ is an associative algebra with the identity element
$\sum_{M\in{\mathcal S}}1\otimes e(M)$.
We define an action of ${\mathcal A}(\grs,{\mathcal S})$ on ${\mathcal M}$ as follows:
For $M,N\in{\mathcal S}, w\in N$ and $g\in \grs$,
we set
\begin{equation}
g\otimes e(M)\cdot w=\delta_{M,N}gw.\label{eq:actionMi}
\end{equation}

For $M\in{\mathcal S}$, define a
subgroup $(\grs)_{M}=\{g\in\grs\ |\ g(M)=M\}$ of $\grs$ and define
subalgebras $s(M)=\spn_{\C}\{g\otimes e(M)\ |\ g\in(\grs)_{M}\}$ and
$D(M)=\spn_{\C}\{g\otimes e(M)\ |\ g\in\grs\}$ of ${\mathcal A}(\grs,{\mathcal S})$. 
Note that $s(M)$ is isomorphic to the
group algebra of $(\grs)_{M}$.
Decompose ${\mathcal S}$ into a disjoint union of $\grs$-orbits 
${\mathcal S}=\cup_{j\in J}{\cal O}_j$. 
Let $M^{(j)}$ be a representative of ${\mathcal O}_j$.

We shall compute some fusion rules for $V_{L^{\oplus {\ell}}}^{\htau}$
in Proposition \ref{prop:fusion-L-l}
by using \cite[Theorem 2]{Tanabe}.
We need the following result which gives
a complete set of representatives of isomorphism classes of
irreducible ${\mathcal A}(\grs,{\mathcal S})$-modules.

\begin{thm} {\rm\cite[Theorem 3.6]{DY}}\label{thm:AGS}
${\mathcal A}(\grs, {\mathcal S})$ is semisimple and the
irreducible ${\mathcal A}(\grs, {\mathcal S})$-modules are precisely
$D(M^{(j)})\otimes_{s(M^{(j)})}U$,
where $U$ ranges over the irreducible $s(M^{(j)})$-modules and $j\in J$.
\end{thm}

Note that $\grs$ acts on $(\Klein^{\ell})_{\equiv\htau}$.
Let $\zero\neq\lambda\in\{0,c\}^{\ell}$ and $\gamma\in\Z_3^{\ell}$ and
set 
${\mathcal R}_{\lambda}$ be the $\grs$-orbit in $(\Klein^{\ell})_{\equiv\htau}$ containing $\lambda$.
Then 
\begin{align}
\label{eq:slg}
{\mathcal S}_{\lambda,\gamma}&=\{V_{L_{(\mu,\gamma)}}\ |\ \mu\in{\mathcal R}_{\lambda}\}
\end{align} 
is an $\grs$-stable set.
We shall describe the irreducible ${\mathcal A}(\grs,{\mathcal S}_{\lambda,\gamma})$-modules
in Proposition \ref{lem:s-module}.
Theorem \ref{thm:AGS} implies that the irreducible ${\mathcal A}(\grs,{\mathcal S}_{\lambda,\gamma})$-modules
are obtained by the irreducible $s(V_{L_{(\lambda,\gamma)}})$-modules.
In order to classify the irreducible $s(V_{L_{(\lambda,\gamma)}})$-modules,
we first investigate the action of $s(V_{L_{(\lambda,\gamma)}})$ on $V_{L_{(\lambda,\gamma)}}$.
We recall the decomposition
$V_{L_{(\lambda,\gamma)}}=\oplus_{\xi\in P(V_{L_{(\lambda,\gamma)}})}
X_{\xi,\gamma}$ in \eqref{eq:deco-L-l}, where $P(V_{L_{(\lambda,\gamma)}})$ is given in \eqref{eq:P-xi}.
For $g\in \grs$,
$g$ is an element in $(\grs)_{V_{L_{(\lambda,\gamma)}}}$
if and only if $g\lambda\equiv_{\htau}\lambda$.
Thus, $(\grs)_{V_{L_{(\lambda,\gamma)}}}$ consists of the 
elements  
\begin{align}
(\otau^{j_1},\ldots,\otau^{j_{\ell}})(\otau,\ldots,\otau)^{-j_{\ell}}\in \grs
\label{eq:fix-H}
\end{align}
with $j_{k}=0$
for all $k\in\supp_{\Klein}(\lambda)$.
Note that $|(\grs)_{V_{L_{(\lambda,\gamma)}}}|=|P(V_{L_{(\lambda,\gamma)}})|=3^{\ell-\wt_{\Klein}(\lambda)}$.
We have
\begin{align*}
gu&=(\otau^{j_1},\ldots,\otau^{j_{\ell}})(\otau,\ldots,\otau)^{-j_{\ell}}u\\
&=\zeta_3^{\sum_{k}j_k\xi_k}(\otau,\ldots,\otau)^{-j_{\ell}}u\in 
(\otau,\ldots,\otau)^{-j_{\ell}}(X_{\xi,\gamma})
\end{align*}
for $\xi=(\xi_i)\in P(V_{L_{(\lambda,\gamma)}})$,
$u\in X_{\xi,\gamma}$ and $g\in(\grs)_{V_{L_{(\lambda,\gamma)}}}$ of the form \eqref{eq:fix-H},
where we define $0 c=0$ in the sum $\sum_{k}j_k\xi_k$.
Note that 
the linear map $X_{\xi,\gamma}\ni u\mapsto
(\otau,\ldots,\otau)^{-j_{\ell}}u\in 
(\otau,\ldots,\otau)^{-j_{\ell}}(X_{\xi,\gamma})$
is an isomorphism of $(V_{L}^{\otau})^{\otimes\ell}$-modules induced by
the isomorphism $(\otau,\ldots,\otau)^{-j_{\ell}}\colon V_{L_{(\lambda,\gamma)}}\rightarrow
 V_{L_{(\htau^{-j_{\ell}}(\lambda),\gamma)}}$ of 
$V_{L^{\oplus\ell}}^{\htau}$-modules.

For $\xi=(\xi_i)\in\{0,1,2,c\}^{\ell}$, $\C e(\xi)$ denotes a vector space
with formal basis $e(\xi)$. 
In view of the above observation, we define an action of $s(V_{L_{(\lambda,\gamma)}})$ on $\C e(\xi)$
by setting 
\begin{align*}
g\otimes e(V_{L_{(\lambda,\gamma)}})\cdot e(\xi)
&=\zeta_3^{\sum_{k}j_k\xi_k}e(\xi)
\end{align*}
for $g\in(\grs)_{V_{L_{(\lambda,\gamma)}}}$ of the form \eqref{eq:fix-H} and 
$\xi\in P(V_{L_{(\lambda,\gamma)}})$, 
where we define $0 c=0$ in the sum $\sum_{k}j_k\xi_k$.

Let $\gamma\in\Z_3^{\ell}$ and $\varepsilon \in\Z_3$.
Then $\{V_{L_{(\zero,\gamma)}}(\varepsilon)\}$ is an $\grs$-stable set.
For the same reason as in the case of ${\mathcal S}_{\lambda,\gamma}$ discussed above,
we define an action of $s(V_{L_{(\zero,\gamma)}}(\varepsilon))$ on $\C e(\xi)$
by setting 
\begin{align*}
g\otimes e(V_{L_{(\zero,\gamma)}}(\varepsilon))\cdot e(\xi)
&=\zeta_3^{\sum_{k}j_k\xi_k}e(\xi)
\end{align*}
for $g=(\otau^{j_1},\ldots,\otau^{j_{\ell-1}},1)\in (\grs)_{V_{L_{(\zero,\gamma)}}(\varepsilon)}
=\grs$ and $\xi\in P(V_{L_{(\zero,\gamma)}}(\varepsilon))$.

We have the following result.

\begin{lem}\label{lem:s-module}
With the above notation, the following assertions hold.

$(1)$
$3^{\ell-\wt_{\Klein}(\lambda)}$ inequivalent irreducible $s({V_{L_{(\lambda,\gamma)}}})$-modules
$\C e(\xi), \xi\in P(V_{L_{(\lambda,\gamma)}})$,
form a complete 
set of irreducible $s({V_{L_{(\lambda,\gamma)}}})$-modules up to isomorphism
and for nonzero $u\in X_{\xi,\gamma}$, $\C u\cong \C e(\xi)$ as $s(V_{L_{(\lambda,\gamma)}})$-modules.
Moreover,
\begin{equation}
\{D(V_{L_{(\lambda,\gamma)}})\otimes_{s(V_{L_{(\lambda,\gamma)}})}\C
e(\xi)\ |\ \xi\in P(V_{L_{(\lambda,\gamma)}})\}
\label{eq:irr-a-1}
\end{equation}
is a complete set of irreducible ${\mathcal A}(\grs,{\mathcal S}_{\lambda,\gamma})$-modules 
up to isomorphism and 
\begin{align}
\dim_{\C}D(V_{L_{(\lambda,\gamma)}})\otimes_{s(V_{L_{(\lambda,\gamma)}})}\C
e(\xi)&=|{\mathcal R}_{\lambda}|=3^{\wt_{\Klein}(\lambda)-1}
\label{eq:dimW-1}
\end{align}
for $\xi\in P(V_{L_{(\lambda,\gamma)}})$.

$(2)$
$3^{\ell-1}$ inequivalent irreducible $s(V_{L_{(\zero,\gamma)}}(\varepsilon))$-modules
$\C e(\xi), \xi\in P(V_{L_{(\zero,\gamma)}}(\varepsilon))$,
form a complete set of irreducible $s(V_{L_{(\zero,\gamma)}}(\varepsilon))$-modules up to isomorphism
and for nonzero $u\in X_{\xi,\gamma}$, $\C u\cong \C e(\xi)$ as $s(V_{L_{(\zero,\gamma)}}(\varepsilon))$-modules.
Moreover,
\begin{equation}
\{D(V_{L_{(\zero,\gamma)}}(\varepsilon))\otimes_{s(V_{L_{(\zero,\gamma)}}(\varepsilon))}\C
e(\xi)\ |\ \xi\in P(V_{L_{(\zero,\gamma)}}(\varepsilon))\}
\label{eq:irr-a-2}
\end{equation}
is a complete set of irreducible 
${\mathcal A}(\grs,\{V_{L_{(\zero,\gamma)}}(\varepsilon)\})$-modules up to isomorphism and 
\begin{align}
D(V_{L_{(\zero,\gamma)}}(\varepsilon))\otimes_{s(V_{L_{(\zero,\gamma)}}(\varepsilon))}\C
e(\xi)&\cong \C e(\xi)
\label{eq:dimW-2}
\end{align}
as vector spaces for $\xi\in P(V_{L_{(\zero,\gamma)}}(\varepsilon))$.
\end{lem}
\begin{proof}
We show the first assertion.
The argument just before the lemma shows 
that $\C u\cong \C e(\xi)$ as $s(V_{L_{(\lambda,\gamma)}})$-modules for nonzero $u\in X_{\xi,\gamma}$.
It is clear that $|{\mathcal R}_{\lambda}|=3^{\wt_{\Klein}(\lambda)-1}$.
We have \eqref{eq:dimW-1} since
$|\grs|=3^{\ell-1}$ and $|(\grs)_{V_{L_{(\lambda,\gamma)}}}|=3^{\ell-\wt_{\Klein}(\lambda)}$.
It follows from $|P(V_{L_{(\lambda,\gamma)}})|=3^{\ell-\wt_{\Klein}(\lambda)}$
that $\{\C e(\xi)\ |\ \xi\in P(V_{L_{(\lambda,\gamma)}})\}$ is 
a complete set of irreducible $s({V_{L_{(\lambda,\gamma)}}})$-modules up to isomorphism.
It follows from Theorem \ref{thm:AGS} that \eqref{eq:irr-a-1} is a complete 
set of irreducible 
${\mathcal A}(\grs,{\mathcal S}_{\lambda,\gamma})$-modules up to isomorphism.

The second assertion can be obtained by a similar argument.
\end{proof}
We want to use the result \cite[Theorem 2]{Tanabe} in 
Proposition \ref{prop:fusion-L-l}.
Let $\lambda,\lambda^1,\lambda^2$ be nonzero elements of $\Klein^{\ell}$ such that
$\lambda^1\not\equiv_{\htau}\lambda^2$ and let $\gamma^1,\gamma^2$ be elements of $\Z_3^{\ell}$.
Set $\gamma^3=\gamma^1+\gamma^2$, ${\mathcal R}_i={\mathcal R}_{\lambda^i}$,
and ${\mathcal S}_i={\mathcal S}_{\lambda^i,\gamma^i}$ for $i=1,2$ (cf. \eqref{eq:slg}).
For $i=1,2$, set
$\xi^i=(\xi^i_j)\in\{0,c\}^{\ell}$ by
\begin{equation}
\label{eq:lambda-xi}
\xi^i_j=\left\{\begin{array}{ll}
0 &\mbox{if }\lambda^i_j=0,\\
c &\mbox{if }\lambda^i_j=a,b,c.
\end{array}\right.
\end{equation}
Note that $\xi^i\in{\mathcal R}_i$ and ${\mathcal R}_i=\{\mu\in(\Klein^{\ell})_{\equiv_{\htau}}\ |\
\supp_{\Klein}(\mu)=\supp_{\Klein}(\lambda^i)\}$ for $i=1,2$.
Set ${\mathcal S}_{3}=
\{V_{L_{(\mu,\gamma^3)}}, V_{L_{({\mathbf 0},\gamma^3)}}(\varepsilon)\ |\
{\mathbf 0}\neq \mu\in(\Klein^{\ell})_{\equiv_{\htau}},
\varepsilon=0,1,2\}$. 
For each $i=1,2,3$, ${\mathcal S}_i$ is an $\grs$-stable set.
Set ${\mathcal T}_i=\{V_{L_{(\xi^i,\gamma^i)}}\}, i=1,2$ and
${\mathcal T}_3=\{V_{L_{(\zero,\gamma^3)}}(\varepsilon)\ |\ \varepsilon=0,1,2\}
\cup \{V_{L_{(\mu,\gamma^3)}}\ |\ {\mathbf 0}\neq \mu\in\{0,c\}^{\ell}\}$.
Then,
${\mathcal T}_i$ is a complete set of representatives of $\grs$-orbits in
${\mathcal S}_i$ for $i=1,2,3$.
We simply write $P_i=P(V_{L_{(\xi^i,\gamma^i)}}),i=1,2$ (cf. \eqref{eq:P-xi}).
Let $P_3=\{0,1,2,c\}^{\ell}$.

We note that
\begin{align}
&\bigcup_{\varepsilon=0}^{2}
\{D(V_{L_{({\mathbf 0},\gamma^3)}(\varepsilon)})
\otimes_{s(V_{L_{({\mathbf 0},\gamma^3)}(\varepsilon)})}\C
e(\xi)\ |\ \xi\in P(V_{L_{({\mathbf 0},\gamma^3)}}(\varepsilon))\}\nonumber\\
&\cup\bigcup_{{\mathbf 0}\neq \lambda\in\{0,c\}^{\ell}}
\{D(V_{L_{(\lambda,\gamma^3)}})\otimes_{s(V_{L_{(\lambda,\gamma^3)}})}\C
e(\xi)\ |\ \xi\in P(V_{L_{(\lambda,\gamma^3)}})\}
\label{eq:irr-a-3}
\end{align}
is a complete set of representatives of isomorphism classes of irreducible ${\mathcal A}(\grs,{\mathcal S}_3)$-modules
by Theorem \ref{thm:AGS} and
\begin{align*}
\{0,1,2,c\}^{\ell}&=\bigcup_{\varepsilon=0}^{2}P(V_{L_{({\mathbf 0},\gamma^3)}}(\varepsilon))\cup
\bigcup_{{\mathbf 0}\neq \lambda\in\{0,c\}^{\ell}}P(V_{L_{(\lambda,\gamma^3)}});
\quad \text{disjoint}.
\end{align*}

Set
${\mathcal M}_i=\oplus_{M\in{\mathcal S}_i}M$ for $i=1,2,3$.
For $i=1,2,3$, we write $W_{i,\xi}$ for an irreducible ${\mathcal
A}(\grs,{\mathcal S}_i)$-module $D(M)\otimes_{s(M)}\C e(\xi)$ in
\eqref{eq:irr-a-1}, \eqref{eq:irr-a-2}, and \eqref{eq:irr-a-3} since they are
parametrized by $\xi\in P_i$.
In \eqref{eq:dimW-1} and \eqref{eq:dimW-2}, we have already seen
\begin{align}
\dim_{\C}W_{i,\xi}=3^{\max\{0, \wt_{\tilde{\Klein}}(\xi)-1\}}.
\label{eq:dim-W}
\end{align}
For $i=1,2$ and $\xi\in P_i$, we note that 
$\supp_{\tilde{\Klein}}(\xi)=\supp_{\Klein}(\lambda^i)$ and 
$\dim_{\C}W_{i,\xi}=|{\mathcal R}_i|$. 
We have
\begin{equation*}
{\mathcal M}_i=\bigoplus_{\xi\in P_i}W_{i,\xi}\otimes
\Hom_{{\mathcal A}(\grs,{\mathcal S}_i)}(W_{i,\xi},{\mathcal M}_i)
\end{equation*}
as an ${\mathcal A}(\grs,{\mathcal S}_i)\otimes_{\C}(V_{L}^{\otau})^{\otimes {\ell}}$-module
for $i=1,2,3$.
Then
$\Hom_{{\mathcal A}(\grs,{\mathcal S}_i)}
(W_{i,\xi},{\mathcal M}_i), \xi\in P_i$ are nonzero inequivalent
irreducible $(V_{L}^{\otau})^{\otimes {\ell}}$-modules by \cite[Theorem 6.14]{DY}.
For any $V_{L^{\oplus\ell}}^{\htau}$-module $M$ in ${\mathcal T}_i$
and any
nonzero $u\in X_{\xi,\gamma^i}, \xi\in P_i$, in the decomposition \eqref{eq:deco-L-l} of $M$,
the ${\mathcal A}(\grs,{\mathcal S}_i)$-submodule of ${\mathcal M}_i$ generated by $u$
is isomorphic to $W_{i,\xi}$ since $\C e(\xi)\cong \C u$ as 
$s(M)$-modules by Lemma \ref{lem:s-module}.
Hence there exists a unique $f_v\in\Hom_{{\mathcal A}(\grs,{\mathcal S}_i)}
(W_{i,\xi},{\mathcal M}_i)$ such that $f_{v}(1\otimes e(\xi))=v$.
In fact, the map $v\mapsto f_v$ is a linear isomorphism.
Therefore we identify $\Hom_{{\mathcal A}(\grs,{\mathcal S}_i)}
(W_{i,\xi},{\mathcal M}_i)$ with $X_{\xi,\gamma^i}$
and we write
\begin{equation}
{\mathcal M}_i=\bigoplus_{\xi\in P_i}W_{i,\xi}\otimes X_{\xi,\gamma^i}.
\label{eq:schur-weyl}
\end{equation}

For any $\xi\in P_i$ and any nonzero 
$v^{\xi}\in X_{\xi,\gamma^i}$, 
we can take a basis $\{w^{ij}\ |\ j=1,\ldots,\dim_{\C}W_{i,\xi}\}$
of $W_{i,\xi}$
such that for $j=1,\ldots,\dim_{\C}W_{i,\xi}$, $w^{ij}\otimes v^{\xi}$
is an element of an irreducible $V_{L^{\oplus {\ell}}}^{\htau}$-module in ${\mathcal S}_i$
and if $\dim_{\C}W_{i,\xi}\geq 2$, which implies $\xi\neq\zero$, then
for $j\neq k$, $w^{ij}\otimes v^{\xi}$ and $w^{ik}\otimes v^{\xi}$
belong to different irreducible $V_{L^{\oplus\ell}}^{\otau}$-modules
by \eqref{eq:deco-L-l-0}.
For $i=1,2$, since $\dim_{\C}W_{i,\xi}=|{\mathcal R}_i|$ by \eqref{eq:dim-W},
there exists a bijection 
$\{1,\ldots,\dim_{\C}W_{i,\xi}\}\ni j\mapsto \mu^{i}_j\in {\mathcal R}_i$
where $\mu^i_j$ is determined by $w^{ij}\otimes v^{\xi}\in V_{L_{(\mu^i_j,\gamma^i)}}\in {\mathcal S}_i$.

To see the above situation, we describe the case of $i=1$ as an example.
Let $\{h_1,\ldots,h_{r}\}$ be a complete set of coset representatives of $(\grs)_{V_{L_{(\lambda^1,\gamma^1)}}}
=\{g\in \grs\ |\ g(\lambda^1)\equiv_{\htau}\lambda^1\}$ in $\grs$ where
$r=|\grs/(\grs)_{V_{L_{(\lambda^1,\gamma^1)}}}|$. We recall 
$r=|{\mathcal R}_1|=\dim_{\C}W_{1,\xi}$ for $\xi\in P_1$.
By \eqref{eq:deco-L-l-0}, we have
\begin{align*}
{\mathcal M}_1&=\bigoplus_{\mu\in {\mathcal R}_1}V_{L_{(\mu,\gamma^1)}}
=\bigoplus_{j=1}^{r}V_{L_{(h_j(\lambda^1),\gamma^1)}}\\
&=\bigoplus_{j=1}^{r}\bigoplus_{\xi\in P_1}X_{h_j(\xi),\gamma^1}
=\bigoplus_{\xi\in P_1}\bigoplus_{j=1}^{r}X_{h_j(\xi),\gamma^1}\\
&=\bigoplus_{\xi\in P_1}W_{1,\xi}\otimes X_{\xi,\gamma^1},
\end{align*}
where we identify $\bigoplus_{j=1}^{r}X_{h_j(\xi),\gamma^1}$ with $W_{1,\xi}\otimes X_{\xi,\gamma^1}$.
Let $\xi$ be an element of $P_1$ and $v^{\xi}$  a nonzero element of $X_{\xi,\gamma^1}$.
For $j=1,\ldots,r$,
we can take $h_j(v^{\xi})\in X_{h_j(\xi),\gamma^1}\subset V_{L_{(h_j(\lambda^1),\gamma^1)}}$
as $w^{1j}\otimes v^{\xi}$ in the argument above.

Set
\begin{equation}
{\mathcal I}=\bigoplus_{(M^1,M^2,M^3)
\in{\mathcal S}_1\times{\mathcal S}_2\times{\mathcal S}_3}
I_{V_{L^{\oplus {\ell}}}^{\htau}}\binom{M^3}{M^1\ M^2}\otimes_{\C}M^1\otimes_{\C} M^2.
\end{equation}
Let $M^i\in{\mathcal S}_i$ for $i=1,2,3$.
For $f\in I_{V_{L^{\oplus {\ell}}}^{\htau}}\binom{M^3}{M^1\ M^2}$ and $g\in \grs$,
we define $_{g}f\in I_{V_{L^{\oplus {\ell}}}^{\htau}}
\binom{g(M^3)}{g(M^1)\ g(M^2)}$ as follows: 
For $u\in g(M^1),v\in g(M^2)$,  set
\begin{equation*}
{}_{g}f(u,x)v=g(f(g^{-1}u,x)g^{-1}(v)).
\end{equation*}
We define an action of ${\mathcal A}(\grs,{\mathcal S}_3)$ on ${\mathcal I}$ as follows: 
Let $M^i\in{\mathcal S}_i$ for $i=1,2,3$.
For $g\otimes e(M)\in {\mathcal A}(\grs,{\mathcal S}_3)$, $v\in M^1, w\in M^2$, and
$f\in I_{V_{L^{\oplus {\ell}}}^{\htau}}\binom{M^3}{M^1\ M^2}$,
set
\begin{align*}
&(g\otimes e(M))\cdot (f\otimes v\otimes w)=
\delta_{M,M^3}\cdot {}_{g}f\otimes g(v)\otimes g(w)\\
&\in I_{V_{L^{\oplus {\ell}}}^{\htau}}
\binom{g(M^3)}{g(M^1)\ g(M^2)}\otimes_{\C}
g(M^1)\otimes_{\C}g(M^2).
\end{align*}

Let $\xi^i\in P_i$ for $i=1,2$.
Fix a nonzero $v^{i0}\in X_{\xi^{i},\gamma^i}$.
Set
\begin{equation}
{\mathcal I}(\xi^1,\xi^2)
=
\spn_{\C}\{f\otimes w^1\otimes v^{10}\otimes w^2\otimes v^{20}\in{\mathcal I}\ |\
w^1\in W_{1,\xi^1}, w^2\in W_{2,\xi^2}\},
\label{eq:st-1}
\end{equation}
which is an
${\mathcal A}(\grs,{\mathcal S}_3)$-submodule of ${\mathcal I}$.
It follows from the comments right after (\ref{eq:schur-weyl}) that
\begin{equation}
\dim_{\C}{\mathcal I}(\xi^1,\xi^2)
=\sum_{\mu^1\in {\mathcal R}_1,\mu^2\in{\mathcal R}_2}
\sum_{M^3\in{\mathcal S}_3}
\dim_{\C}I_{V_{L^{\oplus\ell}}^{\htau}}
\binom{M^3}{V_{L_{(\mu^1,\gamma^1)}}\ V_{L_{(\mu^2,\gamma^2)}}}.
\label{eq:dimI}
\end{equation}
We have the following decomposition of ${\mathcal I}(\xi^1,\xi^2)$
as an ${\mathcal A}(\grs,{\mathcal S}_3)$-module.
\begin{equation}
{\mathcal I}(\xi^1,\xi^2)
=\bigoplus_{\xi\in P_3}W_{3,\xi}\otimes \Hom_{{\mathcal A}(\grs,{\mathcal S}_3)}
(W_{3,\xi},{\mathcal I}(\xi^1,\xi^2)).
\label{eq:st-2}
\end{equation}
By \cite[Theorem 2]{Tanabe}, we have
\begin{align}
\dim_{\C}\Hom_{{\mathcal A}(\grs,{\mathcal S}_3)}
(W_{3,\xi},{\mathcal I}(\xi^1,\xi^2))
\leq
\dim_{\C}I_{(V_L^{\otau})^{\otimes {\ell}}}\binom{X_{\xi,\gamma^3}}
{X_{\xi^1,\gamma^1}\ X_{\xi^2,\gamma^2}}\label{eq:a-inequi}
\end{align}
for $\xi\in P_3$.

Now we compute some fusion rules for $V_{L^{\oplus\ell}}^{\htau}$.

\begin{prop}\label{prop:fusion-L-l}
Let $\lambda,\lambda^1,\lambda^2$ be nonzero elements of $\Klein^{\ell}$ such that
$\lambda^1\not\equiv_{\htau}\lambda^2$, $\gamma,\gamma^1,\gamma^2,\eta\in\Z_3^{\ell}$, $i=1,2$,
and $\varepsilon,\varepsilon_1,\varepsilon_2=0,1,2$.
Then
\begin{align}
V_{L_{({\mathbf 0},\gamma^1)}}(\varepsilon_1)\times V_{L_{({\mathbf 0},\gamma^2)}}(\varepsilon_2)
&=
V_{L_{({\mathbf 0},\gamma^1+\gamma^2)}}(\varepsilon_1+\varepsilon_2)
\label{eq:fusion-L},\\
V_{L_{({\mathbf 0},\gamma^1)}}(\varepsilon)\times V_{L_{(\lambda,\gamma^2)}}
&=
V_{L_{(\lambda,\gamma^1+\gamma^2)}},\label{eq:fusion-Ll-u-0}
\\
V_{L_{(\lambda^1,\gamma^1)}}\times
V_{L_{(\lambda^2,\gamma^2)}}
&=
\sum_{j=0}^{2}V_{L_{(\lambda^1+\htau^j(\lambda^2),\gamma^1+\gamma^2)}},\label{eq:fusion-Ll-u-1}
\\
V_{L_{(\lambda,\gamma^1)}}\times
V_{L_{(\lambda,\gamma^2)}}
&=
\sum_{\rho=0}^{2}V_{L_{({\mathbf 0},\gamma^1+\gamma^2)}}(\rho)+2
V_{L_{(\lambda,\gamma^1+\gamma^2)}},\label{eq:fusion-Ll-u-2}\\
V_{L_{({\mathbf 0},\gamma)}}(\varepsilon_1)\times
\vt{L^{\oplus {\ell}}}{\eta}{\htau^i}{\varepsilon_2}
&=
\vt{L^{\oplus {\ell}}}{\eta-i\gamma}{\htau^i}{i\varepsilon_1+\varepsilon_2},
\label{eq:fusion-Ll-t-0}
\\
V_{L_{(\lambda,\gamma)}}\times
\vt{L^{\oplus {\ell}}}{\eta}{\htau^i}{\varepsilon}
&=
\sum_{\rho=0}^{2}\vt{L^{\oplus {\ell}}}{\eta-i\gamma}{\htau^i}{\rho}.
\label{eq:fusion-Ll-t-1}
\end{align}
\end{prop}

\begin{proof}
We shall show \eqref{eq:fusion-Ll-u-1} and \eqref{eq:fusion-Ll-u-2}. 
We put $\lambda=\lambda^1$ in \eqref{eq:fusion-Ll-u-2} to deal with 
\eqref{eq:fusion-Ll-u-1} and \eqref{eq:fusion-Ll-u-2} simultaneously.
For $\lambda^i,i=1,2$, define
$\xi^i=(\xi^i_j)\in\{0,c\}^{\ell}$ by
\eqref{eq:lambda-xi}.
By (\ref{eq:deco-L-l}), $X_{\xi^i,\gamma^i}$ is a
$(V_L^{\otau})^{\otimes \ell}$-submodule of
$V_{L_{(\lambda^i,\gamma^i)}}$ and 
\begin{equation*}
I_{(V_L^{\otau})^{\otimes \ell}}
\binom{V_{L_{({\mathbf 0},\gamma)}}(r)}
{X_{\xi^1,\gamma^1}\ X_{\xi^2,\gamma^2}}\cong
\bigoplus_{\begin{subarray}{l}
\xi=(\xi_j)\in\Z_3^{\ell}\\
\xi_1+\cdots+\xi_{\ell}=r
\end{subarray}}I_{(V_L^{\otau})^{\otimes \ell}}
\binom{X_{\xi,\gamma}}
{X_{\xi^1,\gamma^1}\ X_{\xi^2,\gamma^2}}
\end{equation*}
for $\gamma\in\Z_3^{\ell}$ and $r\in \Z_3$.
For $\xi\in\Z_3^{\ell}$ and $\gamma\in\Z_3^{\ell}$ such that $\gamma^1+\gamma^2\neq\gamma$,
it follows from Proposition \ref{prop:fusionL} and \cite[Proposition 2.10]{DMZ}
that
\begin{equation*}
I_{(V_L^{\otau})^{\otimes \ell}}
\binom{X_{\xi,\gamma}}
{X_{\xi^1,\gamma^1}\ X_{\xi^2,\gamma^2}}=0.
\end{equation*}

By \cite[Proposition 11.9]{DL}, we obtain
\begin{align*}
&\dim_{\C}I_{V_{L^{\oplus {\ell}}}^{\htau}}\binom{V_{L_{({\mathbf 0},\gamma)}}(r)}
{V_{L_{(\lambda^1,\gamma^1)}}\ V_{L_{(\lambda^2,\gamma^2)}}}\leq
\dim_{\C}I_{(V_L^{\otau})^{\otimes \ell}}
\binom{V_{L_{({\mathbf 0},\gamma)}}(r)}
{X_{\xi^1,\gamma^1}\ X_{\xi^2,\gamma^2}}\\
&=
\sum_{\begin{subarray}{l}
\xi=(\xi_j)\in\Z_3^{\ell}\\
\xi_1+\cdots+\xi_{\ell}=r
\end{subarray}}
\dim_{\C}
I_{(V_L^{\otau})^{\otimes \ell}}
\binom{X_{\xi,\gamma}}
{X_{\xi^1,\gamma^1}\ X_{\xi^2,\gamma^2}}=0.
\end{align*}
For the same reason, we can show easily that
\begin{equation}
I_{V_{L^{\oplus {\ell}}}^{\htau}}\binom{M}
{V_{L_{(\lambda^1,\gamma^1)}}\ V_{L_{(\lambda^2,\gamma^2)}}}=0\label{eq:other}
\end{equation}
for $M\not\in\{V_{L_{(\lambda,\gamma^1+\gamma^2)}}, V_{L_{({\mathbf 0},\gamma^1+\gamma^2)}}(r)\ |\
{\mathbf 0}\neq \lambda\in(\Klein^{\ell})_{\equiv_{\htau}},r=0,1,2\}$.

{From} now on, we use the notation in the preparation just before this proposition.
For example, $\gamma^3=\gamma^1+\gamma^2$, ${\mathcal R}_i=\{h(\lambda^i)\in(\Klein^{\ell})_{\equiv_{\htau}}\ |\ h\in \grs\}$,
and ${\mathcal S}_i=\{V_{L_{(\mu^i,\gamma^i)}}\ |\ \mu^i\in{\mathcal R}_i\}$ for $i=1,2$.
The following symbols are used to describe the fusion rules for $(V_{L}^{\otau})^{\otimes\ell}$:
Set
\begin{equation*}
\Xi(\xi^1,\xi^2)=
\Big\{\xi=(\xi_j)\in\{0,1,2,c\}^{\ell}\ \Big|\
\begin{array}{l}
\xi_j=\xi^1_j+\xi^2_j\\
\mbox{for all }
j\not\in\{k\ |\ \xi^1_k=\xi^2_k=c\}
\end{array}
\Big\}
\end{equation*}
and
\begin{equation*}
\Xi(\xi^1,\xi^2)_k=
\{\xi=(\xi_j)\in\Xi(\xi^1,\xi^2)\ |\
|\{j\ |\ \xi^1_j=\xi^2_j=\xi_j=c\}|=k\}
\end{equation*}
for nonnegative integers $k$.
For example, if $\xi^1=(0,c,1,c)$ and $\xi^2=(1,c,c,2)$ in $\{0,1,2,c\}^{4}$, then
\begin{align*}
\Xi(\xi^1,\xi^2)&=
\{(1,0,c,c),(1,1,c,c),(1,2,c,c),(1,c,c,c)\}
\end{align*}
and $\Xi(\xi^1,\xi^2)_{1}=\{(1,c,c,c)\}$.
Note that
$|\Xi(\xi^1,\xi^2)_k|=
\binom{|\{j\ |\ \xi^1_j=\xi^2_j=c\}|}{k}3^{|\{j\ |\ \xi^1_j=\xi^2_j=c\}|-k}$
and for $\xi\in  \Xi(\xi^1,\xi^2)_k$,
\begin{equation}
\wt_{\tilde{\Klein}}(\xi)=
\wt_{\tilde{\Klein}}(\xi^1)+\wt_{\tilde{\Klein}}(\xi^2)
-2|\{j\ |\ \xi^1_j=\xi^2_j=c\}|+k.\label{eq:wt}
\end{equation}

By Proposition \ref{prop:fusionL} and \cite[Proposition 2.10]{DMZ}, we have
\begin{equation}
\dim_{\C}I_{(V_L^{\otau})^{\otimes {\ell}}}\binom{X_{\xi,\gamma^3}}
{X_{\xi^1,\gamma^1}\ X_{\xi^2,\gamma^2}}
=
\left\{
\begin{array}{ll}
2^k&\mbox{if }\xi\in \Xi(\xi^1,\xi^2)_k,\\
0&\mbox{if }\xi\not\in \Xi(\xi^1,\xi^2)
\end{array}\right.
\label{eq:fusiona}
\end{equation}
for $\xi\in\{0,1,2,c\}^{\ell}$.

By Lemma \ref{lem:fusion-L-l-geq}, we have
\begin{align}
V_{L_{(\mu^1,\gamma^1)}}\times V_{L_{(\mu^2,\gamma^2)}}&\geq
\sum_{j=0}^{2}V_{L_{(\mu^1+\htau^j(\mu^2),\gamma^3)}},\label{eq:fusion-inequiv-1}\\
V_{L_{(\mu,\gamma^1)}}\times V_{L_{(\mu,\gamma^2)}}&\geq
\sum_{\rho=0}^{2}V_{L_{({\mathbf 0},\gamma^3)}}(\rho)+2
V_{L_{(\mu,\gamma^3)}} \label{eq:fusion-inequiv-2}
\end{align}
for $\mu,\mu_1\in{\mathcal R}_1$ and $\mu_2\in{\mathcal R}_2$.
We shall compute the dimension of ${\mathcal I}(\xi^1,\xi^2)$ in
two ways using \eqref{eq:st-1} and \eqref{eq:st-2}.

\medskip
\noindent{\it Case 1.} We deal with the case
$\supp_{\Klein}(\lambda^1)\neq
\supp_{\Klein}(\lambda^2)$. Note that
$g_1(\lambda^1)\not\equiv_{\htau}g_2(\lambda^2)$ for all
$g_1,g_2\in\grs$ and $\wt_{\tilde{\Klein}}(\xi)>0$ for all
$\xi\in\Xi(\xi^1,\xi^2)$. 
We recall that $\dim_{\C} W_{i,\xi^i}=|{\mathcal R}_i|=3^{\wt_{\tilde{\Klein}}(\xi^i)-1}$for $i=1,2$ by \eqref{eq:dim-W}.
By (\ref{eq:dimI}) and
(\ref{eq:fusion-inequiv-1}), we have
\begin{equation}
\dim_{\C}{\mathcal I}(\xi^1,\xi^2)
\geq 3|{\mathcal R}_1|
|{\mathcal R}_2|
= 3^{\wt_{\tilde{\Klein}}(\xi^1)+\wt_{\tilde{\Klein}}(\xi^2)-1}.
\label{eq:ineq-1-1}
\end{equation}

On the other hand, we have 
\begin{align}
\label{eq:ineq-1-2}
&\dim_{\C}{\mathcal I}(\xi^1,\xi^2)\nonumber\\
&=\sum_{\xi\in P_3}\dim_{\C}W_{3,\xi}\dim_{\C}
\Hom_{{\mathcal A}(\grs,{\mathcal S}_3)}
(W_{3,\xi},{\mathcal I}(\xi^1,\xi^2))\nonumber\\
&\leq\sum_{\xi\in\Xi(\xi^1,\xi^2)}
3^{\wt_{\tilde{\Klein}}(\xi)-1}\dim_{\C}I_{(V_L^{\otau})^{\otimes {\ell}}}
\binom{X_{\xi,\gamma^3}}{X_{\xi^1,\gamma^1}\ X_{\xi^2,\gamma^2}}\nonumber\\
&=\sum_{k=0}^{|\{j\ |\ \xi^1_j=\xi^2_j=c\}|}
\sum_{\xi\in\Xi(\xi^1,\xi^2)_k}
3^{\wt_{\tilde{\Klein}}(\xi^1)+\wt_{\tilde{\Klein}}(\xi^2)
-2|\{j\ |\ \xi^1_j=\xi^2_j=c\}|+k-1}2^k\nonumber\\
&=\sum_{k=0}^{|\{j\ |\ \xi^1_j=\xi^2_j=c\}|}
\binom{|\{j\ |\ \xi^1_j=\xi^2_j=c\}|}{k}3^{|\{j\ |\ \xi^1_j=\xi^2_j=c\}|-k}\nonumber\\
&\quad{}\times
3^{\wt_{\tilde{\Klein}}(\xi^1)+\wt_{\tilde{\Klein}}(\xi^2)
-2|\{j\ |\ \xi^1_j=\xi^2_j=c\}|+k-1}2^k\nonumber\\
&=\sum_{k=0}^{|\{j\ |\ \xi^1_j=\xi^2_j=c\}|}
\binom{|\{j\ |\ \xi^1_j=\xi^2_j=c\}|}{k}\nonumber\\
&\quad{} \times
3^{\wt_{\tilde{\Klein}}(\xi^1)+\wt_{\tilde{\Klein}}(\xi^2)
-|\{j\ |\ \xi^1_j=\xi^2_j=c\}|-1}2^k\nonumber\\
&=3^{\wt_{\tilde{\Klein}}(\xi^1)+\wt_{\tilde{\Klein}}(\xi^2)-1}
\end{align}
by \eqref{eq:st-2}, \eqref{eq:a-inequi}, \eqref{eq:wt} and \eqref{eq:fusiona}.
By (\ref{eq:ineq-1-1}) and (\ref{eq:ineq-1-2}), we have
\begin{align*}
\dim_{\C}{\mathcal I}(\xi^1,\xi^2)&=\sum_{M^3\in{\mathcal S}_3}
\sum_{\mu^1\in {\mathcal R}_1,\mu^2\in{\mathcal R}_2}
\dim_{\C}I_{V_{L^{\oplus\ell}}^{\htau}}
\binom{M^3}{V_{L_{(\mu^1,\gamma^1)}}\ V_{L_{(\mu^2,\gamma^2)}}}\\
&= 3^{\wt_{\tilde{\Klein}}(\xi^1)+\wt_{\tilde{\Klein}}(\xi^2)-1},
\end{align*}
and thus it follows from (\ref{eq:other}) that
the equality holds in (\ref{eq:fusion-inequiv-1}).
Setting $\mu^1=\lambda^1$ and $\mu^2=\lambda^2$ in \eqref{eq:fusion-inequiv-1}, we have
\eqref{eq:fusion-Ll-u-1}.

\medskip
\noindent{\it Case 2.} We deal with the case
$\supp_{\Klein}(\lambda^1)=\supp_{\Klein}(\lambda^2)$.
Note that ${\mathcal R}_1={\mathcal R}_2$,
$|\{
(\mu^1,\mu^2)\in{\mathcal R}_1\times{\mathcal R}_1\ |\ \mu^1\not\equiv_{\htau}\mu^2\}|=
|{\mathcal R}_1|(|{\mathcal R}_1|-1)$ and $\xi^1=\xi^2$ in this case. By (\ref{eq:dimI}),
(\ref{eq:fusion-inequiv-1}), and (\ref{eq:fusion-inequiv-2}), we have
\begin{align}
\dim_{\C}{\mathcal I}(\xi^1,\xi^1)
&\geq 5|{\mathcal R}_1|+3|{\mathcal R}_1|(|{\mathcal R}_1|-1)\nonumber\\
&=3^{2\wt_{\tilde{\Klein}}(\xi^1)-1}+
2\cdot 3^{\wt_{\tilde{\Klein}}(\xi^1)-1}.
\label{eq:ineq-2-1}
\end{align}

On the other hand, we have 
\begin{align}
\label{eq:ineq-2-2}
&\dim_{\C}{\mathcal I}(\xi^1,\xi^1)\nonumber\\
&=\sum_{\xi\in P_3}\dim_{\C}W_{3,\xi}\dim_{\C}
\Hom_{{\mathcal A}(\grs,{\mathcal S}_3)}
(W_{3,\xi},{\mathcal I}(\xi^1,\xi^1))\nonumber\\
&\leq\sum_{\xi\in P_3}
3^{\max\{0,\wt_{\tilde{\Klein}}(\xi)-1\}}\dim_{\C}I_{(V_L^{\otau})^{\otimes {\ell}}}
\binom{X_{\xi,\gamma^3}}{X_{\xi^1,\gamma^1}\ X_{\xi^1,\gamma^2}}\nonumber\\
&=
\sum_{
\begin{subarray}{l}
\xi\in\Xi(\xi^1,\xi^1),\\
{\wt_{\tilde{\Klein}}(\xi)\neq 0}
\end{subarray}}
3^{\wt_{\tilde{\Klein}}(\xi)-1}\dim_{\C}I_{(V_L^{\otau})^{\otimes {\ell}}}
\binom{X_{\xi,\gamma^3}}{X_{\xi^1,\gamma^1}\ X_{\xi^1,\gamma^2}}\nonumber\\
&\quad{}+
\sum_{
\begin{subarray}{l}
\xi\in\Xi(\xi^1,\xi^1),\\
\wt_{\tilde{\Klein}}(\xi)= 0
\end{subarray}}
\dim_{\C}I_{(V_L^{\otau})^{\otimes {\ell}}}
\binom{X_{\xi,\gamma^3}}{X_{\xi^1,\gamma^1}\ X_{\xi^1,\gamma^2}}\nonumber\\
&=
\sum_{k=1}^{\wt_{\tilde{\Klein}}(\xi^1)}
\sum_{\xi\in\Xi(\xi^1,\xi^1)_k}
3^{k-1}2^k+3^{\wt_{\tilde{\Klein}}(\xi^1)}\nonumber\\
&=
\sum_{k=1}^{\wt_{\tilde{\Klein}}(\xi^1)}
\binom{\wt_{\tilde{\Klein}}(\xi^1)}{k}3^{\wt_{\tilde{\Klein}}(\xi^1)-k}
3^{k-1}2^k+3^{\wt_{\tilde{\Klein}}(\xi^1)}\nonumber\\
&=
3^{2\wt_{\tilde{\Klein}}(\xi^1)-1}+
2\cdot 3^{\wt_{\tilde{\Klein}}(\xi^1)-1}
\end{align}
by \eqref{eq:st-2}, \eqref{eq:a-inequi}, \eqref{eq:wt}, and \eqref{eq:fusiona}.
By \eqref{eq:ineq-2-1} and \eqref{eq:ineq-2-2}, we have
\begin{align*}
\dim_{\C}{\mathcal I}(\xi^1,\xi^1)&=
\sum_{M^3\in{\mathcal S}_3}
\sum_{\mu^1,\mu^2\in {\mathcal R}_1}
\dim_{\C}I_{V_{L^{\oplus\ell}}^{\htau}}
\binom{M^3}{V_{L_{(\mu^1,\gamma^1)}}\ V_{L_{(\mu^2,\gamma^2)}}}\\
&=
3^{2\wt_{\tilde{\Klein}}(\xi^1)-1}+
2\cdot 3^{\wt_{\tilde{\Klein}}(\xi^1)-1}
\end{align*}
and 
thus it follows from (\ref{eq:other}) that
the equality holds in (\ref{eq:fusion-inequiv-1}) and (\ref{eq:fusion-inequiv-2}).
Setting $\mu^1=\lambda^1, \mu^2=\lambda^2$ in \eqref{eq:fusion-inequiv-1}
and $\mu=\lambda^1=\lambda$ in \eqref{eq:fusion-inequiv-2}, 
we have \eqref{eq:fusion-Ll-u-1} and
\eqref{eq:fusion-Ll-u-2}.

The same argument as above shows \eqref{eq:fusion-Ll-t-1}. We shall sketch the proof 
of \eqref{eq:fusion-Ll-t-1} for $i=1$.
By Proposition \ref{prop:fusionL}, we can show easily that
\begin{equation*}
I_{V_{L^{\oplus {\ell}}}^{\htau}}\binom{M}
{V_{L_{(\lambda,\gamma)}}\ \vt{L^{\oplus {\ell}}}{\eta}{\htau}{\varepsilon}}=0
\end{equation*}
for $M\not\in\{
\vt{L^{\oplus {\ell}}}{\eta-\gamma}{\htau}{\rho}\ |\ \rho=0,1,2\}$.
Take $\grs$-stable sets ${\mathcal S}^{T}_{2}=\{\vt{L^{\oplus {\ell}}}{\eta}{\htau}{\varepsilon}\}$
and ${\mathcal S}^{T}_{3}=\{\vt{L^{\oplus {\ell}}}{\eta-\gamma}{\htau}{\rho}\ |\ \rho=0,1,2\}$
and set ${\mathcal M}^{T}_{i}=
\oplus_{M\in{\mathcal S}^{T}_{i}}M, i=2,3$.
Note that for $M\in{\mathcal S}^{T}_{i}$, $(\grs)_{M}=\{g\in\grs\ |\ g(M)=M\}$
equals $\grs$.
For $\xi\in\{0,1,2\}^{\ell}$,
define an action of $D(M),M\in {\mathcal S}^{T}_{i}$ on $\C e(\xi)$ as follows:
For $g=(\otau^{i_1},\ldots,\otau^{i_{\ell-1}},1)\in\grs$, set
\begin{align*}
g\otimes e(M)\cdot e(\xi)= &
\zeta_3^{\langle(i_1,\ldots,i_{\ell-1},0),\xi\rangle_{\Z_3}}
e(\xi).
\end{align*}
Denote the $D(M)$-module $\C e(\xi)$ by $W^{T}_{i,\xi}$.
Set
\begin{align*}
&P(\vt{L^{\oplus {\ell}}}{\eta}{\htau}{\rho})=
P(\vt{L^{\oplus {\ell}}}{\eta-\gamma}{\htau}{\rho})\\
&=
\{\xi=(\xi_k)\in\{0,1,2\}^{\ell}\ |\ \sum_{k=1}^{\ell}\xi_k\equiv \rho\pmod{3}\}
\end{align*}
for $\rho=0,1,2$. Note that
\begin{equation*}
\{0,1,2\}^{\ell}=\bigcup_{\rho=0}^{2}P(\vt{L^{\oplus {\ell}}}{\eta-\gamma}{\htau}{\rho});\quad \mbox{disjoint.}
\end{equation*}
Set $P^{T}_{2}=P(\vt{L^{\oplus {\ell}}}{\eta}{\htau}{\varepsilon})$
and $P^{T}_{3}=\{0,1,2\}^{\ell}$.
Then
$\{W^{T}_{i,\xi}\ |\ \xi\in P^{T}_{i}\}$
is a complete list of irreducible
${\mathcal A}(\grs,{\mathcal S}^{T}_{i})$-modules by
Theorem \ref{thm:AGS} for $i=2,3$.
For $\xi=(\xi_k),\gamma=(\gamma_k)\in\{0,1,2\}^{\ell}$, set
\begin{equation*}
X^{T}_{\xi,\gamma}=\bigotimes_{k=1}^{\ell}\vt{L}{\gamma_k}{\otau}{\xi_k}.
\end{equation*}
For the same reason as in the proof of (\ref{eq:schur-weyl}),
we have
\begin{align*}
{\mathcal M}^{T}_2=\bigoplus_{\xi\in P^{T}_{2}}W^{T}_{2,\xi}\otimes
X^{T}_{\xi,\eta},\qquad
{\mathcal M}^{T}_3=\bigoplus_{\xi\in P^{T}_{3}}W^{T}_{3,\xi}\otimes
X^{T}_{\xi,\eta-\gamma}
\end{align*}
as an ${\mathcal A}(\grs,{\mathcal S}^{T}_2)$- and ${\mathcal
A}(\grs,{\mathcal S}^{T}_3)$-module, respectively. Set
\begin{align*}
{\mathcal I}^{T}&=\bigoplus_{(M^1,M^2,M^3) \in{\mathcal
S}_1\times{\mathcal S}^{T}_2\times{\mathcal S}^{T}_3}
I_{V_{L^{\oplus {\ell}}}^{\htau}}\binom{M^3}{M^1\ M^2}
\otimes_{\C}M^1\otimes_{\C} M^2.
\end{align*}
Let $\xi^1\in P_1=P(V_{L_{(\lambda,\gamma)}})$ and $\xi^2\in
P^{T}_{2}$. Fix nonzero elements $v^{10}\in X_{\xi^{1},\gamma}$
and $v^{T,20}\in X^{T}_{\xi^{2},\eta}$. Set
\begin{align*}
&{\mathcal I}^{T}(\xi^1,\xi^2)\\
&= \spn_{\C}\{f\otimes w^1\otimes v^{10}\otimes w^2 \otimes
v^{T,20}\in{\mathcal I}^{T}\ |\ w^1\in W_{1,\xi^1}, w^2\in
W^{T}_{2,\xi^2}\}.
\end{align*}
Applying the same arguments as in the case of
\eqref{eq:fusion-Ll-u-1} and \eqref{eq:fusion-Ll-u-2},
we have
\begin{equation*}
\dim_{\C}{\mathcal I}^{T}(\xi^1,\xi^2)=3^{\wt_{\tilde{\Klein}}(\xi^1)}.
\end{equation*}
Therefore, \eqref{eq:fusion-Ll-t-1} holds.

The other formulas can be proved similarly.

\end{proof}

\section{Modules of $V_{L_{\zero\times D}}^{\htau}$}
\label{section:vl0dt}

Let $D$ be a self-orthogonal $\Z_3$-code of length $\ell$.
In this section we discuss $V_{L_{\zero\times D}}^{\htau}$-modules. Note that $V_{L_{\zero\times D}}^{\htau}
= \oplus_{\gamma\in D}V_{L_{({\mathbf 0},\gamma)}}(0)$
as  $V_{L^{\oplus {\ell}}}^{\htau}$-modules. Let
$\gamma^{(1)},\ldots,\gamma^{(\ell)}$ be a basis of $\Z_3^{\ell}$
such that $\gamma^{(1)},\ldots,\gamma^{(d)}$ form a basis of $D$.

For $j=1,\ldots,\ell$, define a linear transformation
$\chi_{j}$ on $V_{(L^{\perp})^{\oplus \ell}}=
\oplus_{\delta\in\Z_3^{\ell}}
V_{L_{\Klein^{\ell}\times \delta}}$ by
$\chi_{j}(u)=\zeta_3^{p_j}u$ for $\delta=\sum_{k=1}^{\ell}p_k\gamma^{(k)}\in\Z_3^{\ell}$ and
$u\in V_{L_{\Klein^{\ell}\times \delta}}$.
The restriction of $\chi_j$ to $V_{L_{\zero\times D}}$ is an automorphism of
$V_{L_{\zero\times D}}$ for $j=1,\ldots,\ell$.
Let $\Phi_D$ be the automorphism group of $V_{L_{\zero\times D}}$
generated by $\chi_1,\ldots,\chi_d$.
Since $\htau$ commutes with $\Phi_D$,
$\Phi_D$ induces an automorphism group of $V_{L_{\zero\times D}}^{\htau}$.
Note that $(V_{L_{\zero\times D}}^{\htau})^{\Phi_D}=V_{L_{({\mathbf 0},{\mathbf 0})}}(0)=
V_{L^{\oplus {\ell}}}^{\htau}$.

For $j=1,\ldots,d$,
$\lambda\in\Klein^{\ell}$, and $\gamma\in\Z_3^{\ell}$,
$V_{L_{\lambda\times(\gamma+D)}}$ is $\chi_j$-invariant and
\begin{equation}
V_{L_{\lambda\times (\gamma+D)}}
=\bigoplus_{\delta\in D}V_{L_{(\lambda,\gamma+\delta)}}
\label{eq:st-D-1}
\end{equation}
is an eigenspace decomposition for $\Phi_D$.  We also have
\begin{align}
V_{L_{\zero\times (\gamma+D)}}(\varepsilon)
&=\bigoplus_{\delta\in D}V_{L_{(\zero,\gamma+\delta)}}(\varepsilon),\quad
\varepsilon\in\Z_3, \gamma\in \Z_3^{\ell},\nonumber\\
V_{L_{\lambda\times (\gamma+D)}}
&=\bigoplus_{\delta\in D}V_{L_{(\lambda,\gamma+\delta)}},\quad
\lambda\in\Klein^{\ell}, \gamma\in \Z_3^{\ell}
\label{eq:st-D-2}
\end{align}
as $V_{L^{\oplus \ell}}^{\htau}$-modules.

For $j=1,\ldots,d$,
$\lambda\in\Klein^{\ell},\gamma\in D^{\perp},
u\in V_{L_{\zero\times D}}$, and $v\in V_{L_{\lambda\times(\gamma+D)}}$,
we have
\begin{equation*}
Y_{V_{L_{\lambda\times(\gamma+D)}}}(\chi_{j}u,x)\chi_{j}v=
\chi_{j}(Y_{V_{L_{\lambda\times(\gamma+D)}}}(u,x)v).
\end{equation*}
Hence
$V_{L_{\lambda\times(\gamma+D)}}\circ \chi_j\cong \chi_j^{-1}
\big(V_{L_{\lambda\times(\gamma+D)}})=V_{L_{\lambda\times(\gamma+D)}}$
as $V_{L_{\zero\times D}}$-modules
and
\begin{align*}
V_{L_{\zero\times(\gamma+D)}}(\varepsilon)\circ \chi_j
&\cong \chi_j^{-1}\big(V_{L_{\zero\times(\gamma+D)}}(\varepsilon)\big)=
V_{L_{\zero\times(\gamma+D)}}(\varepsilon),\quad \gamma\in D^{\perp},\varepsilon\in\Z_3,\\
V_{L_{\lambda\times(\gamma+D)}}\circ \chi_j
&\cong \chi_j^{-1}\big(V_{L_{\lambda\times(\gamma+D)}}\big)
=V_{L_{\lambda\times(\gamma+D)}},
\quad {\mathbf 0}\neq \lambda\in\Klein^{\ell}, \gamma\in D^{\perp}
\end{align*}
as $V_{L_{\zero\times D}}^{\htau}$-modules.

It follows from \eqref{eq:twist-D-L} and the corresponding formula for $\htau^2$-twisted modules
that for $\eta\in D^{\perp}$, $i=1,2$, and $r\in\Z_3$,
\begin{equation}
\vt{L_{\zero\times D}}{\eta}{\htau^i}{r}
\cong\bigoplus_{\gamma\in D}
\vt{L^{\oplus \ell}}{\eta-i\gamma}{\htau^i}{r}
\label{eq:st-D-3}
\end{equation}
as $V_{L^{\oplus \ell}}^{\htau}$-modules.
Using \eqref{eq:st-D-3}, we
define an action of $\chi_{j}$ on $\vt{L_{\zero\times D}}{\eta}{\htau^i}{r}$ for $j=1,\ldots,\ell$
by setting $\chi_{j}(v)=\zeta_3^{-ip_j}v$ for $\delta=\sum_{k=1}^{\ell}p_k\gamma^{(k)}\in\Z_3^{\ell}$ and 
$v\in \vt{L^{\oplus \ell}}{\delta}{\htau^i}{r}$
and extending $\chi_j$ for arbitrary $v\in\vt{L_{\zero\times D}}{\eta}{\htau^i}{r}$ by
\eqref{eq:st-D-3} and linearity.

By Proposition \ref{prop:fusion-L-l}, we have
\begin{equation*}
Y_{\vt{L_{\zero\times D}}{\eta}{\htau^i}{r}}(\chi_{j}u,x)\chi_{j}v=
\chi_{j}(Y_{\vt{L_{\zero\times D}}{\eta}{\htau^i}{r}}(u,x)v)
\end{equation*}
for $u\in V_{L_{\zero\times D}}$ and
$v\in \vt{L_{\zero\times D}}{\eta}{\htau^i}{r}$. Hence 
\begin{equation*}
\vt{L_{\zero\times D}}{\eta}{\htau^i}{r}\circ \chi_j\cong
\chi_{j}^{-1}\big(\vt{L_{\zero\times D}}{\eta}{\htau^i}{r}\big)=
\vt{L_{\zero\times D}}{\eta}{\htau^i}{r}
\end{equation*}
and so we can define an action of $\Phi_D$
on $\vt{L_{\zero\times D}}{\eta}{\htau^i}{r}$. Then it is clear
that (\ref{eq:st-D-3}) is also an eigenspace decomposition of
$\vt{L_{\zero\times D}}{\eta}{\htau^i}{r}$ for $\Phi_D$.

\begin{lem}\label{lem:res-fusion-D}
Let $N$ be an $\N$-graded weak $V_{L_{\zero\times D}}^{\htau}$-module and let $M$ be an irreducible $V_{L^{\oplus{\ell}}}^{\htau}$-submodule of $N$. 
If $M$ is isomorphic to
$V_{L_{(0,\gamma)}}(\varepsilon),\gamma\in\Z_3^{\ell},\varepsilon\in\Z_3$
or $V_{L_{(\lambda,\gamma)}},{\mathbf 0}\neq \lambda\in\Klein^{\ell},
\gamma\in\Z_3^{\ell}$, then $\gamma\in D^{\perp}$. 
If $M$ is isomorphic to $\vt{L^{\oplus \ell}}{\eta}{\htau^i}{\varepsilon},i=1,2,\varepsilon\in\Z_3, \eta\in\Z_3^{\ell}$, then $\eta\in
D^{\perp}$.
\end{lem}
\begin{proof}
Let $\omega_L$ be the Virasoro element of $V_{L}^{\otau}$. For
$i=1,2$ and $j,k,\varepsilon,\in\{0,1,2\}$, let $(W^1,W^2)$ be one
of $(V_{L^{(0,j)}}(\varepsilon),V_{L^{(0,j+k)}}(\varepsilon))$,
$(V_{L^{(c,j)}},V_{L^{(c,j+k)}})$, or
$(\vt{L}{j}{\otau^i}{\varepsilon},
\vt{L}{j+k}{\otau^i}{\varepsilon})$. Let $\lambda_s$ be the
eigenvalue of $(\omega_L)_1$ on the top level of $W^s$ for
$s=1,2$. Note that
\begin{align*}
V_{L^{(0,-k)}}(0)\times W^1 &= W^2 \quad\text{if }(W^1,W^2)=
(\vt{L}{j}{\otau}{\varepsilon}, \vt{L}{j+k}{\otau}{\varepsilon}),\\
V_{L^{(0,k)}}(0)\times W^1 &= W^2 \quad\text{otherwise}
\end{align*}
by Proposition \ref{prop:fusionL} and that
$\lambda_2-\lambda_1\equiv (jk+2k^2)/3\pmod{\Z}$ by
(\ref{eq:eigen}). We have already obtained a decomposition of every
irreducible $V_{L^{\oplus\ell}}^{\htau}$-module as a
$(V_L^{\otau})^{\otimes \ell}$-module in Theorem
\ref{thm:twist-L-st} and (\ref{eq:deco-L-l}).

Now the proof is similar to that of Lemma \ref{lem:res-fusion}
since
$V_{L^{(0,\gamma_1)}}(0)\otimes\cdots\otimes V_{L^{(0,\gamma_{\ell})}}(0)
\subset V_{L_{\zero\times D}}^{\htau}$ for $\gamma=(\gamma_s)_{s=1}^{\ell}\in D$
and $D$ is self-orthogonal.
\end{proof}

Using the same arguments as in the proofs of Lemma \ref{lem:generate}
and Proposition \ref{prop:irr-L-l},
we can show the following theorem. Indeed, we argue for
$V_{L_{\zero\times D}}^{\htau},
V_{L^{\oplus \ell}}^{\htau}$, and $\Phi_D$ in place of
$V_{L^{\oplus \ell}}^{\htau},
(V_{L}^{\otau})^{\otimes \ell}$, and $\grs$ in Section 5,
respectively.

\begin{thm}\label{thm:irr-D}
Let $D$ be a self-orthogonal $\Z_3$-code of length $\ell$.
Then $V_{L_{\zero\times D}}^{\htau}$ is a simple, rational,
$C_2$-cofinite, and CFT type vertex operator algebra.
Let $D^{\perp}/D=\bigcup_{j=1}^{m}(\rho^j+D)$ be a coset decomposition.
The following is a complete set of representatives
of equivalence classes of irreducible $V_{L_{\zero\times D}}^{\htau}$-modules.

$(1)$ ${\displaystyle V_{L_{\zero\times (\rho^j+D)}}(\varepsilon)},\
j=1,\ldots,m,\ \varepsilon=0,1,2$.

$(2)$ $V_{L_{\lambda\times (\rho^j+D)}},\
{\mathbf 0}\neq \lambda\in (\Klein^{\ell})_{\equiv_{\htau}},\
j=1,\ldots,m$.

$(3)$ ${\displaystyle
\vt{L_{\zero\times D}}{\rho^j}{\htau^i}{\varepsilon}}$, $i=1,2$, $j=1,\ldots,m$, $\varepsilon=0,1,2$.

\end{thm}

We compute some fusion rules for $V_{L_{\zero \times D}}^{\htau}$ 
which will be used in Section 7.2.
\begin{prop}\label{prop:fusion-D}
Let $\lambda,\lambda^1,\lambda^2$ be nonzero elements of $\Klein^{\ell}$ such that
$\lambda^1\not\equiv_{\htau}\lambda^2$, $\gamma,\gamma^1,\gamma^2,\eta\in D^{\perp}$,
$i=1,2$, and $\varepsilon,\varepsilon_1,\varepsilon_2=0,1,2$.
Then
\begin{align}
V_{L_{\zero\times (\gamma^1+D)}}(\varepsilon_1)\times
V_{L_{\zero\times (\gamma^2+D)}}(\varepsilon_2)&=
V_{L_{\zero\times (\gamma^1+\gamma^2+D)}}(\varepsilon_1+\varepsilon_2),\\
V_{L_{\zero\times (\gamma^1+D)}}(\varepsilon)\times
V_{L_{\lambda\times(\gamma^2+D)}}&=
V_{L_{\lambda\times(\gamma^1+\gamma^2+D)}},\\
V_{L_{\lambda^1\times(\gamma^1+D)}}\times V_{L_{\lambda^2\times(\gamma^2+D)}}&=
\sum_{j=0}^{2}V_{L_{(\lambda^1+\htau^j(\lambda^2))\times (\gamma^1+\gamma^2+D)}},
\label{eq:fusionD-1}\\
V_{L_{\lambda\times(\gamma^1+D)}}\times V_{L_{\lambda\times(\gamma^2+D)}}&=
\sum_{\rho=0}^{2}V_{L_{\zero\times(\gamma^1+\gamma^2+D)}}(\rho)+
2V_{L_{\lambda\times (\gamma^1+\gamma^2+D)}},\\
V_{L_{\zero\times (\gamma+D)}}(\varepsilon_1)\times
\vt{L_{\zero\times D}}{\eta}{\htau^i}{\varepsilon_2}&=
\vt{L_{\zero\times D}}{\eta-i\gamma}{\htau^i}{i\varepsilon_1+\varepsilon_2},\\
V_{L_{\lambda\times(\gamma+D)}}\times
\vt{L_{\zero\times D}}{\eta}{\htau^i}{\varepsilon}&=
\sum_{\rho=0}^{2}\vt{L_{\zero\times D}}{\eta-i\gamma}{\htau^i}{\rho}.
\end{align}
\end{prop}
\begin{proof}
We shall show (\ref{eq:fusionD-1}).
Restricting intertwining operators for $V_{L_{\zero\times D}}$ in Lemma \ref{lem:fusion-lattice} to 
$V_{L_{\zero\times D}}^{\htau}$-modules,
we have
\begin{equation}
V_{L_{\lambda^1\times(\gamma^1+D)}}\times V_{L_{\lambda^2\times(\gamma^2+D)}}\geq
\sum_{j=0}^{2}V_{L_{(\lambda^1+\htau^j(\lambda^2))\times (\gamma^1+\gamma^2+D)}}.
\label{eq:fusionD-1-1}
\end{equation}

For $k=1,2,r\in\Z_3,\zero\neq \lambda^3\in\Klein^{\ell}$, and
$\gamma^3\in\Z_3^{\ell}$,
\begin{align*}
I_{V_{L^{\oplus \ell}}^{\htau}}
\binom{V_{L_{\zero\times (\gamma^3+D)}}(r)}
{V_{L_{(\lambda^1,\gamma^1)}}\ V_{L_{(\lambda^2,\gamma^2)}}}
&\cong
\bigoplus_{\delta\in D}
I_{V_{L^{\oplus \ell}}^{\htau}}
\binom{V_{L_{(\zero,\gamma^3+\delta)}}(r)}
{V_{L_{(\lambda^1,\gamma^1)}}\ V_{L_{(\lambda^2,\gamma^2)}}},\\
I_{V_{L^{\oplus \ell}}^{\htau}}
\binom{V_{L_{\lambda^3\times (\gamma^3+D)}}}
{V_{L_{(\lambda^1,\gamma^1)}}\ V_{L_{(\lambda^2,\gamma^2)}}}
&\cong
\bigoplus_{\delta\in D}
I_{V_{L^{\oplus \ell}}^{\htau}}
\binom{V_{L_{(\lambda^3,\gamma^3+\delta)}}}
{V_{L_{(\lambda^1,\gamma^1)}}\ V_{L_{(\lambda^2,\gamma^2)}}},
\\
I_{V_{L^{\oplus \ell}}^{\htau}}
\binom{\vt{L_{\zero\times D}}{\gamma^3}{\htau^k}{r}}
{V_{L_{(\lambda^1,\gamma^1)}}\ V_{L_{(\lambda^2,\gamma^2)}}}
&\cong
\bigoplus_{\delta\in D}
I_{V_{L^{\oplus \ell}}^{\htau}}\binom{
\vt{L^{\oplus \ell}}{\gamma^3-\delta}{\htau^k}{r}
}
{V_{L_{(\lambda^1,\gamma^1)}}\ V_{L_{(\lambda^2,\gamma^2)}}}
\end{align*}
as vector spaces by \eqref{eq:st-D-2} and \eqref{eq:st-D-3}.
By \cite[Proposition 11.9]{DL} and Proposition \ref{prop:fusion-L-l},
\begin{align}
&\dim_{\C}I_{V_{L_{\zero\times D}}^{\htau}}
\binom{V_{L_{(\lambda^1+\htau^j(\lambda^2))\times (\gamma^1+\gamma^2+D)}}}
{V_{L_{\lambda^1\times(\gamma^1+D)}}\ V_{L_{\lambda^2\times(\gamma^2+D)}}}
\leq
\dim_{\C}I_{V_{L^{\oplus \ell}}^{\htau}}
\binom{V_{L_{(\lambda^1+\htau^j(\lambda^2))\times (\gamma^1+\gamma^2+D)}}}
{V_{L_{(\lambda^1,\gamma^1)}}\ V_{L_{(\lambda^2,\gamma^2)}}}\nonumber\\
&=
\sum_{\delta\in D}
\dim_{\C}I_{V_{L^{\oplus \ell}}^{\htau}}
\binom{V_{L_{(\lambda^1+\htau^j(\lambda^2),\gamma^1+\gamma^2+\delta)}}}
{V_{L_{(\lambda^1,\gamma^1)}}\ V_{L_{(\lambda^2,\gamma^2)}}}=1
\label{eq:fusionD-1-2}
\end{align}
for $j=0,1,2$ and
\begin{align}
\dim_{\C}I_{V_{L_{\zero\times D}}^{\htau}}
\binom{W}
{V_{L_{\lambda^1\times(\gamma^1+D)}}\ V_{L_{\lambda^2\times(\gamma^2+D)}}}
\leq
\dim_{\C}I_{V_{L^{\oplus \ell}}^{\htau}}
\binom{W}
{V_{L_{(\lambda^1,\gamma^1)}}\ V_{L_{(\lambda^2,\gamma^2)}}}=0
\label{eq:fusionD-1-3}
\end{align}
for any irreducible $V_{L_{\zero\times D}}^{\htau}$-module
$W\not\cong V_{L_{(\lambda^1+\htau^j(\lambda^2)) \times
(\gamma^1+\gamma^2+D)}}, j=0,1,2$.
By (\ref{eq:fusionD-1-1})--(\ref{eq:fusionD-1-3}),
we obtain (\ref{eq:fusionD-1}).

The other formulas can be proved similarly.
\end{proof}

\section{Modules of $V_{L_{C\times D}}^{\htau}$}\label{sec:main}

In this section we shall study $V_{\LCD}^\htau$-modules
for an arbitrary $\htau$-invariant
self-dual $\Klein$-code $C$ with minimum weight at least
$4$ and an arbitrary self-dual $\Z_3$-code $D$.  

Let $N$ be an $\N$-graded weak $V_{\LCD}^\htau$-module.
Since $N$ is a $V_{L^{\oplus\ell}}^{\htau}$-module,
$N$ is a direct sum of irreducible $V_{L^{\oplus\ell}}^{\htau}$-modules listed in 
Proposition \ref{prop:irr-L-l}.
If $N$ contains an irreducible $V_{L^{\oplus\ell}}^{\htau}$-module
which is isomorphic to $V_{L_{(\lambda,\gamma)}}$ for a nonzero $\lambda\in\Klein^{\ell}$
and $\gamma\in\Z_3^{\ell}$, then Theorem \ref{thm:irr-D} implies that
$N$ also contains an irreducible $V_{L^{\oplus\ell}}^{\htau}$-module
which is isomorphic to $V_{L_{(\lambda,\zero)}}$
since $N$ is a $V_{L_{\zero\times D}}^{\htau}$-module and $D$ is self-dual.
This observation is important in the proof of Proposition \ref{prop:CD-untwist}.
Thus, it is necessary to assume $D$ is self-dual.

Recall that for $\mu\in \Klein^{\ell}$, $C(\mu)$ is the $\Klein$-code generated by $\mu$ and $\htau(\mu)$ 
(cf. Section \ref{section:VLt}).
If $\mu\in \Klein^{\ell}$ has positive even weight, then
\begin{equation}
\label{eq:cmu-dec}
V_{L_{C(\mu)\times \zero}}^{\htau}
\cong
V_{L^{\oplus \ell}}^{\htau}\oplus V_{L_{(\mu,\zero)}}
\end{equation}
as $V_{L^{\oplus \ell}}^{\htau}$-modules.
Since
$N$ is also a $V_{L_{C(\mu)\times\zero}}^{\htau}$-module for each $\mu\in C$,
using \cite[Theorem 2.1.2]{Z},
\eqref{eq:cmu-dec}, and the fusion rules for $V_{L^{\oplus\ell}}^{\htau}$ in Proposition \ref{prop:fusion-L-l},
we can obtain information about irreducible $V_{L^{\oplus\ell}}^{\htau}$-modules contained in $N$
(See Proposition \ref{prop:orth} and Proposition \ref{prop:CD-untwist} below).
Thus, we first study $\N$-graded weak $V_{L_{C(\mu)\times\zero}}^{\htau}$-modules with some
conditions for $\mu\in\Klein^{\ell}$ of positive even weight in Section \ref{subsection:muzero}.
Next, we shall classify the
irreducible $V_{\LCD}^\htau$-modules and establish the rationality
of $V_{\LCD}^\htau$ in Theorem \ref{thm:main} 
in Section \ref{subsection:vlcd}.

\subsection{Properties of $V_{L_{C(\mu)\times\zero}}^{\htau}$-modules}\label{subsection:muzero}
Throughout this subsection, $\ell$ and $\tll$ are fixed even positive integers
with $2\leq \tll\leq \ell$.
In this subsection we study
$\N$-graded weak $V_{L_{C(\mu)\times\zero}}^{\htau}$-modules with some
conditions for $\mu\in\Klein^{\ell}$ of positive even weight.
We deal with the case $\mu=(c^{\tll}0^{\ell-\tll})=(c,\ldots,c,0,\ldots,0)$ until Lemma \ref{lem:res-orth}.
We have
\begin{equation}
\label{eq:cm-dec}
V_{L_{C(c^{\tll}0^{\ell-\tll})\times \zero}}^{\htau}
\cong
V_{L^{\oplus \ell}}^{\htau}\oplus V_{L_{((c^{\tll}0^{\ell-\tll}),\zero)}}
\end{equation}
as $V_{L^{\oplus \ell}}^{\htau}$-modules.
We shall fix the following notation.
Let $W^0=\oplus_{i=0}^{\infty}W^0(i)$ be an $\N$-graded weak
$V_{L_{C(c^{\tll}0^{\ell-\tll})\times\zero}}^{\htau}$-module. Let
$M^0=\oplus_{i=0}^{\infty}M^0(i)$ be a $V_{L^{\oplus{\ell}}}^{\htau}$-submodule of $W^0$ such that $M^0(0)\subset
W^0(0)$. Assume that $M^0$ is isomorphic to
$V_{L_{(\lweight,{\mathbf 0})}}$ for some nonzero
$\lweight=(\lweight_1,\ldots,\lweight_{\ell})\in\Klein^{\ell}$.

In Lemma \ref{lem:binom} we shall describe the action of $o(u\circ v)$ 
(see \cite[Definition 2.1.1]{Z})
on the top level of $M^0$
for some elements $u,v\in V_{L_{((c^{\tll}0^{\ell-\tll}),\zero)}}
\subset V_{L_{C(c^{\tll}0^{\ell-\tll})\times \zero}}^{\htau}$
in the decomposition \eqref{eq:cm-dec}.
Since $o(u\circ v)w=0$ for all elements $w$ in the top level of $M^0$,
we shall obtain the relations \eqref{eq:binom-1} and \eqref{eq:binom-2} below,
which play an important role to get information about $\Delta$ in Lemma \ref{lem:res-orth}.
With the help of the action of $G_{\ell}$, 
Lemma \ref{lem:res-orth} immediately induces Proposition \ref{prop:orth}.

For $S\subset \{1,\ldots,\ell\}$, set 
\begin{align}
\label{eq:sstar}
\os
&=\{i\in\{1,\ldots,m\} |\ i\not\in S\}.
\end{align} 
Recall that for each $x \in \K$ we assign $\be(x) \in
L^\perp$ by $\be(0) = 0$, $\be(a) = \be_2/2$, $\be(b) = \be_0/2$,
and $\be(c) = \be_1/2$. 
For $j\in\Klein$, we use $\beta^{(s)}(j)$ to denote the element $\beta(j) \in L^{\perp}$
in the $s$-th entry of $(L^{\perp})^{\oplus\ell}$.
For $p=(p_i)\in\Klein^{\ell}$ and $\epsilon=(\epsilon_i)\in\{1,-1\}^{\ell}$, set
$\beta(p;\epsilon)=\sum_{i=1}^{\ell}\epsilon_i\beta^{(i)}(p_i)$. 
For example,
\begin{align*}
\sum_{i=1}^{\tll}\epsilon_i\beta_{1}^{(i)}/2
&=\sum_{i=1}^{\tll}\epsilon_i\beta^{(i)}(c)
=\beta((c^{\tll}0^{\ell-\tll});\epsilon).
\end{align*}
We simply write  $\beta(p)$ for $\beta(p;(1,\ldots,1))$.
For $\alpha\in (L^{\perp})^{\oplus\ell}$, set
$\e(\alpha)=e^{\alpha}$.
Set 
\begin{align*}
S_j(\lweight)&
=\{i\in\{1,\ldots,m\}\ |\ \lweight_i=j\}
\end{align*}
for $j=a,b,c$.

For a formal Laurent series $p(x)=\sum_{n\in\Z}p_nx^n$ in one
variable $x$ and $i\in\Z$, set
$p(x)|_{x^i}=p_i$. For $p(x,\anx)=\sum_{n,m\in\Z}p_{nm}x^n\anx^m$ and
$i,j\in\Z$, set $p(x,\anx)|_{x^i, \anx^j}=p_{ij}$ similarly.
For homogeneous $u,v\in V$, we shall use the following expression:
\begin{align}
o(u\circ v)
&=\sum_{r=0}^{\wt u}\binom{\wt u}{r}Y_{M}
(Y(u,x)v,\anx)w\big|_{x^{1-r},\anx^{-\wt u-\wt v+r-1}}.\label{eq:ouv}
\end{align}

The following lemma is the key result in this section.

\begin{lem}\label{lem:binom}
$(1)$ 
Let $S$ be a subset of $\{1,\ldots,\tll\}$ with $1\leq |S|\leq {\tll}/2$.
If $\lweight\not\equiv_{\htau}(c^{\tll}0^{\ell-\tll})$ or 
$|S|\leq {\tll}/2-2$, then for $j=0,1,2$ 
and $\epsilon=(\epsilon_i)\in\{1,-1\}^{\ell}$,
we have
\begin{equation}
\delta_{\langle\sum_{i\in S}\epsilon_i\beta_{j}^{(i)}/2,\beta(\lweight)\rangle,-|S|}
\binom{\langle\sum_{i=1}^{\tll}\epsilon_i\beta_{j}^{(i)}/2,\beta(\lweight)\rangle+{\tll}/2}{{\tll}-2|S|+1}
=0.
\label{eq:binom-1}
\end{equation}

$(2)$ Suppose $\lweight\not\equiv_{\htau}(c^{\tll}0^{\ell-\tll})$ or ${\tll}\geq 4$.
Then, for $\epsilon=(\epsilon_i)\in\{1,-1\}^{\ell}$ we have
\begin{equation}
\sum_{j=0}^{2}
\binom{\langle\sum_{i=1}^{\tll}\epsilon_i\beta_{j}^{(i)}/2,\beta(\lweight)\rangle+{\tll}/2}{{\tll}+1}
=0.
\label{eq:binom-2}
\end{equation}
\end{lem}
\begin{proof}
Let $S$ be a subset of $\{1,\ldots,\tll\}$ and set $s=|S|$. 
Let 
\begin{align*}
\mathbf{u}&=\sum_{j=0}^{2}\e(\sum_{i=1}^{\tll}\frac{\beta_{j}^{(i)}}{2}),&
\mathbf{v}&=\sum_{j=0}^{2}
\e(\sum_{i\in S}\frac{\beta_{j}^{(i)}}{2}+
\sum_{i\in \os}\frac{-\beta_{j}^{(i)}}{2}).
\end{align*}
Then
\begin{align*}
\mathbf{u}
&=
\sum_{j=0}^{2}\htau^{j}\e(\sum_{i=1}^{\tll}\frac{\beta_{1}^{(i)}}{2}), &
\mathbf{v}
&=
\sum_{j=0}^{2}
\htau^{j}\e(\sum_{i\in S}\frac{\beta_{1}^{(i)}}{2}+
\sum_{i\in \os}\frac{-\beta_{1}^{(i)}}{2})
\end{align*}
by \eqref{eq:def-tau} and hence 
$\mathbf{u}$ and $\mathbf{v}$ are elements of
$V_{L_{C(c^{\tll}0^{\ell-\tll})\times \zero}}^{\htau}$ of weight $m/2$.
We shall describe the action of
\begin{align}
o(\mathbf{u}\circ \mathbf{v})
&=\sum_{r=0}^{\tll/2}\binom{\tll/2}{r}Y_{W^0}(
Y(\mathbf{u},x)\mathbf{v},\anx)\big|_{x^{1-r},\anx^{-\tll+r-1}}\label{eq:oee}
\end{align} 
on
the top level of $M^0$ (cf. \eqref{eq:ouv}). 
For $j=0,1,2$, set
\begin{align}
\Omega_{1j}&=
\zeta_{36}^{9\tll+18s}x^{-\tll+2s}\exp\big(\sum_{k=1}^{\infty}\frac{(\sum_{i=1}^{\tll}\beta_{j}^{(i)}/2)(-k)}{k}x^k\big)
\e(\sum_{i\in S}\beta_{j}^{(i)}),\nonumber\\
\Omega_{2j}&=\zeta_{36}^{9\tll+18s}x^{(\tll-2s)/2}\exp(\sum_{k=1}^{\infty}\frac{(\sum_{i=1}^{\tll}\beta_{j+1}^{(i)}/2)(-k)}{k}x^k)\nonumber\\
&\quad{}\times
\e(\sum_{i\in S}\frac{-\beta_{j}^{(i)}}{2}+\sum_{i\in \os}\frac{\beta_{j+1}^{(i)}-\beta_{j+2}^{(i)}}{2})\nonumber\\
&\quad{}+\zeta_{36}^{18\tll}x^{(\tll-2s)/2}\exp(\sum_{k=1}^{\infty}\frac{(\sum_{i=1}^{\tll}\beta_{j+2}^{(i)}/2)(-k)}{k}x^k)\nonumber\\
&\quad{}\times
\e(\sum_{i\in S}\frac{-\beta_{j}^{(i)}}{2}+\sum_{i\in \os}\frac{\beta_{j+2}^{(i)}-\beta_{j+1}^{(i)}}{2}).
\label{eq:dec-t}
\end{align}
Using \eqref{eq:mul1} we have
\begin{align}
\label{eq:no-va}
Y(\mathbf{u},x)\mathbf{v}
&=\sum_{j=0}^{2}(\Omega_{1j}+\Omega_{2j})
\end{align}
and hence
\begin{align}
o(\mathbf{u}\circ 
\mathbf{v})
&
=\sum_{j=0}^{2}\sum_{r=0}^{\tll/2}\binom{\tll/2}{r}Y_{W^0}\big(\Omega_{1j},\anx\big)\big|_{x^{1-r},\anx^{-\tll+r-1}}\nonumber\\
&\quad{}+\sum_{j=0}^{2}\sum_{r=0}^{\tll/2}\binom{\tll/2}{r}Y_{W^0}\big(\Omega_{2j},\anx\big)\big|_{x^{1-r},\anx^{-\tll+r-1}}.\label{eq:oee12}
\end{align} 
In the decomposition \eqref{eq:cm-dec} we have
$\Omega_{10}+
\Omega_{11}+
\Omega_{12}\in V_{L^{\oplus {\ell}}}^{\htau}((x))$
and $\Omega_{20}+
\Omega_{21}+
\Omega_{22}\in V_{L_{((c^{\tll}0^{\ell-\tll}),\zero)}}((x))$
since $(\beta_{j+1}-\beta_{j+2})/2=\beta_{j}/2+\beta_{j+1}, j=0,1,2$.

By \cite[Section 4]{TY}, we see that the top level of
$M^0$ is spanned by $\{\e(\lweight ; \epsilon)\ |\ \epsilon\in\{1,-1\}^{\ell}\}$. 
We shall compute 
$\sum_{r=0}^{\tll/2}\binom{\tll/2}{r}Y_{W^0}\big(\Omega_{1j},\anx\big)\big|_{x^{1-r},\anx^{-\tll+r-1}}
\e(\lweight ; \epsilon), j=0,1,2$.
A similar computation as \cite[(8.6.9)]{FLM} shows the following formula:
\begin{align}
&Y_{W^{0}}\big(
\Omega_{1j},\anx
\big)\e(\lweight ; \epsilon)\nonumber\\
&=\zeta_{36}^{9\tll+18s}Y_{W^{0}}
\big(x^{-\tll+2s}\exp\big(\sum_{k=1}^{\infty}\frac{(\sum_{i=1}^{\tll}\beta_{j}^{(i)}/2)(-k)}{k}x^k\big)
\e(\sum_{i\in S}\beta_{j}^{(i)}),\anx
\big)\e(\lweight ; \epsilon)\nonumber\\
&=\zeta_{36}^{9\tll+18s}x^{-\tll+2s}
\anx^{\langle\sum_{i\in S}\beta_{j}^{(i)},\beta(\lweight ; \epsilon)\rangle}(1+\frac{x}{\anx})^{\langle\sum_{i=1}^{\tll}\beta_{j}^{(i)}/2,\beta(\lweight ; \epsilon)\rangle}\nonumber\\
&\quad{}\times\exp\big(\sum_{k=1}^{\infty}
\frac{(\sum_{i\in S}\beta_{j}^{(i)})(-k)}{k}\anx^{k}\big)
\nonumber\\
&\quad{}\times
\exp\big(\sum_{k=1}^{\infty}\big(
(\sum_{i=1}^{\tll}\frac{\beta_{j}^{(i)}}{2k})(-k)(\anx+x)^{k}-(\sum_{i=1}^{\tll}
\frac{\beta_{j}^{(i)}}{2k})(-k)\anx^{k}\big)\big)\nonumber\\
&\quad{}\times\e(\sum_{i\in S}\beta_{j}^{(i)})\e(\lweight ; \epsilon).\label{eq:expand1}
\end{align}
Setting
\begin{align*}
\Psi&=\zeta_{36}^{9\tll+18s}
\exp\big(\sum_{m=1}^{\infty}\frac{(\sum_{i\in S}\beta_{j}^{(i)})(-m)}{m}\anx^{m}\big)
\nonumber\\
&\quad{}\times
\exp\big(\sum_{k=1}^{\infty}\big(
(\sum_{i=1}^{\tll}\frac{\beta_{j}^{(i)}}{2k})(-k)(\anx+x)^{k}-(\sum_{i=1}^{\tll}\frac{\beta_{j}^{(i)}}{2k})(-k)\anx^{j}\big)\big),
\end{align*}
we have
\begin{align}
&Y_{W^{0}}\big(
\Omega_{1j},\anx
\big)\e(\lweight ; \epsilon)\nonumber\\
&=
\sum_{t=0}^{\infty}\binom{\langle\sum_{i=1}^{\tll}\beta_{j}^{(i)}/2,\beta(\lweight ; \epsilon)\rangle}{t}
x^{-\tll+2s+t}
\anx^{-t+\langle\sum_{i\in S}\beta_{j}^{(i)},\beta(\lweight ; \epsilon)\rangle}\nonumber\\
&\quad{}\times\Psi \e(\sum_{i\in S}\beta_{j}^{(i)})\e(\lweight ; \epsilon).\label{eq:expand2}
\end{align}

Let $r$ be an integer with $0\leq r\leq \tll/2$.
To describe the first term of \eqref{eq:oee12},
we need to investigate the coefficient of $x^{1-r}\anx^{-\tll+r-1}$ in (\ref{eq:expand2}).
First,  we shall discuss the case that
there is a nonnegative integer $t$ such that $1-r\geq -\tll+2s+t$ and 
$-\tll+r-1\geq
-t+\langle\sum_{i\in S}\beta_{j}^{(i)},\beta(\lweight ; \epsilon)\rangle$.
Note that if no such $t$ exists, then
the coefficient of $x^{1-r}\anx^{-\tll+r-1}$ in (\ref{eq:expand2}) is equal to zero
since $\Psi\e(\sum_{i\in S}\beta_{j}^{(i)})\e(\lweight ; \epsilon)\in W^{0}[[x,\anx]]$.
Since $\langle \pm\beta_i/2,\pm\beta_j/2\rangle\in\{\pm 1,\pm 1/2\}$ for $0\leq i,j\leq 2$,
we have
$2s+\langle\sum_{i\in S}\beta_{j}^{(i)},\beta(\lweight ; \epsilon)\rangle\geq 0$
and hence
\begin{align*}
-\tll+r-1&\geq
-t+\langle\sum_{i\in S}\beta_{j}^{(i)},\beta(\lweight ; \epsilon)\rangle\\
&\geq -\tll+2s-1+r+\langle\sum_{i\in S}\beta_{j}^{(i)},\beta(\lweight ; \epsilon)\rangle\\
&\geq -\tll+r-1.
\end{align*}
This implies that $2s+\langle\sum_{i\in S}\beta_{j}^{(i)},\beta(\lweight ; \epsilon)\rangle=0$,
$t=\tll-2s-r+1$, and
the coefficient of $x^{1-r}\anx^{-\tll+r-1}$ in (\ref{eq:expand2}) is 
\begin{align}
\zeta_{36}^{9\tll+18s}\delta_{\langle\sum_{i\in S}\beta_{j}^{(i)}/2,\beta(\lweight ; \epsilon)\rangle,-s}
\binom{\langle\sum_{i=1}^{\tll}\beta_{j}^{(i)}/2,\beta(\lweight ; \epsilon)\rangle}{\tll-2s-r+1}
\e(\sum_{i\in S}\beta_{j}^{(i)})\e(\lweight ; \epsilon).\label{eq:coe}
\end{align}

Next, we shall discuss the case that $1-r<-\tll+2s+t$ or $-\tll+r-1<
-t+\langle\sum_{i\in S}\beta_{j}^{(i)},\beta(\lweight ; \epsilon)\rangle$
for all nonnegative integer $t$.
Since $\Psi\e(\sum_{i\in S}\beta_{j}^{(i)})\e(\lweight ; \epsilon)\in W^{0}[[x,\anx]]$,
the coefficient of $x^{1-r}\anx^{-\tll+r-1}$ in (\ref{eq:expand2}) is equal to $0$.
If $\tll-2s-r+1\geq 0$, then by setting $t_0=\tll-2s-r+1$, we have $1-r\geq -\tll+2s+t_0$ and hence 
\begin{align*}
-\tll+r-1&<-t_0+\langle\sum_{i\in S}\beta_{j}^{(i)},\beta(\lweight ; \epsilon)\rangle\\
&=-\tll+2s+r-1+\langle\sum_{i\in S}\beta_{j}^{(i)},\beta(\lweight ; \epsilon).
\end{align*}
Thus, in this case 
$2s+\langle\sum_{i\in S}\beta_{j}^{(i)},\beta(\lweight ; \epsilon)\rangle\neq 0$
and hence the coefficient of $x^{1-r}\anx^{-\tll+r-1}$ in (\ref{eq:expand2}) is also
given by \eqref{eq:coe}.
By \eqref{eq:coe} and $\langle\beta(\lambda^1),\beta(\lambda^2 ; \epsilon)\rangle=
\langle\beta(\lambda^1 ; \epsilon),\beta(\lambda^2)\rangle$
for $\lambda^1,\lambda^2\in\Klein^{\tll}$, we have obtained
\begin{align}
&\sum_{r=0}^{{\tll}/2}\binom{{\tll}/2}{r}Y_{W^0}(\Omega_{1j},y)|_{x^{1-r},\anx^{-\tll+r-1}}\e(\lweight ; \epsilon)\nonumber\\
&=\zeta_{36}^{9\tll+18s}\delta_{\langle\sum_{i\in S}\epsilon_i\beta_{j}^{(i)},\beta(\lweight)\rangle,-s}
\sum_{r=0}^{{\tll}/2}\binom{{\tll}/2}{r}
\binom{\langle\sum_{i=1}^{\tll}\beta_{j}^{(i)}/2,\beta(\lweight ; \epsilon)\rangle}{\tll-2s-r+1}
\e(\sum_{i\in S}\beta_{j}^{(i)})\e(\lweight ; \epsilon)\nonumber\\
&=\zeta_{36}^{9\tll+18s}
\delta_{\langle\sum_{i\in S}\epsilon_i\beta_{j}^{(i)},\beta(\lweight)\rangle,-s}
\binom{\langle\sum_{i=1}^{\tll}\epsilon_i\beta_{j}^{(i)}/2,\beta(\lweight)\rangle+{\tll}/2}{{\tll}-2s+1}
\e(\sum_{i\in S}\beta_{j}^{(i)})\e(\lweight ; \epsilon)
\label{eq:com-un}
\end{align}
for $j=0,1,2$.

We next investigate 
$\sum_{r=0}^{\tll/2}\binom{\tll/2}{r}Y_{W^0}\big(\Omega_{2j},\anx\big)\big|_{x^{1-r},\anx^{-\tll+r-1}}\e(\lweight ; \epsilon)$
for $j=0,1,2$.
We expand $\Omega_{2j}$ as
\begin{align*}
\Omega_{2j}&=\zeta_{36}^{9\tll+18s}x^{(\tll-2s)/2}(1+
(\sum_{i=1}^{\tll}\dfrac{\beta_{j+1}^{(i)}}{2})(-1)x+\cdots)\nonumber\\
&\quad{}\times
\e(\sum_{i\in S}\frac{-\beta_{j}^{(i)}}{2}+\sum_{i\in \os}\frac{\beta_{j+1}^{(i)}-\beta_{j+2}^{(i)}}{2})\nonumber\\
&\quad{}+\zeta_{36}^{18\tll}x^{(\tll-2s)/2}
(1+(\sum_{i=1}^{\tll}\dfrac{\beta_{j+2}^{(i)}}{2})(-1)x+\cdots)\nonumber\\
&\quad{}\times
\e(\sum_{i\in S}\frac{-\beta_{j}^{(i)}}{2}+\sum_{i\in \os}\frac{\beta_{j+2}^{(i)}-\beta_{j+1}^{(i)}}{2}).
\end{align*}
If $0\leq s\leq\tll/2-2$, then $(\tll-2s)/2\geq 2$ and hence
$\Omega_{2j}\in x^2W[[x]]$.
This tells us that
\begin{align}
\sum_{r=0}^{\tll/2}\binom{\tll/2}{r}Y_{W^0}\big(\Omega_{2j},\anx\big)\big|_{x^{1-r},\anx^{-\tll+r-1}}\e(\lweight ; \epsilon)
&=0.\label{eq:ell1}
\end{align}
In the case of $s=\tll/2-1, \tll/2$, 
we do not need explicit expressions of 
$\sum_{r=0}^{\tll/2}\binom{\tll/2}{r}Y_{W^0}\big(\Omega_{2j},\anx\big)
\big|_{
x^{1-r},\anx^{-\tll+r-1}}\e(\lweight ; \epsilon)$
to obtain \eqref{eq:binom-1} and \eqref{eq:binom-2}.

Let $\pr_{M^0} \colon W^0\rightarrow M^0$ be a projection. By
\eqref{eq:ouv}, \eqref{eq:no-va}, \eqref{eq:com-un}, \eqref{eq:ell1} and
\cite[Theorem 2.1.2]{Z},
in the case of $0\leq s\leq {\tll}/2-2$, we have
\begin{align}
0&=\pr_{M^0}
o(\mathbf{u}\circ \mathbf{v})\e(\lweight ; \epsilon)\nonumber\\
&=\zeta_{36}^{9\tll+18s}\sum_{j=0}^{2}
\delta_{\langle\sum_{i\in S}\epsilon_i\beta_{j}^{(i)}/2,\beta(\lweight)\rangle,-s}
\binom{\langle\sum_{i=1}^{\tll}\epsilon_{i}\beta_{j}^{(i)}/2,\beta(\lweight)\rangle+{\tll}/2}{{\tll}-2s+1}\nonumber\\
&\quad{}\times
\e(\sum_{i\in S}\beta_{j}^{(i)})\e(\lweight ; \epsilon).\label{eq:pr1}
\end{align}
In the case of $s={\tll}/2-1$, we have
\begin{align}
0&=\pr_{M^0}
o(\mathbf{u}\circ \mathbf{v})\e(\lweight ; \epsilon)\nonumber\\
&=\zeta_{36}^{9\tll+18s}\sum_{j=0}^{2}\delta_{\langle\sum_{i\in S}\epsilon_i\beta_{j}^{(i)}/2,\beta(\lweight)\rangle,-{\tll}/2+1}
\binom{\langle\sum_{i=1}^{\tll}\epsilon_{i}\beta_{j}^{(i)}/2,\beta(\lweight)\rangle+{\tll}/2}{3}
\nonumber\\
&\quad{}\times \e(\sum_{i\in S}\beta_{j}^{(i)})\e(\lweight ; \epsilon)\nonumber\\
&\quad{}+\pr_{M^0}\sum_{j=0}^{2}
\sum_{r=0}^{\tll/2}\binom{\tll/2}{r}Y_{W^0}\big(\Omega_{2j},\anx\big)
\big|_{
x^{1-r},\anx^{-\tll+r-1}}\e(\lweight ; \epsilon).
\label{eq:pr2}
\end{align}
In the case of $s={\tll}/2$, we have
\begin{align}
0&=\pr_{M^0}
o(\mathbf{u}\circ \mathbf{v})\e(\lweight ; \epsilon)\nonumber\\
&=\zeta_{36}^{9\tll+18s}\sum_{j=0}^{2}\delta_{\langle\sum_{i\in S}\epsilon_i\beta_{j}^{(i)}/2,\beta(\lweight)\rangle,-{\tll}/2}
(\langle\sum_{i=1}^{\tll}\dfrac{\epsilon_i\beta_j^{(i)}}{2},\beta(\lweight)\rangle+{\tll}/2)
\e(\sum_{i\in S}\beta_{j}^{(i)})\e(\lweight ; \epsilon)\nonumber\\
&\quad{}{}+\pr_{M^0}\sum_{j=0}^{2}\sum_{r=0}^{\tll/2}\binom{\tll/2}{r}Y_{W^0}\big(\Omega_{2j},\anx\big)
\big|_{
x^{1-r},\anx^{-\tll+r-1}}\e(\lweight ; \epsilon).\label{eq:pr3}
\end{align}

If $1\leq s\leq m/2-2$, 
then (\ref{eq:binom-1}) follows from (\ref{eq:pr1})
since $\e(\sum_{i\in S}\beta_{0}^{(i)})\e(\lweight ; \epsilon)$,
$\e(\sum_{i\in S}\beta_{1}^{(i)})\e(\lweight ; \epsilon)$,
$\e(\sum_{i\in S}\beta_{2}^{(i)})\e(\lweight ; \epsilon)$ are linearly
independent.
If ${\tll}\geq 4$, then (\ref{eq:binom-2}) follows by taking $S=\varnothing$ in
(\ref{eq:pr1}).

The map $f(\cdot,x)$ defined by $f(u,x)w=\pr_{M^0}(Y_{W^0}(u,x)w)$ for 
$u\in V_{L_{((c^{\tll}0^{\ell-\tll}),{\mathbf 0})}}$
in \eqref{eq:cm-dec} and $w\in V_{L_{(\lweight,{\mathbf 0})}}$
is an element of $I_{V_{L^{\oplus {\ell}}}^{\htau}}
\binom{V_{L_{(\lweight,{\mathbf 0})}}}{V_{L_{((c^{\tll}0^{\ell-\tll}),{\mathbf 0})}}\ V_{L_{(\lweight,{\mathbf 0})}}}$.

Suppose $\lweight\not\equiv_{\htau}(c^{\tll}0^{\ell-\tll})$.
Then,
by \eqref{eq:fusion-Ll-u-1} we have
\begin{align*}
\dim_{\C} I_{V_{L^{\oplus{\ell}}}^{\htau}} 
\binom{V_{L_{(\lweight,{\mathbf 0})}}}
{V_{L_{((c^{\tll}0^{\ell-\tll}),{\mathbf 0})}}\ V_{L_{(\lweight,{\mathbf 0})}}}&=0
\end{align*}
and hence 
in \eqref{eq:pr2} and \eqref{eq:pr3} the second terms are zero:
\begin{align*}
\pr_{M^0}\sum_{j=0}^{2}\sum_{r=0}^{\tll/2}\binom{\tll/2}{r}Y_{W^0}\big(\Omega_{2j},\anx\big)
\big|_{
x^{1-r},\anx^{-\tll+r-1}}\e(\lweight ; \epsilon)
&=0.
\end{align*}
Moreover, if $1\leq s\leq m/2$, 
then \eqref{eq:binom-1} follows from \eqref{eq:pr1}--\eqref{eq:pr3}.
Taking $S=\varnothing$ in \eqref{eq:pr1}--\eqref{eq:pr3},
we have \eqref{eq:binom-2}.
\end{proof}

\begin{rmk}\label{rmk:counter}
In the case of $\ell=\tll=2$, consider the vertex operator algebra
$V_{L_{C((c,c))\times\zero}}$. We see that
$V_{L_{C((c,c))\times\zero}}(1)=\{u\in V_{L_{C((c,c))\times\zero}}\ |\ 
\htau u=\zeta_{3}u\}$ is an irreducible
$V_{L_{C((c,c))\times\zero}}^{\htau}$-module and that
\begin{equation*}
V_{L_{C((c,c))\times\zero}}(1)\cong
V_{L_{(\zero,\zero)}}(1)\oplus V_{L_{((c,c),\zero)}}
\end{equation*}
as $V_{L^{\oplus 2}}^{\htau}$-modules.
Note that the top level of $V_{L_{((c,c),\zero)}}$ is a subspace of the top level of
$V_{L_{C((c,c))\times\zero}}(1)$
and is a subspaces of $(V_{L_{C((c,c))\times\zero}})_1$.
However, we have
\begin{align*}
&\sum_{j=0}^{2}\binom{\langle\beta^{(1)}_{j}/2+\beta^{(2)}_{j}/2,\beta((c,c))\rangle+2/2}{2+1}\\
&=
\binom{\langle\beta^{(1)}_{1}/2+\beta^{(2)}_{1}/2,\beta^{(1)}_{1}/2+\beta^{(2)}_{1}/2\rangle+1}{3}\\
&\quad{}+
\sum_{j=0,2}\binom{\langle\beta^{(1)}_{j}/2+\beta^{(2)}_{j}/2,\beta^{(1)}_{1}/2+\beta^{(2)}_{1}/2\rangle+1}{3}\\
&=\binom{2+1}{3}+2\binom{-1+1}{3}\\
&=1\neq 0.
\end{align*}
Hence formula (\ref{eq:binom-2}) does not hold in this case.
\end{rmk}

\begin{lem}\label{lem:restrict}
Assume that $\lweight\not\equiv_{\htau}(c^{\tll}0^{\ell-\tll})$. For $j=a,b,c$,
the following assertions hold.

$(1)$ If ${\tll}/2\leq |S_j(\lweight)|\leq \tll$, then
$|S_j(\lweight)|={\tll}/2$, $|S_k(\lweight)|=0$ for all $k\neq j$, and $|S_j(\lweight)|$ is an even integer.
In particular, $\langle (k^{\tll}0^{\ell-\tll}),\lweight\rangle_{\Klein}=0$ for all $k=a,b,c$.

$(2)$ If $1\leq |S_j(\lweight)|\leq {\tll}/2$, then
$\sum_{k\in\{a,b,c\},k\neq j}|S_k(\lweight)|$ is an even integer.
In particular, $\langle (j^{\tll}0^{\ell-\tll}),\lweight\rangle_{\Klein}=0$.
\end{lem}
\begin{proof}
Suppose ${\tll}/2\leq|S_j(\lweight)|\leq \tll$. 
We use the notation defined just after \eqref{eq:sstar}.
Take $S\subset
S_j(\lweight)$ such that $|S|={\tll}/2$ and set
$\epsilon=(\epsilon_i)\in\{1,-1\}^{\ell}$ by
\begin{equation}
\label{eq:epsilon-lweight}
\epsilon_i = \left\{\begin{array}{rl}
-1 &\text{if }i\in S_j(\lweight),\\
1 &\text{otherwise}.
\end{array}\right.
\end{equation}
Then,
$\langle\sum_{i\in S}\epsilon_i\beta^{(i)}(j),\beta(\lweight)\rangle=-|S|$
and by (\ref{eq:binom-1}),
\begin{align*}
0&=
\binom{\langle\beta((j^{\tll}0^{\ell-\tll});\epsilon),\beta(\lweight)\rangle+{\tll}/2}{{\tll}-2|S|+1}\\
&=
\binom{-|S_j(\lweight)|-\sum_{k\in\{a,b,c\},k\neq j}|S_k(\lweight)|/2+{\tll}/2}{{\tll}-2\cdot {\tll}/2+1}\\
&=-|S_j(\lweight)|-\sum_{k\in\{a,b,c\},k\neq j}|S_k(\lweight)|/2+{\tll}/2\\
&\leq -{\tll}/2-\sum_{k\in\{a,b,c\},k\neq j}|S_k(\lweight)|/2+{\tll}/2\\
&=-\sum_{k\in\{a,b,c\},k\neq j}|S_k(\lweight)|/2\leq 0.
\end{align*}
Thus $|S_j(\lweight)|={\tll}/2$ and $|S_k(\lweight)|=0$ for all
$k\in\{a,b,c\}, k\neq j$. By (\ref{eq:binom-2}),
\begin{align*}
0&=\sum_{k=a,b,c}
\binom{\langle\beta((k^{\tll}0^{\ell-\tll});\epsilon),\beta(\lweight)\rangle+{\tll}/2}{{\tll}+1}\\
&=\binom{-|S_j(\lweight)|+{\tll}/2}{{\tll}+1}+\sum_{k\in\{a,b,c\},k\neq j}
\binom{|S_j(\lweight)|/2+{\tll}/2}{{\tll}+1}\\
&=\binom{0}{{\tll}+1}+2\binom{\tll/4+{\tll}/2}{{\tll}+1}\\
&=2\binom{3\tll/4}{{\tll}+1}.
\end{align*}
Hence $|S_j(\lweight)|={\tll}/2$ is even. In particular,
$\langle(k^{\tll}0^{\ell-\tll}),\lweight\rangle_{\Klein}=0$ for $k=a,b,c$.
Therefore, (1) holds.

Suppose $1\leq |S_j(\lweight)|\leq {\tll}/2$. Set
$\epsilon=(\epsilon_i)\in\{-1,1\}^{\ell}$ by \eqref{eq:epsilon-lweight}.
Then
$\langle\sum_{i\in S_j(\lweight)}\epsilon_i\beta^{(i)}(j),
\beta(\lweight)\rangle =-|S_j(\lweight)|$. By (\ref{eq:binom-1}),
we have
\begin{align*}
0&=
\binom{\langle\beta((j^{\tll}0^{\ell-\tll});\epsilon),\beta(\lweight)\rangle+{\tll}/2}{{\tll}-2|S_j(\lweight)|+1}\\
&=\binom{-|S_j(\lweight)|-\sum_{k\in\{a,b,c\},k\neq j}|S_k(\lweight)|/2+{\tll}/2}{{\tll}-2|S_j(\lweight)|+1}.
\end{align*}
Since $|S_j(\lweight)|\leq {\tll}/2$, $\sum_{k\in\{a,b,c\},k\neq
j}|S_k(\lweight)|$ is an even integer. Hence
$\langle(j^{\tll}0^{\ell-\tll}),\lweight\rangle_{\Klein}=0$. This proves (2).
\end{proof}

\begin{lem}\label{lem:res-orth}
$(1)$
$\langle(c^{\tll}0^{\ell-\tll}),\lweight\rangle_{\Klein}=0$.

$(2)$
If $\tll\geq 4$, then
$|\supp_{\Klein}(\lweight)\cap\{1,\ldots,m\}|<{\tll}$.
\end{lem}
\begin{proof}
We may assume
$\supp_{\Klein}(\lweight)\cap\{1,\ldots,m\}\neq \varnothing$.
First, we shall show that $\langle(c^{\tll}0^{\ell-\tll}),\lweight\rangle_{\Klein}=0$.
If $\lweight\equiv_{\htau}(c^{\tll}0^{\ell-\tll})$, then the assertion is clear from the 
definition of $\langle\,\cdot\,,\,\cdot\,\rangle_{\Klein}$. Assume that
$\lweight\not\equiv_{\htau}(c^{\tll}0^{\ell-\tll})$. We may also assume that
$|S_c(\lweight)|=0$ and $0\leq|S_a(\lweight)|,|S_b(\lweight)|\leq
{\tll}/2$ by Lemma \ref{lem:restrict}. If $1\leq
|S_a(\lweight)|, |S_b(\lweight)|$, then
$|S_b(\lweight)|=|S_b(\lweight)|+|S_c(\lweight)|$ and
$|S_a(\lweight)|=|S_a(\lweight)|+|S_c(\lweight)|$ are even integers
by Lemma \ref{lem:restrict} (2). Hence, we have
$\langle(c^{\tll}0^{\ell-\tll}),\lweight\rangle_{\Klein}=0$ in this case.
Suppose $|S_a(\lweight)|=0$. Then $|S_b(\lweight)|>0$ since 
$\supp_{\Klein}(\lweight)\cap\{1,\ldots,m\}\neq \varnothing$.
Set $\epsilon=(\epsilon_i)\in\{-1,1\}^{\ell}$ by \eqref{eq:epsilon-lweight} with $j=b$.
Note that $\binom{-|S_b(\lweight)|+{\tll}/2}{{\tll}+1}=0$ since
$1\leq |S_b(\lweight)|\leq \tll/2$. By (\ref{eq:binom-2}),
\begin{align*}
0&
=\sum_{j=a,b,c}
\binom{\langle\beta((j^{\tll}0^{\ell-\tll});\epsilon),\beta(\lweight)\rangle+{\tll}/2}{{\tll}+1}\\
&=\binom{-|S_b(\lweight)|+{\tll}/2}{{\tll}+1}+\sum_{j=a,c}
\binom{|S_b(\lweight)|/2+{\tll}/2}{{\tll}+1}\\
&=2\binom{|S_b(\lweight)|/2+{\tll}/2}{{\tll}+1}.
\end{align*}
Hence $|S_b(\lweight)|$ is an even integer. In particular,
$\langle(c^{\tll}0^{\ell-\tll}), \lweight\rangle_{\Klein}=0$. In the case of
$|S_b(\lweight)|=0$, we can show that $\langle(c^{\tll}0^{\ell-\tll}),
\lweight\rangle_{\Klein}=0$ similarly.

Next, we shall show that if $\tll\geq 4$ then $|\supp_{\Klein}(\lweight)\cap\{1,\ldots,m\}|<{\tll}$.
Suppose by contradiction that $|\supp_{\Klein}(\lweight)\cap\{1,\ldots,m\}|
= |S_a(\lweight)|+|S_b(\lweight)|+|S_c(\lweight)|=\tll$.

\medskip\noindent\emph{Case 1.}
Suppose $|S_k(\lweight)|=\tll$ for some $k\in\{a,b,c\}$.
Setting $\epsilon=(1,\ldots,1)$ in (\ref{eq:binom-2}),
we have
\begin{align*}
0&=
\sum_{j=a,b,c}\binom{\langle\beta((j^{\tll}0^{\ell-\tll})),\beta(\lweight)\rangle+{\tll}/2}{{\tll}+1}\\
&=
\binom{\langle\beta((k^{\tll}0^{\ell-\tll})),\beta(\lweight)\rangle+{\tll}/2}{{\tll}+1}\\
&\quad{}+\sum_{j\neq k}\binom{\langle\beta((j^{\tll}0^{\ell-\tll})),\beta(\lweight)\rangle+{\tll}/2}{{\tll}+1}\\
&
=\binom{{\tll}+{\tll}/2}{{\tll}+1}+2\binom{-\tll/2+{\tll}/2}{{\tll}+1}\\
&=\binom{{\tll}+{\tll}/2}{{\tll}+1}\neq 0.
\end{align*}
This is a contradiction.

\medskip\noindent\emph{Case 2.}
Suppose $|S_k(\lweight)|<\tll$ for all $k=a,b,c$.
Note that $(c^{\tll}0^{\ell-\tll})\not\equiv_{\htau}\Delta$
in this case.
There exists $j\in\{a,b,c\}$ such that
$1\leq |S_j(\lweight)|\leq {\tll}/2$ since 
$\supp_{\Klein}(\lweight)\cap\{1,\ldots,m\}\neq \varnothing$.
Set $\epsilon=(\epsilon_i)\in\{-1,1\}^{\ell}$ by \eqref{eq:epsilon-lweight}.
Then
$\langle\sum_{i\in S_j(\lweight)}\epsilon_i\beta^{(i)}(j),\beta(\lweight)\rangle
=-|S_j(\lweight)|$. By (\ref{eq:binom-1}), we have
\begin{align*}
0&=
\binom{\langle\beta((j^{\tll}0^{\ell-\tll});\epsilon),\beta(\lweight)\rangle+{\tll}/2}{{\tll}-2|S_j(\lweight)|+1}\\
&=\binom{-|S_j(\lweight)|-\sum_{k\in\{a,b,c\},k\neq j}|S_k(\lweight)|/2+{\tll}/2}{{\tll}-2|S_j(\lweight)|+1}\\
&=\binom{-|S_j(\lweight)|-({\tll}-|S_j(\lweight)|)/2+{\tll}/2}{{\tll}-2|S_j(\lweight)|+1}\\
&=\binom{-|S_j(\lweight)|/2}{{\tll}-2|S_j(\lweight)|+1}\neq 0.
\end{align*}
This is a contradiction.
Therefore, we conclude that $\wt_{\Klein}(\lweight)<{\tll}$.
\end{proof}

\begin{prop}\label{prop:orth}
Let $\mu = (\mu_k)$ be a nonzero element of $\K^{\ell}$ such that
$\wt_{\Klein}(\mu)$ is even and $\wt_{\Klein}(\mu)\geq 4$. Let
$W=\oplus_{i=0}^{\infty}W(i)$ be an $\N$-graded weak
$V_{L_{C(\mu)\times\zero}}^{\htau}$-module. Let
$M=\oplus_{i=0}^{\infty}M(i)$ be an irreducible $V_{L^{\oplus{\ell}}}^{\htau}$-submodule 
of $W$ such that $M(0)\subset W(0)$.
Assume that $M$ is isomorphic to $V_{L_{(\lweight,{\mathbf 0})}}$
for some nonzero $\lweight=(\lweight_k)\in\Klein^{\ell}$. Then
$\langle \mu,\lweight\rangle_{\Klein}=0$ and 
$|\supp_{\Klein}(\mu)\cap\supp_{\Klein}(\Delta)|<\wt_{\Klein}(\mu)$.
\end{prop}
\begin{proof}
There exists $g\in\tgrl$ such that $g(\mu)=(c^{\tll}0^{\ell-\tll})$,
where $\tll=\wt_{\Klein}(\mu)$. Consider a vertex operator
algebra $V_{L_{{C(g(\mu))\times\zero}}}^{\htau}$ and a
$V_{L_{{C(g(\mu))\times\zero}}}^{\htau}$-module $W\circ g^{-1}$ defined by
$W\circ g^{-1}=W$ as vector spaces and
$Y_{W\circ g^{-1}}(u,x)=Y_{W}(g^{-1}u,x)$ for $u\in V_{L_{{C(g(\mu))\times\zero}}}^{\htau}$.
Note that $M\circ g^{-1}$ is a $V_{L^{\oplus{\ell}}}^{\htau}$-submodule
of $W\circ g^{-1}$ which is isomorphic
to $V_{L_{(g(\lweight),{\mathbf 0})}}$. 
Since $g$ is an automorphism of $\Klein^{\ell}$,
it is sufficient to
show that $\langle g(\mu), g(\lweight)\rangle_{\Klein}=\langle
(c^{\tll}0^{\ell-\tll}), g(\lweight)\rangle_{\Klein}=0$ and
$|\supp_{\Klein}(g(\lweight))\cap \{1,\ldots,\tll\}|<\tll$ for
$V_{L_{{g(C(\mu))\times\zero}}}^{\htau}$ and a
$V_{L_{{g(C(\mu))\times\zero}}}^{\htau}$-module $W\circ g^{-1}$.
These results hold by Lemma \ref{lem:res-orth}.

\end{proof}

\subsection{Modules of $V_{L_{C\times D}}^{\htau}$}\label{subsection:vlcd}
In this subsection we shall classify the
irreducible $V_{\LCD}^\htau$-modules and establish the rationality
of $V_{\LCD}^\htau$ for arbitrary $\htau$-invariant
self-dual $\Klein$-code $C$ with minimum weight at least
$4$ and arbitrary self-dual $\Z_3$-code $D$.

For any nonzero $\mu\in \Klein^{\ell}$ of even weight and
any self-orthogonal $\Z_3$-code $D$, we have
\begin{equation}
\label{eq:mu-D-dec}
V_{L_{C(\mu)\times D}}^{\htau} \cong V_{L_{\zero\times D}}(0)\oplus V_{L_{\mu\times D}}
\end{equation}
as $V_{L_{\zero\times D}}^{\htau}$-modules. The following lemma
will be used in Lemma \ref{lem:CD-un-st} and Proposition \ref{prop:CD-twist}.

\begin{lem}\label{lem:prod}
For any nonzero $\mu\in \Klein^{\ell}$ of even weight and
any self-orthogonal $\Z_3$-code $D$, we have
$V_{L_{\mu\times D}}\cdot V_{L_{\mu\times D}}=V_{L_{C(\mu)\times D}}^{\htau}$ in \eqref{eq:mu-D-dec}.
\end{lem}

\begin{proof}
We may assume that $\mu=(c^{\tll}0^{\ell-\tll}) ,\tll>0$ by the action of $\tgrl$ (see Proof of 
Proposition \ref{prop:orth}). Then the
assertion follows from (\ref{eq:no-va}).
\end{proof}

For the remainder of this paper, $C$ is a $\htau$-invariant
self-dual $\Klein$-code of length ${\ell}$ with minimum weight at
least 4 and $D$ is a self-dual $\Z_3$-code of the same length.
Let $C_{\equiv_{\htau}}$ be the set of all orbits of
$\htau$ in $C$. Note that
\begin{equation*}
V_{L_{C\times D}}(\varepsilon) \cong V_{L_{\zero\times
D}}(\varepsilon)\oplus \bigoplus_{{\mathbf 0}\neq\lambda\in
C_{\equiv_{\htau}}}V_{L_{\lambda\times D}}, \quad\varepsilon=0,1,2
\end{equation*}
as $V_{L_{\zero\times D}}^{\htau}$-modules by Lemma \ref{lem:lattice-st}.

By
Proposition \ref{prop:fusion-D} and Lemma \ref{lem:prod},
the same argument as in the proof of \cite[Theorem 5.4]{KMY}
shows the following lemma.

\begin{lem}\label{lem:CD-un-st}
Let $(N^1,Y^1)$ and $(N^2,Y^2)$ be irreducible
$\VLCD^{\htau}$-modules and let $\varepsilon\in\Z_3$. Suppose for
each $i=1,2$, there is a $V_{L^{\oplus {\ell}}}^{\htau}$-submodule
of $N^i$ which is isomorphic to
$V_{L_{(\zero,\zero)}}(\varepsilon)$. Then, $N^1$ and $N^2$ are
isomorphic $V_{L_{C\times D}}^{\htau}$-modules.
\end{lem}

As it was mentioned at the beginning of this section, 
we need to assume that $D$ is self-dual 
to show the following proposition.

\begin{prop}\label{prop:CD-untwist}
Let $N$ be an $\N$-graded weak $V_{L_{C\times D}}^{\htau}$-module
which has a $V_{L^{\oplus {\ell}}}^{\htau}$-submodule isomorphic
to $V_{L_{(\lambda,\gamma)}}$ for some nonzero
$\lambda\in\Klein^{\ell}$ and $\gamma\in\Z_3^{\ell}$. Then there
exists a $V_{L^{\oplus {\ell}}}^{\htau}$-submodule $M$ of $N$
which is isomorphic to $V_{L_{(\zero ,\zero )}}(\varepsilon)$ for some
$\varepsilon\in\Z_3$. Consequently, there exists a
$V_{L_{\zero\times D}}^{\htau}$-submodule of $N$ which is
isomorphic to $V_{L_{\zero\times D}}(\varepsilon)$. The
$V_{L_{C\times D}}^{\htau}$-submodule of $N$ generated by $M$ is
isomorphic to $V_{L_{C\times D}}(\varepsilon)$.
\end{prop}

\begin{proof}
Let $W^1$ be an irreducible
$V_{L^{\oplus {\ell}}}^{\htau}$-submodule of
$N$ which is isomorphic to
$V_{L_{(\lambda,\gamma)}}$ for a nonzero $\lambda\in\Klein^{\ell}$
and $\gamma\in\Z_3^{\ell}$. Since $N$ is a
$V_{L_{\zero\times D}}^{\htau}$-module,
$\gamma\in D^{\perp}=D$ by Theorem
\ref{thm:irr-D} and consequently, there exists a
$V_{L^{\oplus {\ell}}}^{\htau}$-submodule $W^2$ of
$N$ which is isomorphic to
$V_{L_{(\lambda,{\mathbf 0})}}$.

Suppose for any $\varepsilon\in\Z_3$, there is no
$V_{L^{\oplus {\ell}}}^{\htau}$-submodule of
$N$ which is isomorphic to
$V_{L_{({\mathbf 0},{\mathbf 0})}}(\varepsilon)$. Let
$N^1=\oplus_{n=0}^{\infty}N^1(n)$ be the
$V_{L_{C\times\zero}}^{\htau}$-submodule of $N$ generated by $W^2$. Note that
every irreducible $V_{L^{\oplus {\ell}}}^{\htau}$-submodule of
$N^1$ is isomorphic to $V_{L_{(\lambda^1,{\mathbf 0})}}$ for a
nonzero $\lambda^1\in \Klein^{\ell}$ by Proposition \ref{prop:fusion-D} and the assumption.
Let $M=\oplus_{n=0}^{\infty}M(n)$ be an irreducible
$V_{L^{\oplus {\ell}}}^{\htau}$-submodule of
$N^1$ such that $M(0)\subset
N^1(0)$. There exists a nonzero $\lweight\in\Klein^{\ell}$ such that
$M$ is isomorphic to $V_{L_{(\lweight,{\mathbf 0})}}$ as
$V_{L^{\oplus {\ell}}}^{\htau}$-modules.
Since $N^1$ is a $V_{L_{C(\mu)\times\zero}}^{\htau}$-module for all $\mu\in C$, we have
$\langle\mu,\lweight\rangle_{\Klein}=0$ by Proposition \ref{prop:orth} and 
hence $\lweight\in C^{\perp}=C$. By Proposition \ref{prop:orth} again,
$\wt_{\Klein}(\lweight)= |\supp_{\Klein}(\lweight)\cap\supp_{\Klein}(\lweight)|<\wt_{\Klein}(\lweight)$.
This is a contradiction.
Thus, there exists an irreducible
$V_{L^{\oplus {\ell}}}^{\htau}$-module
$M$ isomorphic to
$V_{L_{(\zero,\zero)}}(\varepsilon)$ for some
$\varepsilon\in\Z_3$. By
Proposition \ref{prop:irr-L-l} and Theorem \ref{thm:irr-D}, the
$V_{L_{\zero\times D}}^{\htau}$-submodule of $N$ generated by $M$ is
isomorphic to $V_{L_{\zero\times D}}(\varepsilon)$.

Let $N^2$ be the $V_{L_{C\times D}}^{\htau}$-submodule of $N$
generated by $M$. By Proposition \ref{prop:fusion-D},
\begin{equation*}
N^2 \cong V_{L_{\zero\times D}}(\varepsilon)\oplus
\bigoplus_{{\mathbf 0}\neq \lambda\in C_{\equiv_{\htau}}}
V_{L_{\lambda\times D}}
\end{equation*}
as $V_{L_{\zero\times D}}^{\htau}$-modules (cf. Proof of Lemma
\ref{lem:generate}). Since any nonzero
$V_{L_{C\times D}}^{\htau}$-submodule of
$N^2$ must contain $V_{L_{\zero\times D}}(\varepsilon)$
by the argument above, $N^2$ is irreducible. By
Lemma \ref{lem:CD-un-st}, $N^2$ is isomorphic to
$V_{L_{C\times D}}(\varepsilon)$ as
$V_{L_{C\times D}}^{\htau}$-modules.
\end{proof}

\begin{prop}\label{prop:CD-twist}
Let $N$ be an $\N$-graded weak $V_{L_{C\times D}}^{\htau}$-module.
Suppose $N$ has a $V_{L^{\oplus {\ell}}}^{\htau}$-submodule $M$
which is isomorphic to $\vt{L^{\oplus
{\ell}}}{\eta}{\htau^i}{\varepsilon}$ for some
$\eta\in\Z_3^{\ell}$ and $\varepsilon\in\Z_3$. Then $M$ is a
$V_{L_{C\times D}}^{\htau}$-submodule of $N$ which is isomorphic
to $\vt{L_{C\times D}}{{\mathbf 0}}{\htau^i}{\varepsilon}$.
\end{prop}

\begin{proof}
Note that the $V_{L_{\zero\times D}}^{\htau}$-submodule of $N$ generated by
$M$ is isomorphic to $\vt{L_{\zero\times D}}{\zero }{\htau^i}{\varepsilon}$
by Theorem \ref{thm:irr-D}.
Take any nonzero $\lambda\in C$ and consider a vertex operator
subalgebra $V_{L_{C(\lambda)\times D}}^{\htau}$ of
$V_{L_{C\times D}}^{\htau}$. Let $N^1$ be the
$V_{L_{C(\lambda)\times D}}^{\htau}$-submodule of
$N$ generated by $M$. Note that for
$\varepsilon_1\in\Z_3$ with $\varepsilon_1\neq \varepsilon$, the
difference of the minimal eigenvalues of $\omega_1$ in
$\vt{L_{\zero\times D}}{\zero}{\htau^i}{\varepsilon_1}$ and in
$\vt{L_{\zero\times D}}{{\mathbf 0}}{\htau^i}{\varepsilon}$ is not
an integer, where $\omega$ is the Virasoro element of
$V_{L_{C(\lambda)\times D}}^{\htau}$. By Theorem
\ref{thm:irr-D} and Proposition \ref{prop:fusion-D},
$N^1$ is a direct sum of $V_{L_{\zero\times D}}^{\htau}$-modules,
each of which is isomorphic to $\vt{L_{\zero\times D}}{{\mathbf 0}}{\htau^i}{\varepsilon}$.
We write $N^1=\oplus_{j\in {\mathcal J}}M^j$, $M^j\cong
\vt{L_{\zero\times D}}{{\mathbf 0}}{\htau^i}{\varepsilon}$.
We can take $M^{j_1}=M$ for some
$j_1\in {\mathcal J}$.
For each $j\in {\mathcal J}$, let $\varphi_j \colon M^j \rightarrow
\vt{L_{\zero\times D}}{{\mathbf 0}}{\htau^i}{\varepsilon}$ be an
isomorphism of $V_{L_{\zero\times D}}^{\htau}$-modules and let
$\pr_{j} \colon N^1\rightarrow M^{j}$ be a projection.

We want to show that $|{\mathcal J}|=1$. Suppose ${\mathcal J}$
contains at least two elements and take $j_2\in{\mathcal J},
j_2\neq j_1$. For any $j\in{\mathcal J}$,
$v\in V_{L_{\lambda\times D}}$, and $w\in M$, define
$f_{j}(v,x)w=\varphi_j(\pr_j(Y_{N}(v,x)w))$. Then,
$f_{j}\in I_{V_{L_{\zero\times D}}^{\htau}}
\binom{\vt{L_{\zero\times D}}{{\mathbf 0}}{\htau^i}{\varepsilon}} {V_{L_{\lambda\times D}}\
\vt{L_{\zero\times D}}{\zero}{\htau^i}{\varepsilon}}$.
Note that
at most one $f_{j}$ is not zero since
$\dim_{\C}I_{V_{L_{\zero\times D}}^{\htau}}
\binom{\vt{L_{\zero\times D}}{{\mathbf 0}}{\htau^i}{\varepsilon}} {V_{L_{\lambda\times D}}\
\vt{L_{\zero\times D}}{{\mathbf 0}}{\htau^i}{\varepsilon}}=1$
(cf. \cite[Proof of Lemma 5.6]{TY}). Since $N^1$ is generated by $M$,
we have $f_{j_2}\neq 0$. Consequently, ${\mathcal J}=\{j_1,j_2\}$
and $f_{j_1}=0$. Namely,
\begin{equation*}
N^1=M^{j_1}\oplus M^{j_2}
\end{equation*}
and $V_{L_{\lambda\times D}}\cdot M^{j_1}=M^{j_2}$.
For any $k=1,2$, $v\in V_{L_{\lambda\times D}}$, and $w\in
M^{j_2}$, define
$f_{2,j_k}(v,x)w=\varphi_{j_k}(\pr_{j_k}(Y_{N}(v,x)w))$.
Then, $f_{2,j_k}\in I_{V_{L_{\zero\times D}}^{\htau}}
\binom{\vt{L_{\zero\times D}}{{\mathbf 0}}{\htau^i}{\varepsilon}}
{V_{L_{\lambda\times D}}\ \vt{L_{\zero\times D}}{{\zero}}{\htau^i}{\varepsilon}}$.
By Lemma \ref{lem:prod}, we have
\begin{align*}
V_{L_{\lambda\times D}}\cdot M^{j_2}&=
V_{L_{\lambda\times D}}\cdot (V_{L_{\lambda\times D}}\cdot M^{j_1})
=
(V_{L_{\lambda\times D}}\cdot V_{L_{\lambda\times D}})\cdot M^{j_1}\\
&=
(V_{L_{\zero\times D}}^{\htau}\oplus V_{L_{\lambda\times D}})\cdot M^{j_1}
=M^{j_1}\oplus M^{j_2}.
\end{align*}
Hence $f_{2,j_1}$ and $f_{2,j_2}$ are linearly independent (cf.
\cite[Proof of Lemma 5.6]{TY}). This contradicts the fact that
$\dim_{\C}I_{V_{L_{\zero\times D}}^{\htau}}
\binom{\vt{L_{\zero\times D}}{{\mathbf 0}}{\htau^i}{\varepsilon}} {V_{L_{\lambda\times D}}\
\vt{L_{\zero\times D}}{{\mathbf 0}}{\htau^i}{\varepsilon}}=1$.

Therefore, $M$ is a $V_{L_{C\times D}}^{\htau}$-submodule of $N$. By
Theorems \ref{thm:twist-L-st} and \ref{thm:irr-D},
$ \vt{L_{C\times D}}{\zero}{\htau^i}{\varepsilon}\cong
\vt{L_{\zero\times D}}{\zero}{\htau^{i}}{\varepsilon}$ as
$V_{L_{\zero\times D}}^{\htau}$-modules. 
The same arguments as in \cite[Lemma C.2]{LYY} can show that any irreducible
$V_{L_{C\times D}}^{\htau}$-module which is isomorphic to
$\vt{L_{\zero\times D}}{\zero}{\htau^i}{\varepsilon}$ as
$V_{L_{\zero\times D}}^{\htau}$-modules must be isomorphic to
$\vt{L_{C\times D}}{\zero}{\htau^i}{\varepsilon}$. Hence the assertion holds.
\end{proof}

\begin{thm}\label{thm:main}
Let $C$ be a $\htau$-invariant self-dual $\Klein$-code of length
${\ell}$ with minimum weight at least $4$ and let $D$ be a
self-dual $\Z_3$-code of the same length. Then
$V_{L_{C\times D}}^{\htau}$ is a simple, rational,
$C_2$-cofinite, and CFT type
vertex operator algebra. There are exactly  $9$ equivalence
classes of irreducible $V_{L_{C\times D}}^{\htau}$-modules which
are represented by the following ones.

$(1)$ $V_{L_{C\times D}}(\varepsilon)$, $\varepsilon=0,1,2$.

$(2)$ $\vt{L_{C\times D}}{{\mathbf 0}}{\htau^i}{\varepsilon}$,
$i=1,2,\varepsilon=0,1,2$.
\end{thm}

\begin{proof}
The
simplicity of $V_{L_{C\times D}}^{\htau}$ is a consequence of \cite[Theorem 4.4]{DM}.
Since $V_{L_{C\times D}}^{\htau}$ is a direct sum of finitely many
irreducible $V_{L_{\zero\times D}}^{\htau}$-modules,
$V_{L_{C\times D}}^{\htau}$ is $C_2$-cofinite by \cite{Buhl}. The
classification of irreducible $V_{L_{C\times D}}^{\htau}$-modules
follows from Propositions \ref{prop:CD-untwist} and \ref{prop:CD-twist}.

We shall show that $V_{L_{C\times D}}^{\htau}$ is rational. Let
$N$ be an $\N$-graded weak $V_{L_{C\times D}}^{\htau}$-module. Let
$M$ be the sum of irreducible
$V_{L_{\zero\times D}}^{\htau}$-submodules of
$N$, each of which is isomorphic to any of
$V_{L_{\zero\times D}}(\varepsilon)$,
$\vt{L_{\zero\times D}}{{\mathbf 0}}{\htau^i}{\varepsilon}$,
$\varepsilon\in\Z_3$,
$i=1,2$. We denote by $W$ the $V_{L_{C\times D}}^{\htau}$-submodule
of $N$ generated by $M$. By Propositions
\ref{prop:CD-untwist} and \ref{prop:CD-twist}, $W$ is a
completely reducible $V_{L_{C\times D}}^{\htau}$-module. If the
$V_{L_{C\times D}}^{\htau}$-module $N/W$ is not zero, then $N/W$
has a $V_{L_{\zero\times D}}^{\htau}$-submodule
isomorphic to one of $V_{L_{\zero\times D}}(\varepsilon)$,
$\vt{L_{\zero\times D}}{{\mathbf 0}}{\htau^i}{\varepsilon}$,
$\varepsilon\in\Z_3$, $i=1,2$ by Propositions
\ref{prop:CD-untwist} and \ref{prop:CD-twist}.
This contradicts our choice
of $W$. Hence $N=W$. This implies that
$V_{L_{C\times D}}^{\htau}$ is rational.
\end{proof}

\begin{rmk}\label{rmk:Leech}
In \cite{KLY1},
it is shown that there exist a $\Klein$-code $C$ of length $12$ and
a $\Z_3$-code $D$ of the same length,
which satisfy the conditions in Theorem \ref{thm:main},
and such that $L_{C\times D}$ is isomorphic to the Leech lattice $\Leech$.
In this case $\htau$ corresponds to a unique fixed-point-free isometry of $\Leech$ of order $3$
up to conjugacy (cf. \cite{Atlas}).
Hence, as a special case of Theorem \ref{thm:main},
we obtain the classification of irreducible modules,
the rationality, and the $C_2$-cofiniteness for $V_{\Leech}^{\htau}$.
\end{rmk}

\begin{rmk}\label{rmk:e8}
For $\ell = 4$, let $C$ and $D$ be a $\K$-code and a $\Z_3$-code
with generating matrices
\begin{equation*}
\begin{pmatrix}
a & a & 0 & 0\\
b & b & 0 & 0\\
0 & 0 & a & a\\
0 & 0 & b & b
\end{pmatrix}, \qquad
\begin{pmatrix}
1 & 1 & 1 & 0\\
1 & -1 & 0 & 1
\end{pmatrix},
\end{equation*}
respectively. It is clear that $C$ is $\htau$-invariant self-dual
and $D$ is self-dual.
The lattice $\LzeroD$ is a
$\sqrt{2}$\,$($$E_8$-lattice$)$ and $\LCD$ is an $E_8$-lattice. Note
that $D$ is the $[4,2,3]$ ternary tetra code.

We can not apply Theorem \ref{thm:main} to $V_{L_{C\times D}}$
since the minimum weight of $C$ equals $2$.
\end{rmk}

\section{List of Notations}

\begin{longtable}{lp{12cm}}
$\zeta_n$&$\exp(2\pi\sqrt{-1}/n)$.\\
$\langle\,\cdot\,,\,\cdot\,\rangle$&
the ordinary inner product of the Euclidean space $\R^{\ell}$.\\
$L$ & $\sqrt{2}(A_2$-lattice).\\
$L^{\perp}$ & the dual lattice of $L$.\\
$\beta_1,\beta_2$&a $\Z$-basis of $L$ such that $\langle\beta_1,\beta_1\rangle=\langle\beta_2,\beta_2\rangle=4$
and $\langle\beta_1,\beta_2\rangle=-2$.\\
$\beta_0$&$\beta_0=-\beta_1-\beta_2$.\\
$\tilde{\beta}_1,\tilde{\beta}_2$&
the basis of $L^{\perp}$ defined by
$\tilde{\beta}_1=\beta_1/2$ and $\tilde{\beta}_2=(\beta_1-\beta_2)/6$.\\
$\otau$ & an isometry of $L$ induced by the permutation $\beta_1\mapsto \beta_2\mapsto \beta_0\mapsto \beta_1$.\\
$\tgrs$&the direct product of  ${\ell}$ copies of the group
$\langle\otau\rangle$ generated by $\otau$.\\
$\htau$&$\htau=(\otau_,\ldots,\otau)\in\tgrs$
(For simplicity of notation, we denote $(\otau,\ldots,\otau)$ by $\htau$ also).\\
$\grs$&$\{(\otau^{i_1},\ldots,\otau^{i_{\ell}-1},1)\in \tgrs\ |\ i_1,\ldots,i_{\ell-1}\in\Z\}$.\\
$\tgrl$&$\tgrs\rtimes \mathfrak{S}_{\ell}$, where
$\mathfrak{S}_{\ell}$ is the symmetric group of degree $\ell$.\\
$\Klein$ & $\Klein=\{0,a,b,c\}\cong\Z_2\times\Z_2$ is Klein's four-group.\\ 
$C$ & a code over $\Klein$.\\
$D$ & a code over $\Z_3$.\\
$\supp_{\Klein}(\lambda)$&
$\supp_{\Klein}(\lambda)=\{i\ |\ \lambda_i\neq 0\}$ where $\lambda=(\lambda_1,\ldots,\lambda_{\ell})\in \Klein^{\ell}$.\\
$\supp_{\Z_3}(\gamma)$&
$\supp_{\Z_3}(\gamma)=\{i\ |\ \gamma_i\neq 0\}$ where $\gamma=(\gamma_1,\ldots,\gamma_{\ell})\in \Z_3^{\ell}$.\\
$\wt_{\Klein}(\lambda)$&the cardinality of $\supp_{\Klein}(\lambda)$.\\
$\wt_{\Z_3}(\gamma)$&the cardinality of $\supp_{\Z_3}(\gamma)$.\\
$\langle\lambda,\mu\rangle_{\Klein}$
&
$\langle\lambda,\mu\rangle_{\Klein}=\sum_{i=1}^{\ell}\lambda_i\mu_i\in\Klein$
where 
$\lambda=(\lambda_i), \mu=(\mu_i)\in\Klein^{\ell}$.\\
$\langle\gamma,\delta\rangle_{\Z_3}$
&
$\langle\gamma,\delta\rangle_{\Z_3}=\sum_{i=1}^{\ell}\gamma_i\delta_i\in\Z_3$
where 
$\gamma=(\gamma_i), \delta=(\delta_i)\in\Z_3^{\ell}$.\\
$C(\lambda)$&the $\Klein$-code generated by $\lambda\in\Klein^{\ell}$ and $\htau(\lambda)$.\\
$D(\gamma)$&the $\Z_3$-code generated by $\gamma\in\Z_3^{\ell}$.\\
$\beta(x)$&$\beta(0)=0,\beta(a)=\beta_2/2,\beta(b)=\beta_0/2,\beta(c)=\beta_1/2$.\\
$L^{(x,i)}$&$L^{(x,i)}=\beta(x)+i(-\beta_1+\beta_2)/3+L$ where $x\in\Klein$ and $i\in\Z_3$.\\
$L_{(\lambda,\gamma)}$&
$L_{(\lambda,\gamma)}=L^{(\lambda_1,\gamma_1)}\oplus\cdots\oplus L^{(\lambda_{\ell},\gamma_{\ell})}\subset (L^{\perp})^{\oplus\ell}$
where $\lambda\in \Klein^{\ell}$ and $\gamma\in \Z_3^{\ell}$.\\
$L_{P\times Q}$&
$L_{P\times Q}=\cup_{\lambda\in P,\gamma\in Q}L_{(\lambda,\gamma)}$.\\
$T$ & a subgroup in the center of $(\widehat{L^{\perp}})^{\ell}$ generated by
$\kz^{(r)}(\kz^{(s)})^{-1}$, $1\leq r,s\leq \ell$, where $\kz^{(s)}$
denotes $\kz $ in the $s$-th entry of
$(\widehat{L^{\perp}})^{\ell}$. \\
$K_0$ & $K_0 = \{ a\times_{\htau}\htau(a)^{-1}\,|\, a \in \hLCzerotau\}$.\\
$K$ & $K = \{ a\times_{\htau}\htau(a)^{-1}\,|\, a \in \hat{L}_{C\times D,\htau}\}$.\\
$\VLCDTeta(\htau^i)$ & irreducible
$\htau^i$-twisted $V_{\LCD}$-module where 
$\eta \in D^\perp$ and $i=1,2$. \\
$\tilde{\Klein}$ & $\tilde{\Klein}=\{0,1,2,a,b,c\}$ (cf. Section 5).\\
$\supp_{\tilde{\Klein}}(\lambda)$&
$\supp_{\tilde{\Klein}}(\lambda)=\{i\ |\ \lambda_i\in\{a,b,c\}\}$ 
where $\lambda=(\lambda_1,\ldots,\lambda_{\ell})\in \tilde{\Klein}^{\ell}$.\\
$\wt_{\tilde{\Klein}}(\lambda)$&the cardinality of $\supp_{\tilde{\Klein}}(\lambda)$.\\
$X_{i,j}$ & 
$X_{i,j}=\left\{\begin{array}{ll}V_{L^{(0,j)}}(i)&\mbox{if }i=0,1,2,\\
V_{L^{(i,j)}}&\mbox{if }i=a,b,c.
\end{array}\right.$\\
$X_{\xi,\gamma}$&$X_{\xi,\gamma}=\otimes_{i=1}^{\ell}X_{\xi_i,\gamma_i}$
where $\xi=(\xi_1,\ldots,\xi_{\ell})\in \tilde{\Klein}^{\ell}$
and $\gamma=(\gamma_1,\ldots,\gamma_{\ell})\in\Z_3^{\ell}$.\\
$P(V_{L_{(\zero,\gamma)}}(\varepsilon))$
&$P(V_{L_{(\zero,\gamma)}}(\varepsilon))
=
\{\xi=(\xi_k)\in\Z_3^{\ell}\ |\ \sum_{k=1}^{\ell}\xi_k\equiv \varepsilon\pmod{3}\}$ where
$\gamma\in\Z_3^{\ell}$ and $\varepsilon\in\Z_3$.\\
$P(V_{L_{(\lambda,\gamma)}})$&
$P(V_{L_{(\lambda,\gamma)}})
=\{\xi\in\{0,1,2,c\}^{\ell}\ |\
\supp_{\tilde{\Klein}}(\xi)=\supp_{\Klein}(\lambda)\}$
where
$0\neq \lambda\in\Klein^{\ell}$ and $\gamma\in\Z_3^{\ell}$.\\
$\beta^{(i)}(j)$ & $\beta^{(i)}(j)$ denotes $\beta(j)\in L^{\perp}$ in the $i$-th entry of $(L^{\perp})^{\ell}$ where $j=a,b,c$.\\
$\beta(p;\epsilon)$&
$\beta(p;\epsilon)=\sum_{i=1}^{\ell}\epsilon_i\beta^{(i)}(p_i)$
where $p=(p_i)\in\Klein^{\ell}$ and $\epsilon=(\epsilon_i)\in\{1,-1\}^{\ell}$. \\
$\beta(p)$ & $\beta(p)=\beta(p;(1,\ldots,1))$.\\
$\e(\alpha)$ & $\e(\alpha)=e^{\alpha}$
where $\alpha\in (L^{\perp})^{\oplus\ell}$.\\
$\os$ & $\os=\{i\in\{1,\ldots,m\} |\ i\not\in S\}$
where $S$ is a subset of $\{1,\ldots,\ell\}$.\\
$S_j(\lambda)$ &$S_j(\lambda)=\{k\in \{1,\ldots,m\}\ |\ \lambda_k=j\}$ where $\lambda\in\Klein^{\ell}$ and 
$j=a,b,c$.
\end{longtable}

\end{document}